\documentclass[11pt,a4paper]{article}
\usepackage{arxiv}
\usepackage[utf8]{inputenc} 
\usepackage[T1]{fontenc}    
\usepackage{hyperref}       
\usepackage{url}            
\usepackage[shortlabels]{enumitem}
\usepackage{amsmath,amsfonts,amsthm,amssymb}      
\usepackage{float}
\usepackage{rotating}     
\usepackage{graphicx}
\usepackage{microtype}
\def\allk{\forall \ k \in V_{n-1} \setminus V_{3}}
\usepackage[classfont=bold,langfont=roman,funcfont=italic]{complexity}
\usepackage{natbib}
\usepackage{xcolor}
\newtheorem{definition}{Definition}
\newtheorem{lemma}{Lemma}
\newtheorem{problem}{Problem}
\newtheorem{example}{Example}

\newtheorem{remark}{Remark}

\newtheorem{theorem}{Theorem}
\newtheorem{corollary}{Corollary}
\newtheorem{algorithm}{Algorithm}
\newtheorem{prop}{Property}

\newcounter{casenum}
\newenvironment{caseof}{\setcounter{casenum}{1}}{\vskip.5\baselineskip}
\newcommand{\case}[2]{\vskip.5\baselineskip\par\noindent {\bfseries Case \arabic{casenum}:} #1\\#2\addtocounter{casenum}{1}}

\def\fab{f_{\alpha \beta}}
\title{On the Importance of Studying the Membership Problem for Pedigree Polytopes}
\date{}					

\author{Tiru Arthanari, External Collaborator/ Visitor, \\ Department of ISOM, University of Auckland, \\ Auckland, New Zealand  \\e-mail: t.arthanari@auckland.ac.nz}

\date{}

\renewcommand{\undertitle}{}

\hypersetup{
pdftitle={On the Importance of Studying the Membership Problem for Pedigree Polytopes}
pdfsubject={Optimization}
pdfauthor={Tiru Arthanari}
pdfkeywords={Symmetric Traveling
Salesman Problem, Pedigree Polytope, Multistage Insertion
formulation, membership problem, computational Complexity}
}
\begin{document}
\maketitle

\begin{abstract}
Given $n \geq  3 $,  a combinatorial object called a \textit{ pedigree } is defined using $3$-element subsets from  $\{1,\ldots, n\}$ obeying certain conditions.  The convex hull of pedigrees is called the pedigree polytope for $n$. Pedigrees are in $1-1$ correspondence with Hamiltonian cycles. Properties of pedigrees, pedigree polytopes, adjacency structure of the graph of the pedigree polytope and their implication on the adjacency structure of the Symmetric Traveling Salesman problem (STSP)  polytope have been studied earlier in the literature by the author. 

The question:  Given $X$, does it belong to the pedigree polytope for $n$? is called the membership problem. This article provides proof that the membership problem for pedigree polytopes can be solved efficiently. Due to the pedigree's stem property, we can check the membership problem sequentially for $ k = 4, \ldots, n$. One constructs a layered network, recursively, to check membership in the pedigree polytope. Proof of the proposed framework's validity is given. This article's significant and far-reaching contribution is that the membership problem has a strongly polynomial-time framework.
 
Since the polynomial solvability of the membership problem implies that one can solve efficiently any linear optimisation problem over the pedigree polytope. And a  specific linear optimisation over the pedigree polytope (the multistage insertion formulation) solves the STSP.  
 The consequence of this result is that we have proof of $NP = P$. A recent book by the author entitled \textit{Pedigree Polytopes} brings together published results on pedigrees and some new results, mainly in Chapters 5 and 6. The primary purpose of this article is to present the latest results from that book in a self-contained fashion so that experts can vet the same. Some of the proofs and presentation of concepts in this article are new. 

\end{abstract}


\setcounter{figure}{0}
\newpage

\tableofcontents
\newpage
\section{Introduction}\label{intro}

 Given a polyhedral convex set C and a point X, deciding whether the point is a member of the convex set is called the membership problem.
Let $n$ be an integer, $n \geq 3$. Let $V_{n}$ be a set of {\em
vertices}. Assuming, without loss of generality, that the vertices
are numbered in some fixed order, we write $V_{n} = \{ 1, \ldots ,
n \}$. Let $E_{n} = \{ (i,j) | i,j \in V_{n}, i<j \}$  be the set
of {\em  edges}. The cardinality of $E_{n}$ is denoted by $ p_{n}
= n(n-1)/2$. Let $K_{n} = (V_{n} , E_{n}) $ denote the complete
graph of $n$ vertices.
We denote the elements of $E_{n}$ by $e$ where $e = ( i , j ) $.
We also use the notation $ij$ for $ ( i , j ).$ 
Unlike the usual practice, an edge is assumed to be written with
$i < j$.
Given $ n \in Z_{+}$, let $[n]$ denote $\{1, \ldots, n\}$ and $[r, n]$ denote $ \{r, \ldots, n\}.$  For  , $n \geq 3$ consider the set of triplets/triangles  given by $\Delta^{k} = \{ \{i, j, k\} \mid 1 \leq  \  i < j < k \},  k \in [3, n].$  We define a generator of a triangle as follows: $\{1, 2, 3\}$ has no generator. For any other triangle $v = \{i, j, k\}$,  if  $j >3,$  $u = \{r, i , j\}$ or $u = \{i, s , j\},$ is called a generator of $v$, otherwise $u =\{1, 2, 3\} $ is the only generator for $v$. And  $e = \{i, j\}$ is called the \textit{common edge} for $ v.$ 
For example, $v = \{2, 4, 6\}$ has the set of generators $G(v) =\{ \{1, 2, 4\}, \{2, 3, 4\}\};$ $v = \{2, 3, 7\}$ has the set of generators $G(v) =\{ \{1, 2, 3\}\}.$ Notice that $G(v =\{i, j, k\})$ is same for all $ k \in [j+1, n].$ We can now define a pedigree.

\begin{definition}[Pedigree]\label{pedig}
    A pedigree $P$ is a sequence  of $n-2$ triangles, 
$$P =( \{1, 2, 3\}, \{i_{4}, j_{4}, 4\}, \ldots, \{i_{n}, j_{n}, n\})$$ such that 
\begin{itemize}[label=]
\item{1)} $\{i_{k}, j_{k}, k\} \in \Delta^{k}$ has a generator in $P$ for all $k \in [4, n]$  and
\item{2)} the common edges, $\{i_{4}, j_{4}\}, \ldots, \{i_{n}, j_{n}\}$, are all distinct.
\end{itemize}
\end{definition}

For example, $P = ( \{1, 2, 3\}, \{1, 3, 4\}, \{3, 4, 5\})$ is a pedigree as $\{1, 3, 4\}$ has $\{1, 2, 3\}$ as  generator and the common edge is $\{1, 3\}$ and $\{3, 4, 5\}$ has $\{1, 3, 4\}$ as  generator and the common edge is $\{3, 4\}.$ And the common edges are distinct. But $P = ( \{1, 2, 3\}, \{1, 2, 4\}, \{3, 4, 5\})$ is not a pedigree as $\{3, 4, 5\}$ has no generator. Or $P = ( \{1, 2, 3\}, \{1, 2, 4\}, \{2, 4, 5\}, \{1, 2, 6\} )$ is not a pedigree as the common edges are not distinct.

Observe that the sequence of common edges can equivalently give a pedigree $P$, that is, $(e_{4}, \ldots, e_{n} )$, where $e_{k} = \{i_{k}, j_{k}\}, k\in [4, n].$
For example, $P = ( \{1, 2, 3\}, \{1, 3, 4\}, \{3, 4, 5\}, \{1, 4, 6\})$ can be given by $
(\{1, 3\}, \{3, 4\}, \{1, 4\} ).$ We shall use this as an equivalent definition of pedigree when needed. 
Other equivalent definitions of a pedigree are given in~\cite{TSADMpaper, arthanari2023pedigree}. 
Given a pedigree $P$ for $n$, and a $k \in [4,n-1]$ let $P/k$ denote the first $k-2$ elements of $P$, that is, $( \{1, 2, 3\}, \{i_{4}, j_{4}, 4\}, \ldots, \{i_{k}, j_{k}, k\}).$ We then have,

\begin{lemma}[Stem Property]
If $P$ is a pedigree for $n$, then $P/k$  for $k \in [4, n-1]$  is a pedigree for $k$. 
\end{lemma}

Given a pedigree $P \in P_{k}$, we say $P^{\prime} \in P_{k+1}$ is an \textit{extension}  of $P$ in case there exits a $u = \{i, j, k+1\}$, such that $P^{\prime} = (P, u)$ and $u$ has a generator in $P$ and $e =\{i, j\}$ does not appear in the set of common edges of $P$. 
Assume that the triangles $\{i, j, k\}$ for a $k$ are given in a fixed order. (Say, according to the \textit{ edge label }, $l_{ij} = i+ p_{j-1}$, where $p_{j} = j (j-1)/2.$  So, $\{1, 3, 6\}, \{2, 5, 6\} $ are the second and the eighth triangles, respectively, among triangles $\{i, j, 6\}\in \Delta^{6}$.) Let $E_{k-1}  = \{ \{i, j\} | 1 \leq i < j < k\}.$ Denote by $|E_{k-1}|$ the number of elements in $E_{k-1}.$ Notice that $|\Delta^{k}| = |E_{k-1}|$. Let $B^{\binom{k}{3}}$ denote the set of $0-1$ vectors having $\binom{k}{3}$ coordinates. Let us denote the characteristic vector corresponding to a pedigree $P$ by $X = ({\bf x}_{3}, \ldots, {\bf x}_{n} )$ that is, ${\bf x}_{k} \in B^{|E_{k-1}|}, x_{k}(u) = 1 $ if $u = \{i, j, k\}$ is a triangle in $P$ and $x_{k}(u) = 0$ otherwise, for $k \in [3, n].$ The set of characteristic vectors (CV) of pedigrees for $n$ is denoted by $P_{n}.$ We denote by $conv(P_{n})$ the convex hull of  $P_{n}$, and  is called the \textit{ pedigree polytope } for $n$. 
For example, $P = ( \{1, 2, 3\}, \{1, 3, 4\}, \{3, 4, 5\}, \{1, 4, 6\})$ has the corresponding $X = (1, 0 1 0, 0 0 0 0 0 1, 0 0 0 1 0 0 0 0 0 0) \in R^{20}.$  [Note: We may say pedigree instead of characteristic vector of a pedigree when it is clear from the context that we are referring to the CV of a pedigree. Let $\tau_{k} = \binom{k}{3} -1$. We may give $X \in R^{\tau_{n}}$ leaving out $\{1, 2, 3\}$.]

Given $X \in conv(P_{k}),$ let $\Lambda_{k}(X)$ denote the set of all weight vectors $\lambda$ such that $\sum_{i \in I(\lambda)}\lambda_{i}X^{i} = X$ and $\sum_{i \in I(\lambda)} \lambda_{i} = 1$, where $I(\lambda)$ is the set of indices for which $\lambda_{i} $ is positive. And $X^{i} \in P_{k},$ are distinct pedigrees given in a fixed order. Given $X \in conv(P_{k})$, any pedigree that receives a positive weight in any convex combination of pedigrees yielding $X$ is called an \textit{active} pedigree for $X$.

\begin{lemma}
  Given $X \in conv(P_{n})$   implies $X/k \in conv(P_{k}), \forall \  k\in [4, n-1].$ 
\end{lemma} 
\textbf{Proof:}
   Let  $k < n, \textrm{ be the largest  } k \textrm{ for which } X/k \notin conv(P_{k})$.  Then $X/k$ cannot be a convex combination of pedigrees in $P_{k}$ alone. However, we have $X/k+1 \in conv(P_{k+1}),$  by our supposition. So, there exits a weight vector $\lambda,$ with $\sum_{i}\lambda_{i} = 1 \ \& \ \lambda_{i} >0 $ such that $ X/k+1 = \sum_{i}\lambda_{i}X^{i}$ where  $X^{i} \in P_{k+1}.$ Now $X^{i}/k \in P_{k}$ for all $i \in I(\lambda)$. Partition $I(\lambda)$ such that $I_{l}  = \{i \in I(\lambda)| X^{i}/k = Y^{l} \in P_{k} \textrm{ for some } l \}.$  Thus, $X/k = \sum_{l}(\sum_{i \in I_{l}}\lambda_{i})Y^{l}$. We have a contradiction. \qed
   
    This lemma shows that checking the membership of $X$ in the pedigree polytope for $n$, $conv(P_{n})$, can be done recursively by checking whether $X/k \in conv(P_{k})$ for $k \in [4, n].$ If $X/k \notin conv(P_{k})$ for some $k$, we conclude $X \notin conv(P_{n}).$ For a given $n, \  X~\in~R^{\binom{n}{3}}$, consider $X/4 =({\bf x}_{3},{\bf x}_{4}).$ If $X/4$ is not a non-negative vector, or $\sum_{u}x_{3}(u) \neq 1,  \textrm{ or } \sum_{u}x_{4}(u) \neq 1$, definitely $X/4 \notin conv(P_{4})$. Otherwise $X/4 \in conv(P_{4}).$ How do we check $X/k \in conv(P_{k})$ in general? Answering this question is the main objective of this article. 

A subset $H$ of $E_n$ is called a {\em Hamiltonian cycle} in $K_n,$ if it is the edge set of a simple cycle in $K_n$, of length $n$.
We also call such a Hamiltonian cycle a $ n-tour$ in $K_n.$ At
times we  represent  $H$ by the vector $(1 i_{2} \ldots i_{n} 1)$
where $(i_{2} \ldots i_{n})$ is a permutation of $(2 \ldots n),$
corresponding to $H$.
Let $Q_{n}$ denote the standard $STSP$ polytope, given by
$$ Q_{n} = conv(\{ X_H: X_H  \textrm{ is the characteristic vector of } H \in \mathcal{H}_{n} \})$$
where $\mathcal{H}_{n}$ denotes the set of all {\em Hamiltonian
cycles} ( or \  $n - tours$ ) in $K_{n}.$ \

In polyhedral combinatoric approaches, generally,  $Q_{n}$ is
studied while solving $STSP$ (see \cite{TSPbook}). 

Given  $H \in \mathcal{H}_{k-1}$, the operation  {\em insertion}
is defined as follows: Let $e= (i,j) \in H $. Inserting $k$ in $e$
is equivalent to replacing $e$ in $H$ by $\{(i,k), (j,k)\}$
obtaining a $k-tour.$ When we denote $H$ as a subset of $E_{k-1}$,
then inserting $k$ in $e$ gives us a $H' \in \mathcal{H}_{k}$ such
that,
$$ H' = (H \cup \{(i,k),(j,k)\}) \setminus \{e\}.$$
We write $H \overrightarrow{  \ \ e,k \ \  } H'.$

Similarly, the inverse operation, shrinking, can be defined.
Pedigree can be equivalently defined using Hamiltonian cycles explicitly as follows:
\begin{definition}{[Pedigree]}
$W = (e_{4},\ldots,e_{n}) \in E_{3} \times \ldots
\times E_{n-1}$  is called a {\em  pedigree}  if and only if there
exists a $H \in \mathcal{H}_{n}$ such that $H$ is obtained from
the $3-tour$ by the sequence of insertions, viz., $$ 3-tour
\overrightarrow{ \ \ e_{4},4 \ \ } H^{4} \ \ldots \ H^{n-1}
\overrightarrow{  \ \ e_{n},n \ \  } H .$$
\end{definition}
For example, $W = ((2,3), (1, 2), (3,4))$ is a pedigree as it yields, by insertion operations, the $6-tour$,
$H = \{(1,3), (3,6),(4,6),(2,4),(2,5), (1,5)\}$. 
But $W = ((2,3), (1, 2), (1,4))$ is not a pedigree, as the decision to insert $6$ in $(1,4)$ is not possible as the $5-tour$ that results from the insertion decisions for $4$ and  $5$ is $\{(1,3),(3,4),(2,4), (2,5), (1,5)\}$; The pedigree $W$ is referred to as the pedigree of $H$.
 
 The pedigree of $H$ can be obtained  by shrinking $H$ sequentially to the $3-tour$
 and noting the edge created at each stage. We then write the
 edges obtained in the reverse order of their occurrence to obtain the corresponding pedigree.
Motivation for the definition of pedigree and its connection
to $MI-$-formulation ~\cite{Alt} are given in ~\cite{TSADMpaper}. We denote by  $P_{MI}(n)$,  the set of solutions to the integer relaxation of the  $MI-$formulation. (This LP formulation, called $ MI-$relaxation, is discussed in Subsection~\ref{MI_Ins}.

Consider the decision problem, called Membership in Pedigree Polytope Problem ($M3P$) stated as:
\begin{problem}\label{problemain}
   Input: n > 3,  $X \in Q^{\tau(n)}$
   Question: Is $X$ in $conv(P_{n})$ ?
\end{problem}  
Our main task of this article is to examine the complexity of this problem, $M3P$. We assert that $M3P$ is polynomially solvable.

Next, we outline the structure of this paper.  

Section~\ref{preli} gives preliminaries and notation used.  Section~\ref{conslay}  describes the construction of the Layered Network ($N_k, R_k, \mu$). In Section~\ref{Fkfeasibility}  the feasibility of the problem called $F_{k}$ is shown to be necessary for membership in the pedigree polytope. Section~\ref{multi} defines and studies the solutions to a multicommodity flow Problem to check the membership in the pedigree polytope. This section proves that a necessary and sufficient condition for membership in pedigree polytope is the optimal total flow in the multicommodity flow problem is equal to the maximum possible flow, $z_{max}$. Section~\ref{compcomplex}  analyses the computational complexity of checking the necessary and sufficient conditions discussed in the previous section. Finally,  Section~\ref{conclusions} provides concluding remarks, implications and discussions on future research paths.
\section{Preliminaries \& Notations}\label{preli}
We repeat some notations and preliminaries from \cite{TSADMpaper}
for convenience.  Let $ R$ denote the set of reals. Similarly, $Q,
\ Z, \ N$ denote the rationals, integers, and natural numbers
respectively, and $B$ stands for the binary set of \{0, 1\}. Let $
R_{+}$ denote the set of nonnegative reals. Similarly, the
subscript ${}_{+}$ is understood with rationals. Let $R^{d}$
denote the set of $d-$tuples of reals. Similarly the superscript
${}^{d}$ is understood with rationals, etcetera. Let  $R^{m \times n}$
denote the set of $m \times n$ real matrices.

For a subset $F \subset E_{n} $ we write the {\em characteristic}
vector of $F$ by $x_{F} \in R^{p_{n}}$ where $ x_{F}(e) = 1 $ if  $e \in F$, and $0$ otherwise.

We assume that the edges in $E_n$ are ordered in increasing order of the edge labels.
For a subset $S \subset V_n $ we write

$E(S) =  \{ij| ij \in E , i,j \in S \}.
$. Given $u \in R^{p_{n}}, \ \  F \subset E_{n}$ ,we define,
$ u(F) \ = \sum_{e \in F} u(e) .$
For any $S \subset V_n,$ let $\delta(S)$ denote
the set of edges in $E_n$ with one end in $S$ and the other in
$S^{c} = V_n \setminus S.$ For $S = \{ i \}$ , we write $\delta (
\{ i \} ) = \delta (i) $.

 In addition to the notations and preliminaries introduced, we require a few definitions and concepts concerning bipartite flow problems.
For details on graph-related terms, see any standard text on graph theory such as \cite{BandM}.

\subsection{Rigid, Dummy arcs in a capacited Transportation  Problem}
Consider a balanced transportation problem, in which some arcs
called the \emph{forbidden }arcs are not available for
transportation. We call the problem of finding whether a feasible
flow exists in such an incomplete bipartite network,  a Forbidden
Arcs Transportation ($FAT$) problem \cite{KGMurty3}. This could be
viewed as a capacited transportation problem, as well. In general
a $FAT$ problem is given by $O = \{ O_{\alpha}, \alpha = 1,
\ldots, n_1 \}$, the set of origins, with availability
$a_{\alpha} \geq 0$ at $O_{\alpha}$, $D = \{ D_{\beta}, \beta = 1,
\ldots, n_2 \}$, the set of destinations with requirement
$b_{\beta} \geq 0$ at $D_{\beta}$ and $\mathcal{A} =\{(O_{\alpha},
D_{\beta})| \textrm{ arc } (O_{\alpha}, D_{\beta}) \textrm{ is
\emph{ not forbidden} } \}$, the set of arcs. We may also use
$(\alpha,\beta)$ to denote an arc. Since the problem is balanced, we have, $\sum_{\alpha}a_{\alpha} = \sum_{\beta}b_{\beta}.$
When any arc $(\alpha,\beta)$ has the capacity restriction, it is specified by $f_{\alpha \beta} \leq c_{\alpha \beta}$, as an additional restriction. This problem is stated as Problem~\ref{FAT} and has many appearances in this paper. 
\begin{problem}[$FAT$- Problem]\label{FAT}
Find $f$ satisfying
\begin{eqnarray}
\sum_{O_{\alpha} \in O}f_{\alpha \beta}& = & b(\beta), \textrm{ for all } D_{\beta} \in D \\
\sum_{D_{\beta} \in D}f_{\alpha \beta}& = & a(\alpha), \textrm{ for all } O_{\alpha} \in O \\
f_{\alpha \beta} & \leq & c_{\alpha \beta}, \textrm{ for all arc } (\alpha, \beta)  \textrm{ not forbidden and has a capacity restriction }\\
f_{\alpha, \beta} & \geq & 0, \textrm{ for all arc} (\alpha, \beta)  \textrm{ not forbidden }
\end{eqnarray}
\end{problem}
We state Lemma~\ref{theorem:lemma3}  on such
a flow feasibility problem arising with respect to non-empty partitions of a finite set.

\begin{lemma}\label{theorem:lemma3}
Suppose $\mathcal{D} \neq \emptyset$ is a finite set and  $
g:\mathcal{D} \rightarrow \ Q_{+}$, is a nonnegative function such that, $g(\emptyset)=0$, for any subset $S$ of $\mathcal{D}$, $g(S) = \sum_{d \in S}g(d)$.
Let $\mathcal{D}^{1}=\{D_{\alpha}^{1}, \alpha = 1, \ldots ,
n_{1}\},\mathcal{D}^{2}= \{D_{\beta}^{2}, \beta = 1, \ldots ,
n_{2}\}$ be two non empty partitions of $\mathcal{D}.$ ( That is,
$\bigcup_{\alpha =1}^{n_1}{D}^1_{\alpha} = \mathcal{D}$
and $D^1_{s}\bigcap D^1_{r} = \emptyset, r\neq
s.$ Similarly $\mathcal{D}^2$ is understood.)
 Consider the $FAT$ problem defined as follows:

Let the origins correspond to $D_{\alpha}^{1}$, with availability
 $a_{\alpha}  = g(D_{\alpha}^{1}),  \alpha = 1,\ldots, n_1$ and the
destinations correspond to $D_{\beta}^{2}$, with requirement
 $b_{\beta}  = g(D_{\beta}^{2}), \beta = 1,\ldots, n_2$. Let the set of arcs be given by
$$
 \mathcal{A} =\{(\alpha,\beta)| D_{\alpha}^{1} \cap D_{\beta}^{2} \neq \emptyset
 \}.
$$
Then \mbox{$\fab = g(D^{1}_{\alpha} \cap D^{2}_{\beta}) \geq 0$} is a feasible solution for the $FAT$ problem considered.
\end{lemma}
\begin{proof}
Since $\mathcal{D}^{i}, i=1,2,$ are partitions, the problem is balanced. It is easy to see that $\sum_{\alpha}\fab =g(D^{2}_{\beta}).$ That is, demand at every  
destination is met, without violating any availability restrictions. So, $f$ as defined is a feasible solution to the $FAT$ problem.
\end{proof}
However, if there are capacity restrictions on the arcs, then Lemma~\ref{theorem:lemma3} need not be true always.
Several other $FAT$ problems are defined and studied in the later sections of this paper. $FAT$ problems can be solved using any efficient bipartite maximal flow algorithm (see~\cite{KandV}). If the maximal flow is equal to the maximum possible flow, namely $a(O)$, we have a feasible solution to the problem.
\begin{definition}{[Rigid Arcs]}
Given a  $FAT$ problem with a feasible solution $f$ we say
$(\alpha , \beta) \in \mathcal{A}$ is a {\em rigid } arc in case
$f_{ \alpha , \beta}$ is same in all feasible solutions to the
problem. Rigid arcs have {\em frozen flow}.
\end{definition}

\begin{definition}{[Dummy Arc]}
A rigid arc with zero frozen flow is called a {\em dummy} arc.
\end{definition}
The set of rigid arcs in a $FAT$ problem is denoted by
$\mathcal{R}$. Identifying $\mathcal{R}$  is the {\em frozen flow
finding} problem ($FFF$ problem). Interest in this arises in
various contexts.(see \cite{Networks}) Application in statistical
data security is discussed by Guesfield \cite{Gus}. The problem of
protecting sensitive data in a two-way table, when the
insensitive data and the marginal are made public, is studied
there. A sensitive cell is unprotected if its exact value can be
identified by an adversary. This corresponds to finding rigid
arcs and their frozen flows.
\begin{example}\label{RigidEx}
[a] Consider the $FAT$ problem given by $O = \{1, 2, 3\}$, with $a(1) = 0.3, \ a(2) = 0.3, \ a(3) = 0.4,$ $D = \{4, 5, 6, 7\}$, with $b(4) = 0.4,\ b(5) = 0.2,\ b(6) = 0.1, \ b(7) = 0.3.$ The set of forbidden arcs, $F = \{ (1,6), (1,7), (3,4),(3,6)\}.$ A feasible solution to this $FAT$ problem is shown in figure~\ref{RigidFig}
\begin{figure}
    \centering
    \includegraphics[width =0.7\textwidth]{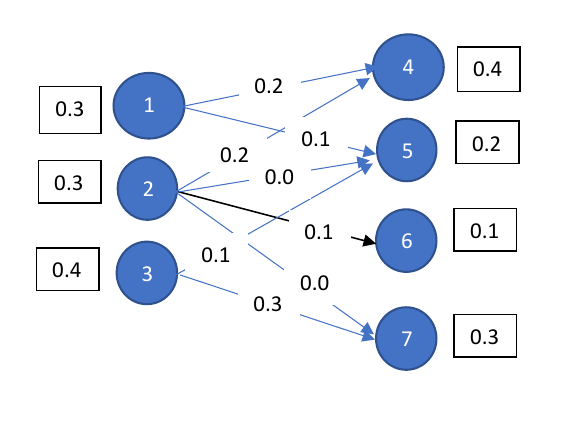}
    \caption{ $FAT$ problem for Example~\ref{RigidEx}~[a]\\ {\bf Legend:} Flows shown along the arcs, $a(i)$ and $b(j)$ are shown in boxes, for origin $i$ and destination $j$.}
    \label{RigidFig}
\end{figure}
Notice that, since $(2,6)$ is the only arc entering $6$, it is a rigid arc with frozen flow $f_{2,6} = 0.1.$ Every other arc flow can be altered without losing feasibility, and so they are not all rigid. Also, there are no dummy arcs as the only rigid arc has a positive frozen flow.
[b] Consider the $FAT$ problem discussed in [a] with the capacity restriction that $c_{1,4} = 2.$ That is $f_{1,4} \leq c_1{1,4}.$ Notice that $f$ given in [a] is feasible for this problem as well. However, the set of rigid arcs will be different. As $f_{1,4} = c_{1,4},$ and in no feasible solution $f_{1,4}$ can be less than $2$ as $b(4) = 4$ and $f_{2,4}$ can not be greater than $2$. This is so because $f_{2,6} =1$ in all feasible solutions to meet the demand at $6$. Easy to check that all arcs are rigid in this case.  Also, $(2,5)$ and $(2,7)$ are dummy arcs.
And $f$ is the unique feasible solution for this problem. (See Figure~\ref{RigidFig_b}, where rigid and dummy arcs are shown with black and red arrows, respectively.) 
\end{example}
\begin{figure}[htb]
    \centering
    \includegraphics[width = 0.7\textwidth]{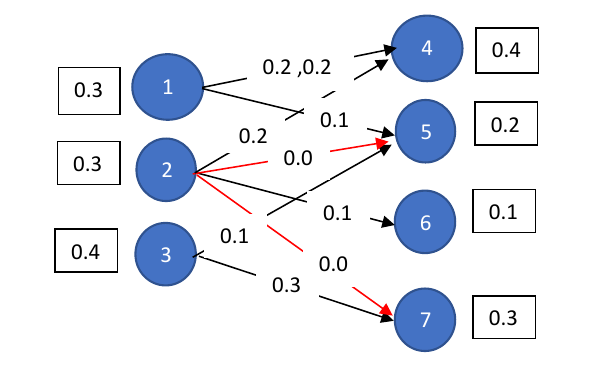}
    \caption{ $FAT$ problem for Example~\ref{RigidEx}~[b]\\ {\bf Legend:} Capacity where present, flow are shown along any arc; $a(i)$ and $b(j)$ are shown in boxes, for origin $i$ and destination $j$.}
    \label{RigidFig_b}
\end{figure}
Even though this problem can be posed as a linear programming
problem we provide the graph algorithm developed in \cite{Gus}, in Appendix~\ref{Appendixfrozen}, with required definitions. We call this algorithm as {\em frozen flow finding} or the {\em FFF} algorithm, and is a linear time algorithm, $O(|G_{f}|)$, where $G_{f}$ is a bipartite graph defined on the node set of the $FAT$ problem. 

\begin{definition}{[Layered Network]}
A network $N = (\mathcal{V}, \mathcal{A})$ is called a {\em
layered network} if the node set of $N$ can be partitioned into
$l$ sets $V_{[1]}, \ldots , V_{[l]}$ such that if $(u,v) \in
\mathcal{A}$ then $u \in V_{[r]}, v \in V_{[r+1]}$ for some $r =
1, \ldots , l-1.$ We say $N$ has $l$ layers. Nodes in $V_{[1]}$,
$V_{[l]}$ are called sources and sinks respectively.
\end{definition}
In the cases we are interested we have {\em capacities both on
arcs and nodes}. Any flow in the layered network should satisfy
apart from nonnegative and flow conservation,  the capacity
restrictions on the nodes. Of course, any such problem can be
recast as a flow problem with capacities only on arcs. (see
\cite{FF})

Given $n >3$, let $\Delta = \{\{i, j, k\}, 1 \leq i < j< k \leq n \}$. We have $\binom{n}{3}$ triangles in $\Delta$. The edge $e = \{i, j\}$ is also denoted by $(i, j)$. We may use equivalently $(i, j, k)$ or $(e, k)$ or $ [k: (i,j)]$ for a triangle. Let $\tau_{n} = \binom{k}{3} -1$, denote the number of triangles in $\Delta$ leaving out $\{ 1,2, 3\}$. Notice that $\tau_{n} = \sum_{j \in [4, n]}p_{j-1}.$ We have from~\cite{TSADMpaper, TSACompoly}, 1) the dimension of the pedigree polytope is $\tau_{n} - (n-3),$  and 2)  the polytope is combinatorial in the sense of Naddef \& Pulleyblank~\cite{naddef1981hamiltonicity}, that is, every two non-adjacent pedigrees have their midpoint equal to the midpoint of two other pedigrees. We require a few definitions and concepts concerning bipartite graphs and maximal flow in network problems. For details on graph-related terms, see any standard text on graph theory such as \cite{BandM}. Given a feasible flow $f$ to a bipartite flow feasibility problem, the arcs for which the flow remains constant in all feasible solutions to the problem are called {\emph rigid arcs}. If that arc flow is zero, such arcs are called dummy arcs. Subsection~\ref{Appendixfrozen}(Appendix A) provides definitions and concepts needed for the frozen flow finding (FFF) algorithm given by Gusfield~\cite {Gus} to find rigid and dummy arcs in a bipartite flow feasibility problem.

\section{Polytope containing the Pedigree Polytope}\label{mi_formulation}
This subsection introduces a combinatorial linear programming problem that is used in some of the proofs given in later sections, and also brings out the $1-1$ correspondence between pedigrees and Hamiltonian cycles when we have integer solutions to the problem. 

We are interested in a special class of linear programming (LP) problems, called \emph{Combinatorial LP}, introduced by~\cite{StrPoly}. Recall the definition of an LP problem.
Given $A \in R^{m \times n}, b  \in R^{m},c \in R^{n}$, a linear programming problem in standard form is stated as:
\begin{align*}
\min_x \quad c^T x \\
Ax & = b \\
x_j &\geq 0 \quad \forall j \in \{1, \dots, n\}.
\end{align*}

\begin{definition}[Combinatorial LP]
A class of linear programming problems is called  \emph{combinatorial} if the size of the entries of the matrix $A$ is polynomially bounded in the dimension of the problem. Where the dimension of the input is the number of data items in the input (that is, each number is considered to add one to the dimension of the input). The  size  of  a  rational  number  $p/q$ is the length of its binary description (i.e., 
$ size (p/q) = \lceil log_{2}(p + 1) \rceil + \lceil log_{2}(q+ 1) \rceil $ where $\lceil x \rceil$  denotes the smallest integer not less than $x$). The size of a rational vector or matrix is the sum of the sizes of its entries.
 \citep{StrPoly}
\end{definition}

\cite{StrPoly} was the first to provide a strongly polynomial algorithm that solves a combinatorial LP. That is, the number of arithmetic steps used by Tardos's algorithm depends only on the dimension of the matrix $A$, but is independent of both $b$ and $c$. 
Tardos' approach is to identify the coordinates equal to zero at optimality by solving several auxiliary dual problems via a basic algorithm for LP, such as an ellipsoid or interior-point method. References to further research in this direction, triggered by Tardos's celebrated work, see 
\citep{mizuno2018enhanced, dadush2020revisiting}.
Tardos's breakthrough theoretical result is useful in Section~\ref{compcomplex} in establishing results on the complexity of the approach outlined in the paper for the main problem.

Researchers have extensively studied the properties of the convex hull of tours, called the STSP polytope, $Q_{n}$. Dantzig, Fulkerson and Johnson’s formulation is the standard formulation of the STSP problem~\cite{DFJ}. Several classes of inequalities that define facets of $Q_{n}$ have been discovered, after the subtour elimination inequalities were introduced in the standard formulation. Branch and bound, Branch and cut and partial/intelligent search approaches use such inequalities to find an optimal tour. A vast literature covers the symmetric travelling salesman problem (STSP), different approaches to solving the problem and studying its complexity.~\cite{DFJ,TSPbook, GLS,GandJ}. It is one of the difficult combinatorial optimisation problems~\cite {Karp}.
\subsection{Pedigree Optimization Problem}
Let $P$ be a pedigree for $n.$ 
Let $x_{i j k}$ correspond to the decision to include a triangle $u = \{i, j , k\} \in \Delta^{k}$, in the pedigree, $P$. 
That is,  $x_{i j k} = 1$ if $\{i, j, k\}$ is included in $P$, $0$ otherwise. We may use $x_{u}$ as well for $x_{i j k}$ when we are not particularly interested in $i, j \ \& \ k.$
Since we select one triangle for each $k$, we have, $\sum_{\{i, j, k\} \in \Delta^{k}}x_{ijk} =1, \ \    k \in [3, n].$ 
The common edges need to be distinct for $k \in [4, n] $, which requires, an edge $\{i, j\}$ can be a common edge at most once. That is,  $\sum_{k \in [4,  n]}x_{ijk} \leq 1,$ for  $ \{i, j\} \in E_{3}$, and  we require, $\sum_{k \in [j+1,  n]}x_{ijk} \leq 1, \textrm{ for }   j >3, \{i, j\} \in E_{n-1}$. Since we require a generator for each $v =\{i,j,k\} \in P, k >3$, $ \sum_{u \in \Delta^{j} \  \& \ i \in u }x_u \geq 1, \textrm{ if } x_{i,j,k} = 1 , j \in [4, n-1], k\geq j+1.$ Or together we require, $ -\sum_{u \in \Delta^{j} \  \& \ i \in u }x_u +\sum_{k \in [ j+1, n]}x_{ijk}   \leq 0 \ \textrm{ for } \{i, j\} \in E_{n-1} \setminus E_{3}.$

Equivalently, we can check that, given a cost function $\mathcal{C}:\Delta\rightarrow R$, we can define the \textbf{Pedigree Optimisation} problem as a $0-1$ integer programming problem~\ref{pop}.
\begin{problem}[Pedigree Optimization]\label{pop}
$$minimise~{\displaystyle \sum_{u \in \Delta} \mathcal{C}_{u} x_{u}}$$
subject to
\begin{eqnarray}
\sum_{u \in \Delta^{k}}x_{u} = 1, & & k \in [3, n],\label{eq:6} \\
\sum_{k \in [4, n ]}x_{ijk} \le 1, & & (i, j) \in  E_{3}, \label{eq:7} \\
 -\sum_{u \in \Delta^{j} \  \& \ i \in u }x_u  + \sum_{k \in [j+1, n]}x_{ijk} \le 0, & &  (i, j) \in E_{n-1} \setminus E_{3}, \label{eq:8} \\
x_{u}  = 0 \textrm{ or } 1,  & &    k \in [4, n],  u \in \Delta^{k}.  \label{eq:9} 
\end{eqnarray}
\end{problem}

\subsection{Multistage Insertion and Related Results}\label{MI_Ins}

Here we will discuss the multistage insertion formulation ($ MI$-formulation)~\cite{TSA1} for its connection to the pedigree optimisation problem. (Other formulations are  discussed in~\cite{  }Insertion heuristics for STSP problem starts with the $3- tour = (1 2 3 1)$ and considers insertion of $4$, in one of the available edges, obtaining a $4-tour$, for instance, $(1 4 2 3 1)$ is obtained by inserting $4$ in the edge $(1, 2).$ Proceeding this way we obtain an $n-tour$ by inserting $n$ in an edge available in the $(n-1)-tour$ that results from the earlier insertions. These insertion decisions correspond to the selection of triangles in the pedigree optimisation problem. Selection of triangle  $x_{ijk}$ in a pedigree $P$ corresponds to the insertion decision of $k$ in $e =(i, j)$ for $k >3$. 
Let $c_{ij}$ be the cost of travelling from city $i$ to city $j$. We assume $c_{ij} = c_{ji}$  (symmetric). Let $\mathcal{C}_{123} = c_{12}+c_{23}+c_{13}$, the cost of the $3-tour$.  Let $\mathcal{C}_{ijk} = c_{ik} + c_{jk} - c_{ij}.$ Here  $\mathcal{C}_{ijk}$ gives the incremental cost of inserting $k$ between $(i, j).$ Pedigree optimisation problem, Problem~\ref{pop}, with these objective coefficients, is the $MI$- formulation of the STSP. 

Now, add to the constraints of Problem~\ref{pop} the non-binding constraints for $(i, n)$, for $i\in [1, n-1],$
\begin{equation}
    -\sum_{r = 1}^{i-1}x_{rin} - \sum_{s = i+1}^{n-1}x_{isn} \leq 0,   i = 1,\ldots ,n-1. \label{eq:10}
\end{equation}

Observe that for each $(i,j) \in E_{n}$, we have an inequality constraint (equations~\ref{eq:6},~\ref{eq:7}, and ~\ref{eq:10}) and there are $ p_{n}$ of them. 
Let the slack variable vector for these inequality constraints be given by ${\bf u}  = (u_{ij}).$ From Problem~\ref{pop}, replacing the $0 -1$ restriction on variables by $x_{ijk} \geq 0$, we obtain the $ MI$-relaxation problem, denoted by $MIR(n)$, where $n$ refers to the number of cities. The set of feasible solutions to $MIR(n)$ is denoted by the polytope $P_{MI}(n).$ Similarly we understand  $MIR(k)$, and $P_{MI}(k), k \in [4, n].$ 

\begin{lemma}\label{oneone}
Given $n \geq 4, X,$ an integer solution to $MIR(n)$, the slack variable vector, ${\bf u} \in B^{p_{n}}$ gives the edge-tour incident vector of the corresponding $n$-tour.\footnote{
Lemma~\ref{oneone} is proved in~\cite{Alt} as Lemma 4.2.}
\end{lemma}
\textbf{Proof:} For any integer solution to $MIR(n)$, [1] the slack variable for  inequality~\ref{eq:7}  is either $0$ or $1$ depending on $\sum_{k \in [4, n]} x_{ijk}$ is $1$ or $0$ respectively;  [2] the left hand side of inequality~\ref{eq:8} is $0$ if $\sum_{u \in \Delta^{j} \  \& \ i \in u }x_u = 1$ and  $\sum_{k \in [j+1, n]}x_{ijk} = 1$; and [3] the left hand side of inequality~\ref{eq:8} is $0$ if $\sum_{u \in \Delta^{j} \  \& \ i \in u }x_u = 0$. Let the integer solution  $X$ be given by $x_{i_{k}j_{k}k} =1, k \in [4, n]$. Initially, we have three edges $(i, j) \in E_3$ from the $3-tour$. Selecting $(i_{4}, j_{4})$ to insert $4$ leaves us with two of the remaining edges from $E_{3}$. In addition, two new edges $(i_{4}, 4), (j_{4}, 4)$ are created by the insertion decision resulting in a $4-tour$. And so on, with the edges available after the insertion of $n$ in $(i_{n}, j_{n})$ we have an $n-tour$. Notice that [1] whenever an available edge $(i, j)$ is used for insertion, or  [2] an edge $(i, j)$ is not created by the insertion decisions,  the corresponding inequality holds as an equality. Therefore, the inequalities that correspond to available edges, eventually, have $u_{ij} = 1$. Or ${\bf u}$ gives the edge-tour incident vector of the $n-tour$ obtained from the insertion decisions given by $X$. \qed

Given a tour we can obtain the corresponding pedigree by shrinking a $n-tour$ to an $(n-1)-tour$ by noting the edge  $(i_{n}, j_{n})$ where $n$ is inserted or choosing the corresponding triangle, and so on until we reach the $3-tour$. We have  the pedigree given by the common edges $((i_{4}, j_{4}), \ldots, (i_{n}, j_{n}))$. Thus, we have one-to-one correspondence between pedigrees and tours. All $X$ in  $P_n$ are extreme points of $P_{MI}(n)$ polytope. But
$P_{MI}(n)$ has fractional extreme points as well. Thus,  $conv(P_n) \subset P_{MI}(n).$ \citep{TSADMpaper}. 
 Therefore, it is sufficient to consider $X \in P_{MI}(n)$ while checking for membership in $conv(P_{n})$. 
\begin{theorem}
    Given $X \in P_{MI}(n)$ and $X/k \in conv(P_{k})$, if
\mbox{$x_{k+1}(e^{\prime}) = 1$} for some $e^{\prime}$, then  $X/k+1 \in conv(P_{k+1}).$
\end{theorem}\label{packability_cor}
{\bf Proof:} Since $X/k$ is in $conv(P_{k})$, we have a $\lambda \in \Lambda_{k}(X)$.  Choose any such $\lambda$. For any $r \in I(\lambda), X^{r}$ is an integer solution to $MIR(k)$, and the edge-tour incident vector of the corresponding $k-tour, T^{r},$ is given by the slack variable vector, say  $U^{r} \in B^{p_{k}}$.  $X/k$ is feasible for $MIR(k)$. Let the slack variable vector be denoted by ${\bf u}$. Since  $X/k= \sum_{r \in I(\lambda)}\lambda_{r}X^{r}$,  we can write $ {\bf u} =\sum_{r \in I(\lambda)}\lambda_{r}U^{r}.$ 
As $X/k+1 \in P_{MI}(k+1)$, $X/k+1$ is feasible for $MIR(k+1)$  which implies $x_{k+1} \leq {\bf u.}$  Or $x_{k+1}(e^{\prime}) = 1 \implies u_{e^{\prime}} = 1.$ This implies $e^{\prime} \in T^{r}, \forall \ r \in I(\lambda).$ 
Let $Y^{r} = (X^{r}, ind(\{e^{\prime}, k+1\}) ),$ where $ind(v)$ is the indicator of $v \in \Delta^{k+1}$. Therefore, $Y^{r}$ is a pedigree in $P_{k+1}$, for  $r \in I(\lambda)$, and $X/k+1 = \sum_{r \in I(\lambda)}\lambda_{r}Y^{r}.$ \qed

Theorem~\ref{packability_cor} is useful in proving a sufficient condition for membership in $conv(P_{n})$, in Section~\ref{multi}.

\begin{definition}{[$P_{MI}(n)$ Polytope]}
The polytope corresponding  to $MI$-relaxation is called the
$P_{MI}(n)$ polytope, where $n$ refers to the number of cities.
\end{definition}
All $X$ in  $P_n$ are extreme points of $P_{MI}(n)$ polytope. But
$P_{MI}(n)$ has fractional extreme points as well. Thus the pedigree polytope $conv(P_n)$ is contained in $P_{MI}(n)$ polytope \citep{TSADMpaper}. 

\begin{definition} \label{AnEn}
In general, let $E_{[n]}$ denote the matrix corresponding to
equation~(\ref{eq:6}); let $A_{[n]}$ denote the matrix
corresponding to the inequalities~(\ref{eq:6}, \ref{eq:7} \&
\ref{eq:9}). Let ${\bf 1'}_{r}$ denote the row vector of $r \  1's.$ 
Let $I_{r}$ denote the identity matrix of size $r \times r$.
\[
E_{[n]} = \left(
\begin {array}{cccc}
{\bf 1}_{p_{3}} & \ldots & 0 &  0 \\
\vdots & \vdots & \vdots & \vdots \\
0 &   \ldots & {\bf 1}_{p_{n-2}} & 0\\
 0  & 0 & \ldots & {\bf 1}_{p_{n-1}} \\
\end{array}
\right)
=\left(
\begin {array}{cc}
E_{[n-1]}  &  0 \\
 {\bf 0}  & {\bf 1}_{p_{n-1}} 
\end{array}
\right).
\]
To derive a recursive expression for $A_{[n]}$ we first define
\[
A^{(n)}= \left(
\begin {array}{c}
 I_{p_{n-1}}   \\
 -  M_{n-1} \\
\end{array}
\right)
\]
where $M_{i}$  is the $ i \times p_{i}$  {\em node-edge incidence} matrix of $K_{i}$, the complete graph on $V_{i}$, and the edges $(i,j)$ are ordered in the ascending order of the edge labels, $l_{ij}$.

Then
\[
A_{[n]}= \left(
\begin {array}{ccccccc}
 A^{(4)}  & | & A^{(5)} &  | &        & | &  A^{(n)} \\
          & | &         &  | & \ddots & | &  \\
\mathbf{0} & | & \mathbf{0}       &  | &        & | &  \\
\end{array}
\right) =
 \left( \begin{array}{ccc}
A_{[n-1]}& | & A^{(n)}\\
\mathbf{0}& | & \\
\end{array}\right).
\]
Observe that  $A^{(n)}$ is the submatrix of $A_{[n]}$
corresponding to ${\bf x}_{n}.$ The number of rows of $0's$ is
decreasing from left to right.
\end{definition}

We can state Problem~$MIR(l)$ in matrix notation for a given $l \geq 4$ as Problem~\ref{probMat}.  Let $U^{(3)}
= \left( \begin{array}{c}
     {\bf 1}_{3}  \\
     \bf{0}
  \end{array}
  \right),$ where the  $ \bf{0}$ is a column vector of $p_{l}-3$ zeros. 
     
     \begin{problem}[Problem $MI(l)$]\label{probMat}
     Find $X \in R^{\tau_{l}}, U^{(l-3)} \in R^{p_{l}}$ solution to
     $$minimise~~\mathcal{C}X$$
     subject to
    \begin{eqnarray*}
    E_{[l]}X & = & {\bf 1}_{l-3} \\
    A_{[l]}X + U^{(l)}& = & U^{(3)} \label{ineqmat} \\
  X & \geq & {\bf 0}\\
  U^{(l)} & \geq & {\bf 0}
    \end{eqnarray*}
     \end{problem}

Expanding the l.h.s of the equations in problem~\ref{probMat} for $l$ using the recursive structure of $E_{[l]}$ and $A_{[l]}$ and replacing $X = (X/l-1, {\bf x}_{l})$ we can observe,
\begin{eqnarray}
    E_{[l-1]}X/l-1  & = & {\bf 1}_{l-4} \\
     {\bf 1}_{p_{l-1}}{\bf x}_{l} & = & 1 \\
  \left( \begin{array}{c}
A_{[l-1]}\\
\mathbf{0}\\
\end{array}\right)X/l-1 + A^{(l)}{\bf x}_{l} + U^{(l)}& = & U^{(3)} 
    \end{eqnarray}
Or from the definition of $U^{(l-1)}$ from Problem~$MIR(l-1)$,we have
 \begin{equation}
U^{(l-1)} - A^{(l)}{\bf x}_{l} = U^{(l)}, ~ for ~all ~ l, 4 \le
l \le n.
\end{equation}

\begin{prop}[Stem Property]\label{propMI}
Let  $n \geq 5$. If $X \in R^{\tau_{n}}$ is feasible for Problem~$MIR(n)$ implies $X/l$ is feasible for Problems~$MI(l)$ for  $l \in \{4,\ldots, n-1\}.$ 
\end{prop}
\begin{proof}
Suppose the property is not true. Look at the $l$ for which $X/l$ is not feasible for Problem~$MIR(l).$ Then there exits a $U^{(l)}_{ij} < 0$ for some $ij$. This implies $U^{(n)}_{ij} < 0$ as well since $\sum_{k =l+1}^{n}x_{ijk} \geq 0$. So $X$ is not feasible for Problem~$MIR(n)$, leading to a contradiction. So there is no such $l$, or the stem property is true.
\end{proof}

Given $n, Y \in R^{\tau_n}$, using $U^{(l-3)},$ the slack variable vectors in Problem~$MIR(l)$, for $l, \in \{ 4, \ldots, n\}$ (with trailing zeros where needed), we can reformulate the constraints of the $MI$ relaxation problem in matrix notation as follows :

\begin{eqnarray}\label{reform}
E_{[n]}Y  & = & {\bf 1}_{n-3}\label{EnX} \\
U^{(0)} & = & \left( \begin{array}{c}
     {\bf 1}_{3}  \\
     \bf{0}
  \end{array}\label{U0}
  \right) \\
 U^{(l-4)} - A^{(l)}{\bf y}_{l} & = & U^{(l-3)},~ for ~all ~ l,~ 4 \le l \le n.
 \label{UAk+1X}\\
 U & \ge & 0, \label{unonneg} \\
 Y & \ge  & 0. \label{Xnonneg}
\end{eqnarray}

\begin{example}
Consider $n=5$, then $p_{4} = 6, p_{5} = 10,$  and $\tau_{n} = 9.$ Let $X$ be given by,
\[
X'= (1\!/2 , 0 ,  1\!/2;\ 1\!/2 , 0, 0, 0,  0,  1\!/2) \in Q^{\tau_{5}}.
\]
where $X'$ denotes transpose of $X$. 
$X$ satisfies the equality restriction~(\ref{eq:6}) of the $MI$ relaxation. Also, the corresponding slack variable  $u$ obtained from the rest of the inequalities~(\ref{eq:6},\ref{eq:7}  ~\& ~\ref{eq:9}) is $ \geq 0.$  We have,
\[
u' = \left( \begin{array}{c}
											{\bf 1}_{3}\\
0	
\end{array} \right) - A_{[5]}X = (0 ,  1 , 1\!/2 ,  1\!/2 , 1 , 0 , 1\!/2 ,  1\!/2 , 1\!/2 , 1\!/2 ).\nonumber
\]
(Matrices $E_{[5]}$ and $A_{[5]}$ are shown in Figure~\ref{matrices}.)  

So $(X,u)$ is a feasible solution to the $MI$ relaxation. Here, $U^{(1)}$ is given by $(1\!/2, 1,1\!/2, 1\!/2,1, 1\!/2 )'.$ and $U^{(2)} =  u$.
\end{example}\label{ex:one}

\begin{figure}[ht]
\[
E_{[5]} = \left(
\begin {array}{ccccccccc}
1 & 1 & 1 &  0 & 0 & 0 & 0 & 0 & 0  \\
0 & 0 & 0 &  1 & 1 & 1 & 1 & 1 & 1  \\
\end{array}
\right),
\]
\[
A_{[5]} = \left(
\begin {array}{ccccccccc}
\ 1  & \ 0 &  \ 0 &  \ 1 & \ 0 & \ 0 & \ 0 & \ 0 & \ 0  \\
\ 0  & \ 1 &  \ 0 &  \ 0 & \ 1 & \ 0 & \ 0 & \ 0 & \ 0  \\
\ 0  & \ 0 &  \ 1 &  \ 0 & \ 0 & \ 1 & \ 0 & \ 0 & \ 0 \\
-1   & -1  &  \ 0 &  \ 0 & \ 0 & \ 0 & \ 1 & \ 0 & \ 0  \\
-1   & \ 0 & -1   &  \ 0 & \ 0 & \ 0 & \ 0 & \ 1 & \ 0  \\
\ 0  & -1  &  -1  &  \ 0 & \ 0 & \ 0 & \ 0 & \ 0 & \ 1  \\
\ 0  & \ 0 &  \ 0 &  -1  & -1  & \ 0 & -1  & \ 0 & \ 0  \\
\ 0  & \ 0 &  \ 0 &  -1  & \ 0 & -1  & \ 0 & -1  & \ 0  \\
\ 0  & \ 0 &  \ 0 &  \ 0 & -1  & -1  & \ 0 & \ 0 & -1 \\
\ 0  & \ 0 &  \ 0 &  \ 0 & \ 0 & \ 0 & -1  & -1  & -1 \\
\end{array}
\right).
\]

\caption{Matrices $E_{[5]}$ and $A_{[5]}$ }\label{matrices}
\end{figure}

\begin{lemma}\label{obvious}
If $Y \in P_{MI}(n)$ then there exists $U$ satisfying Equations~\ref{U0} through ~\ref{unonneg}. 
Conversely, if $U$ and $Y$ satisfy Equations~\ref{EnX} through ~\ref{Xnonneg} then $Y \in P_{MI}(n)$.
\end{lemma}
\begin{proof}
Follows from Property~\ref{propMI}.
\end{proof}
\begin{remark}
     Given $n$, the $MI-relaxation$ problem is in the combinatorial LP class, and the dimension of the problem is less than $n^{2} \times n^{3}$. 
\end{remark}
\section{Characterization of Pedigree Polytope}\label{char}
In this section,  given a $X \in P_{MI}(n)$, we wish to know
whether $X$ is indeed in $conv(P_{n})$. Results
presented here are from \cite{TSADMpaper}; some proofs are omitted. Let $|P_k|$
denote the cardinality of $P_k$. Assume that the pedigrees in
$\mathcal{P}_k$ are numbered (say,  according to the
lexicographical ordering of the edge labels of the edges appearing in a pedigree.

\begin{definition}\label{realx}
Given $X =({\bf x}_{4},\ldots , {\bf x}_{n}) \in P_{MI}(n)$ recall that we
denote by $X/k = ({\bf x}_{4},\ldots , {\bf x}_{k})$,  the
restriction of $X$, for $4 \leq k \leq n$.

 Given $X \in P_{MI}(n)$ and $X/k \in conv(P_{k})$, consider \mbox{ $\lambda \in
R^{|P_{k}|}_{+}$} that can be used as a weight to express $X/k$ as
a convex combination of $X^{r} \in P_{k}.$  Let $I(\lambda )$
denote the index set of positive coordinates of $\lambda .$  Let
$\Lambda _{k}(X) $ denote the set of all possible {\em weight
vectors}, for a given $X$ and $k$, that is, $$\Lambda _{k}(X) = \{
\lambda \in R^{|P_{k}|}_{+} | \sum_{r \in I(\lambda),X^{r} \in
P_{k}} \lambda_r X^{r} = X/k, \  \sum_{r \in I}\lambda_r =1 \}.$$
\end{definition}

\begin{definition}
Consider a $X \in P_{MI}(n)$ such that $X/k \in conv(P_{k})$. We
denote  the $k-tour$ corresponding to a pedigree $X^{\alpha} $ by
$H^{\alpha}.$ Given a weight vector $\lambda \in \Lambda _{k}(X)$,
we define a $FAT$ problem with the following data:

$$
\begin{array}{lrclr}
O --   &Origins ] & : & \alpha, \alpha \in I(\lambda ) \\
a --            &Supply     ]   &   :   &   a_{\alpha} = \lambda_{\alpha} \\
D --   &Destinations        ] & :   & \beta, e_{\beta} \in E_{k}, x_{k+1}(e_{\beta}) >0 \\
b --        &Demand     ]   &   :   &   b_{\beta} = x_{k+1}(e_{\beta}) \\
\mathcal{A} -- &Arcs            ] & :&  \{ (\alpha,  \beta) \in O \times D | e_{\beta}  \in H^{\alpha}  \} \\
\end{array}
$$

We designate this problem as $FAT_{k}(\lambda )$. Notice that arcs
$(\alpha,  \beta)$ not satisfying $e_{\beta}  \in H^{\alpha}$ are
the forbidden arcs. We also say $FAT_{k}$ is feasible if problem
$FAT_{k}(\lambda )$ is feasible for some $\lambda \in \Lambda
_{k}(X).$
\end{definition}

Equivalently, the arcs in $\mathcal{A}$ can be interpreted as
follows: If $W^{\alpha}$ is the pedigree corresponding to
$X^{\alpha}\in P_k$ for an $\alpha \in I(\lambda )$ then the arcs
$(\alpha, \beta)\in \mathcal{A}$ are such that 
a generator of $e_{\beta}$ is available in the pedigree $W^{\alpha}$ and $e_{\beta}$ does not occur in $W^{\alpha}$. Or $(W^{\alpha},e_{\beta})$ is an extension of $W^{\alpha}.$

\begin{figure}[htb]
    \centering
    \includegraphics[width = 0.6\textwidth]{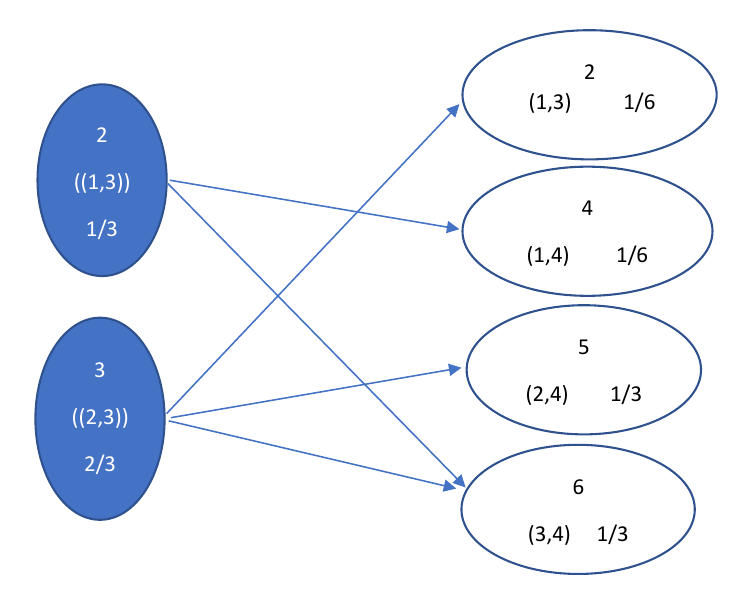}
    \caption{$FAT_4(\lambda)$
Problem for Example~\ref{exampleFat} \\ {\bf Legend:} An origin node gives 
$\begin{array}{c}
pedigree \  no.\\
pedigree\\
availability
\end{array}
$ ; and a destination node gives $\begin{array}{ccc}
  & edge \ label & \\
edge &  \ & demand
\end{array}
$.}
    \label{figure1}
\end{figure}

\begin{example}\label{exampleFat}
Consider $X =(0 , \frac{1}{3} , \frac{2}{3} , 0 , \frac{1}{6} , 0
, \frac{1}{6} , \frac{1}{3} ,  \frac{1}{3})$. We wish to check
whether $X$ is in $conv(P_5)$. It is easy to check that $X$ indeed
satisfies the constraints of $P_{MI}(5)$. Also $X/4 =(0 ,
\frac{1}{3} , \frac{2}{3})$ is obviously in $conv(P_4)$. And
$\Lambda_4(X) = \{(0 , \frac{1}{3} , \frac{2}{3})\}$. Assume that
the pedigrees in $\mathcal{P}_4$ are numbered  such that,
$X^{1}=(1 , 0 , 0), X^{2}=(0 ,  1 , 0) \textrm{ and } X^{3}=(0 , 0
, 1)$ and the edges in $E_4$ are numbered according to their edge
labels. Then $I(\lambda) = \{2, 3 \}$. Here $k=4$ and the
$FAT_4(\lambda)$ is given by a problem with origins, $O = \{ 2,
3\}$ with supply $a_2 = \frac{1}{3}, a_3 = \frac{2}{3}$ and
destinations, $D = \{2, 4, 5, 6 \}$ with demand $b_2 =  b_4 =
\frac{1}{6}, b_5 = b_6 = \frac{1}{3}$. Corresponding to origin $2$
we have the pedigree $W^2 = ((1, 3))$. And the edge corresponding
to destination $2$ is $e_2 = (1, 3)$. As $(W^2, e_2)$ is not an
extension of $W^2$, we do not have an arc from origin $2$ to
destination $2$. Similarly, $(W^3, e_4)$, $(W^2, e_5)$ are not
extensions of $W^3$ and $W^2$ respectively, so we do not have arcs
from origin $3$ to destination $4$ and origin $2$ to destination
$5$. We have the set of arcs given by,
$$\mathcal{A} = \{( 2, 4), (2, 6), (3, 1),
(3, 5) \textrm{ and } (3, 6) \}.$$  Notice that $f$ given by
$f_{24} = f_{26}= f_{32}=  f_{36}= \frac{1}{6}, f_{35}
=\frac{1}{3}$ is feasible to $FAT_4(\lambda)$.  (See
Figure~\ref{figure1}, where the nodes have pedigree numbers or edge labels. shown, we only need the pedigrees and edges for solving the $FAT_{k}(\lambda)$ problem.)  This $f$ gives a weight vector to
express $X$ as a convex combination of the vectors in $P_5$, which
are the extensions corresponding to arcs with positive flow. This
role of $f$  is, in general, true, and we state this as
Theorem~\ref{flowisweightvector}.

Thus, we have expressed $X$ as a convex combination of the incidence vectors of the pedigrees $W^{7}((1, 3)(1, 4)),
W^{8}=((1, 3)( 3, 4))$ , $W^{10}= ((2, 3)( 1, 3))$, and  $W^{11}=((2, 3)(3, 4))$, each of them receive a weight of $\frac{1}{6}$,  and $W^{12}=((2, 3)( 2, 4))$, which receives a weight of
$\frac{1}{3}$.

It is easy to check that $f$ is the unique feasible flow in this
example, so no other weight vector exists to certify $X$ in
$conv(P_5)$. In other words, all the arcs are rigid, and there are no dummy arcs.
\end{example}

\begin{theorem} \label{flowisweightvector}
Let $k \in [4, n-1]$. Suppose $\lambda \in \Lambda
_{k}(X) $ is such that  $FAT_{k}(\lambda )$ is feasible. Consider
any feasible flow $f$ for the problem.  Let $W^{\alpha \beta}$
denote the extension $(W^{\alpha} , e_{\beta})$, a pedigree in
$\mathcal{P}_{k+1}$, corresponding to the arc $( \alpha,\beta )$.
Let $\mathcal{W}_f$ be the set of such pedigrees, $W^{\alpha
\beta}$ with positive flow $f_{{\alpha \beta}}$. Then $f$ provides
a weight vector to express $X/k+1$ as a convex combination of
pedigrees in $\mathcal{W}_f$.
\end{theorem}

Next we observe that $conv(P_n)$ can be characterized using
a sequence of flow feasibility problems as stated in the following
theorems:

\begin{theorem} \label{FF}
If $X \in conv(P_{n })$ then    $FAT_{k}$ is feasible $\allk .$
\end{theorem}

\begin{theorem} \label{fat}
Let $k \in [4, n-1]$. If  $\lambda \in \Lambda
_{k}(X) $ is such that  $FAT_{k}(\lambda )$ is feasible, then
$X/(k+1) \in conv(P_{k+1}).$
\end{theorem}

Proofs of these theorems are given in the appendix of~\cite{TSADMpaper}. In general,
we do not have to explicitly deal with the set $\Lambda_k(X)$. The set
is used in the proofs only.  

Thus, we can state the following conclusion from the above discussion as Theorem~\ref {nandccondition1}:

\begin{theorem}\label{nandccondition1}
Given $n, X \in P_{MI}(n), X/n-1 \in conv(P_{n-1}),$ let the $FAT$ problem called $FAT_{n-1}(\lambda),$ for a $\lambda \in \Lambda_{n-1}(X)$,  be defined as follows:
The set of origins $\mathcal{O}$ = $I(\lambda),$ and set of destinations,$\mathcal{D} =  V_{[n-3]}.$ And the set of arcs\footnote{$\mathcal{A}$ can be defined equivalently using $H^{r}$, the $n-1$ tour of $P^{r}.$ } is given by $\mathcal{A}  = \{ (P^{r} , e), r \in I(\lambda) | \textrm{ a generator of }  e \in P^{r} and e \notin P^{r} \}.$

A necessary and sufficient condition for the membership of $X$ in the pedigree polytope, $conv(P_{n})$ is $FAT_{n-1}$ is feasible.
\end{theorem}

In Theorem~\ref{fat}, we have a procedure to check whether a given
$X \in P_{MI}(n)$, is in the pedigree polytope, $conv(P_{n})$.
Since feasibility of a $FAT_{k}(\lambda)$ problem for a weight
vector $\lambda$ implies $X/(k+1)$ is in $conv(P_{k+1})$, we can
sequentially solve $FAT_{k}(\lambda_{k} )$ for each $k = 4,
\ldots, n-1$ and if $FAT_{k}(\lambda_{k} )$ is feasible, we set $k = k+1$ and while  $k < n$  we repeat; at any stage, if the problem is infeasible, we stop. So if we have reached $k=n$, we have proof that $X \in conv(P_{n}).$ But finding a suitable $\lambda_{k}$ such that $FAT_{k}(\lambda_{k})$ is feasible to show that $X/k \in conv(P_{k}$ is not easy and the complexity of finding such a $\lambda_{k}$ has not been studied in the literature. However, we devote the rest of this article to  explaining the construction of a layered network, stating and proving a necessary and sufficient condition answering the membership question of  $M3P$, and  establishing the complexity of checking that condition.

\section{Construction of  the Layered Network}\label{conslay}
In this section we define the layered network $N_k$ and a set of pedigrees $R_k$ with fixed positive weight $\mu_P$ for $P \in R_k$,
with respect to a given $X \in P_{MI}(n)$ and for $k \in
V_{n-1}\setminus V_{3}.$ Given  $X/k \in conv(P_{k})$,  this
network is used in showing whether or
not $X/k+1 \in conv(P_{k+1}).$

We denote the node set of $N_{k}$ by $\mathcal{V} (N_{k})$ and the
arc set by $\mathcal{A}(N_{k})$. Let $v = [k:e] $ denote a node
in the $(k-3)^{rd}$ layer corresponding to an edge $e \in
E_{k-1}$. Let $x(v) = x_k(e)$ for $v = [k:e]$.
\subsection{ Construction of  the Network for \texorpdfstring{ $k = 4$}{k = 4}}\label{consfour}
  However, as shown in section~\ref{char}, $X \in P_{MI}(n),$ implies $X/n-1  \in P_{MI}(n-1).$ Therefore, if $X/5 \notin conv(P_{5})$ we can conclude that $X \notin conv(P_{n}).$ $FAT_{4}$ feasibility is both necessary and sufficient for $X/5 \in conv(P_{5}).$ 

Let
$$V_{[r]} = \{ v | v =[r+3: e], e \in E_{r+2}, x(v) > 0 \}.$$
The node name $[r+3: e]$ alludes to the insertion decision
corresponding to the stage $r$; that is, the edge $e$ is used for the insertion of $r+3$. First, we define the nodes in the network
$N_{4}$.
$$\mathcal{V}(N_{4}) = V_{[1]} \bigcup V_{[2]}.$$
And
$$\mathcal{A}(N_{4}) = \{ (u,v)| u =[4: e_{\alpha}] \in V_{[1]}, v = [5: e_{\beta}] \in V_{[2]}, e_{\alpha} \in G(e_{\beta}) \}.$$

Let $x(w)$ be the capacity on a node $w \in V_{[r]}, r =1,2$.
Capacity on an arc $(u,v) \in \mathcal{A}(N_{4}) $ is $x(u)$.

{\bf Caution:} Notice that nodes in $V_{[1]}, V_{[2]}$ are the same as the origin and destination nodes in $FAT_{4}$ problem, for consistency with subsequent development of the layered network for other $k >4$, we choose these notations, and refer to this problem as $F_{4}$.

Given this network, we consider a flow feasibility problem of
finding a nonnegative flow defined on the arcs that saturates all
the node capacities and violates no arc capacity, that is, solving the $FAT_{4}$ problem.  
Thus, the construction of the layered network proceeds from solving a bipartite flow feasibility problem for stage $1$, involving layers $1$ and $2$. If $FAT_{4}$ is infeasible, we conclude $X \notin conv(P_{n})$ and the construction of the network stops. 
When $FAT_{4}$ is feasible we use the $FFF$ algorithm (or any such)
and identify the set of rigid arcs, $\mathcal{R}$. If there are any dummy arcs in
$\mathcal{R}$, we delete them from the set of arcs and
update. For every rigid arc with a positive frozen
flow, we replace the capacity of the arc with the frozen flow. And
include the corresponding arc (pedigree, $P \in conv(P_{5})$  in the set $R_4$ with $\mu_P$, the frozen flow as the weight for $P$. 
We reduce the node capacities of the head and tail nodes of the rigid arc $P$ by $\mu_P$, for each $P \in R_4$, and delete nodes with zero capacity, along with the arcs incident with them. Consider the sub-network obtained as $N_{4}$. Let $\bar{x}(v)$ be the updated capacities of node $v$ in the network. Now we say $(N_{4}, R_{4}, \mu)$ is {\em well-defined}. Observe that either the node set in $N_{4}$ is $\emptyset$ or $R_4 = \emptyset$ , but not both. 

\begin{example}
Consider $X \in P_{MI}(5)$ given by
$
\begin{array}{rcl}
{\bf x}_{4} & = & (0, 3/4, 1/4); \\
{\bf x}_{5} & = & (1/2, 0, 0, 0, 0, 1/2); \\
\end{array}$ It can be verified that $FAT_{4}$ is
feasible and $f$ given by
$f_{([4:1,3],[5:1,2])} = 1/4, f_{([4:2,3],[5:1,2])} = 1/4, \textrm{ and } f_{([4:1,3],[5:3,4])} = 1/2$ does it. And so, we conclude $X \in conv(P_{5}).$ But no arc is rigid in $FAT_{4}$, is obvious.
Therefore, $(N_{4} = F_{4}, R_{4} = \emptyset, \mu = 0).$
However, $X^{\prime}$ given by
$$
\begin{array}{rcl}
{\bf x^{\prime}}_{4} & = & (0, 3/4, 1/4); \\
{\bf x^{\prime}}_{5} & = & (0, 1/4, 0, 0, 1/4, 1/2); \\
\end{array}
$$
is such that, $X^{\prime} \in P_{MI}(5)$ but $FAT_{4}$ is not feasible. 
So $X^{\prime} \notin conv(P_{5}).$  
\end{example}
\begin{example}
Consider Example~\ref{exampleFat}, we have $FAT_{4}$ feasible but all arcs are rigid and so
$(N_{4} = \emptyset, R_{4} = \{ ((1, 3)(1, 4)),
((1, 3)( 3, 4)) , ((2, 3)( 1, 3)), ((2, 3)(
3, 4)),((2, 3)( 2, 4)) \}, \mu = (\frac{1}{6}, \frac{1}{6},\frac{1}{6}, \frac{1}{6},\frac{1}{3} )
$
\end{example}
\subsection{Overview of the  Membership Checking in  Pedigree Polytope}
If $(N_{k-1},R_{k-1}, \mu)$ is well defined and $k < n$ we proceed further.  
Next, consider a bipartite flow feasibility problem for stage $k-3$, called $F_{k}$. If $F_{k}$ is infeasible, we conclude $X \notin conv(P_{n}),$ and the construction of the network stops. Otherwise, we proceed to construct $(N_{k}, R_{k}, \mu)$. Next, we seek a feasible solution to a multicommodity flow problem. If we have a feasible solution, we have evidence that $X/k \in conv(P_{k}),$  and we say the network is well-defined. If  $(N_{k}, R_{k}, \mu)$ is not well-defined, we conclude $X \notin conv(P_{n}),$ and stop. Otherwise proceed until we reach $k= n-1,$  constructing $(N_{k},R_{k},\mu),$ and concluding $X \in conv(P_{n})$ or not.

\begin{figure}[ht]
    \centering
    \includegraphics[scale =0.9]{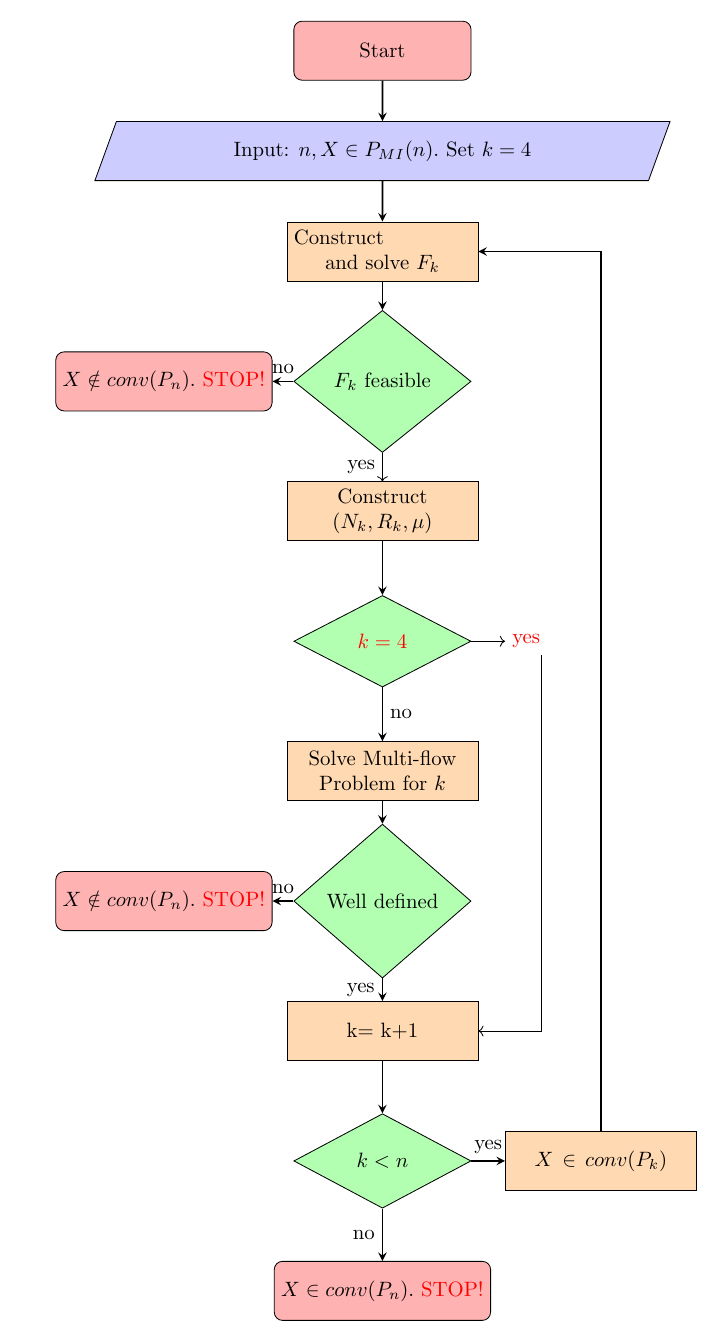}
    \caption{Flowchart showing the Membership Checking in $conv(P_{n})$}
    \label{fig1}
\end{figure}

This construction process is explained in Figure~\ref{fig1}, and the following discussion elaborates on these steps 

and proves necessary results to ensure that the conclusion made on the status of $X$’s membership in the pedigree polytope is valid.
\subsection{\texorpdfstring{Construction of the Layered Network: $N_{k}, k >  4$ }{Construction  of the Layered Network:  Nk, k > 4}}
Primarily, we are checking whether some pedigrees in $P_{k}$ that are active for $P \in X/k$ can be extended to pedigrees in $P_{k+1}.$  This is easy for a pedigree in $R_{k-1}$, all that we need to check is whether $(P, v)$ for a $v \in \Delta^{k+1}$ is a pedigree in $P_{k+1}.$ On the other hand, with respect to $N_{k-1}$, we have paths that correspond to pedigrees, and some of them are active for $X/k$. We need to ensure only such active pedigrees are considered for extending to pedigrees in $P_{k+1}.$
Here we have two tasks to perform, [1] constructing $F_{k}$ and checking feasibility of $F_{k},$ (if infeasible, provides a sufficient condition for $X/k+1 \notin P_{k+1}$) and [2] if $F_{k}$ is feasible, constructing $(N_{k}, R_{k}, \mu)$ and checking whether the multicommodity flow problem $Mcf(k)$ is feasible, completes the membership checking for $ X/k+1$ in $conv(P_{k+1})$. 
\subsubsection{\texorpdfstring{Construction of the FAT problem  $F_{k}, k  >  4$}{ Construction of the FAT problem  Fk, k  >  4}}
The set of origins, $Origins$, for the bipartite graph of $F_{k}$ includes the last stage nodes in $N_{k-1}$ with updated node capacities as availabilities. In addition, for each $P \in R_{k-1}$ we introduce a \textit{ shrunk } node $P$ as origin, with availability equal to $\mu_{P}.$ The set of destinations, $Destinations$, for the bipartite graph of $F_{k}$ include the triangles $v$ with demands  $x_{v} > 0, v \in \Delta^{k+1}.$ 

The arcs in this problem are constructed depending on whether we can extend a pedigree in $P_{k}$ active for $X/k$ to a pedigree in $P_{k+1}.$ 

Consider a link $L = (u = P, v= \{r,s,k+1\}\})$ connecting a pedigree in $R_{k-1}$  and a destination $v$. If there is a generator of $v $ in $P$ and $\{r,s\} $ is not one of the common edges of $P$, we have an arc $(u, v),$ with capacity $\mu_{P}.$

Let the set of nodes in $N_{k-1}$ be denoted by $\mathcal{V}(N_{k-1}).$ and the arcs in $N_{k-1}$ be $\mathcal{A}(N_{k-1}).$ Consider a link $L = (u =\{r, s, k\}, v= \{i,j,k+1\})$ connecting an origin $u$ and a destination $v$.
In this case,  paths that are to be \textit{ avoided } in  $N_{k-1}$ have the 
\textbf{Requirements:} 
\begin{itemize}[label=]
\item{1)}  do not contain a generator $u \textrm{ or }v$, 
\item{2)}  contain $u \textrm{ or }v,$  and 
\item{3)}  do not end in $u$. 
 \end{itemize} 
 We define a \textit{ restricted network } $N_{k-1}(L), $ deleting certain nodes from $N_{k-1}$ and pedigrees in $R_{l}, l\in [4, k-2]$, to achieve the above requirements.fin

\begin{definition}{[Restricted Network
$N_{k-1}(L)$]}\label{restrict} 
Given $k \in[5, n-1]$, a link $L= (u= \{r,s, k\}, v =\{i,j, k+1\})$.
$N_{k-1}(L)$ is the subnetwork induced by the subset of nodes
$\mathcal{V}(N_{k-1}) \setminus \bf{D}$, where 
$\bf{D}$, the set of deleted nodes, is constructed as
follows: 

Let $\bf{D}= \emptyset.$
\begin{itemize}[label=]
    \item{ a)}  Include $\{i, j,l\}$ in $\bf{D}$, 
    for $l \in [\max(4, j) , k-1].$ 
    \item{ b)}  Include $\{r,s,l\}$ in $\bf{D}$,
    for $l \in [\max(4, s) , k-1].$ 
    \item{ c)}  Include $w \in \Delta^{s}, w \notin
    G(u)$ in $\bf{D}.$
    \item{ d)}   Include $w \in \Delta^{j}, w \notin
    G(v)$ in $\bf{D}.$ 
     \item{ e)}
    Include   all nodes in $\Delta^{k} \setminus \{u\}$ in
    $\bf{D}$. %
    \item{ f)}
    If any undeleted node $\{y, w, l\}, l > 4$ is such that all its generators have been deleted, include that node in $\bf{D}$. 

    Repeat this step until no more nodes are deleted. Set $\mathcal{V}(N_{k-1}(L)) = \mathcal{V}(N_{k-1}) \setminus \bf{D}.$
    \item{ g)} If a node from $\bf{D}$ appears in any of the rigid pedigrees for any  $l$, such pedigrees are deleted. This means the corresponding (shrunk) nodes are deleted from $\mathcal{V}(N_{k-1}(L))$.  
\end{itemize}

The sub-network induced by
$\mathcal{V}(N_{k-1}(L))$ is called the \emph{Restricted Network}
$N_{k-1}(L)$. 
\end{definition}
Recall that a vertex-induced sub-network is a subset of the nodes of a network together with any edges whose both endpoints are in this subset. 

\begin{remark}
A few observations that are easy to make are listed below:

\begin{itemize}[label=]
\item{[1] }  Deletion rule [a], [b] ensure requirement (2). Deletion rule [c]
and [d] ensure requirement (1). Deletion rule [e] ensures
requirement (3). Rule [f] ensures each node in the restricted
network has at least one generator. 
Lastly, [g] ensures that in any stage $l$, a rigid pedigree is deleted from $R_{l}$, if it (the corresponding pedigree path) contains a deleted node. 
So the corresponding (shrunken) node is deleted from the restricted network.
\item{[2] }  
Construction of the restricted network for any link can be done by applying deletion rules $[a]$ through $[e]$, building the set  $\bf{D}$. 
Then we need only to check whether any node in $\mathcal{V}(N_{k-1}) \setminus \bf{D}$ is deleted using the deletion rule $[f]$. We check the undeleted nodes layer by layer. If any node is deleted as per rule $[f]$, it can only affect the generator availability status of nodes in the subsequent layers and not in the current layer or the layers already checked applying rule $f$. So while applying rule $[f]$ whenever $\bf{D}$ is updated, it is sufficient to check rule $[f]$ for the remaining undeleted nodes. 

Therefore, to obtain finally $\mathcal{V}(N_{k-1})(L)$, deletion rule $f$ is repeated at most  $|\mathcal{V}(N_{k-1})|$ times. 
\end{itemize}
\end{remark}

Once $N_{k-1}(L)$ is created, we consider the problem of finding the maximal flow in $N_{k-1}(L)$ satisfying non-negativity of flow, flow conservation at nodes, and capacity restrictions on the available nodes and arcs. The only sink in the network is $\{r,s,k\}$, and the sources
are the undeleted nodes in $\Delta^{4}$,  and the undeleted pedigrees in $R_{l},  l \in [4, k-2]$. Let $C(L)$ be the value of the maximal flow in the restricted network $N_{k-1}(L)$. 
We find $C(L)$ for each link $L$. 
Besides,  if the path that brings this flow $C(L)$ for any $L$ is unique, we save the corresponding path, named $P_{unique}(L)$, against $L$. 
Notice that the complexity of finding a maximal flow in the network  $N_{k-1}(L)$ is strongly polynomial in the number of nodes and edges in the network \citep{KandV}.

\begin{problem}[FAT problem, $F_{k}$]
Given $k \in [5, n-1]$, the forbidden arc transportation problem maximising the flow in the bipartite network~\footnote{We also refer to the bipartite network by the same name, $F_{k}$.} with $
\begin{array}{lrclr}\label{defFk}
 & & \\
O --   &Origins ] & : & u = \{r,s,k\} \textrm{ and } P \in R_{k-1} \\
a-- &Availability ] & : & \bar{x}(u) \textrm{ for } u \textrm{ and } \mu_P \textrm{ for } P \\
D --   &Sinks           ] & :   & v = \{i, j, k+1\} \\
b-- &Demand ] & : & x(v) \textrm{ for } v \\
\mathcal{A} -- &Arcs            ] & :&  \{ (u, v) \textrm{ such that, either} \  L = (u, v) \textrm{ is a link and } \ C(L) > 0 \\
  & ] & : & \textrm{ or }  u = P \textrm{ and }  (P, v) \in P_{k+1} \}  \\
C -- &Capacity] &:& C_{u, v} = C(L) \textrm{ or } \mu_P.\\
\end{array}
$
 
 is called the problem $F_{k}$. 

\end{problem}

If $F_{k}$ is feasible we apply the $Frozen \ Flow
  \ Finding$ algorithm and identify the set $\mathcal{R}$ and the dummy subset
of arcs in that. Update $\mathcal{A}$ by deleting the dummy arcs. Update the capacity of any rigid arc,  $C(L)$, to the frozen flow. 
\subsection{Completing the construction of the Layered Network}\label{DefRk}
Here we are in a position to complete the construction of the network $(N_{k}, R_{k}, \mu)$. Construction of $R_{k}$ is driven by the rigid arcs in $F_{k}$. If the rigid flow is from a unique path, then that contributes a rigid pedigree to $R_{k}$. Similarly, if the rigid flow is from a shrunken node, $P \in R_{k-1}$, then that contributes a rigid pedigree to $R_{k}$. Let $L = (u = [k:e_{\alpha}], v = [k+1:e_{\beta}])$ be a rigid arc in $F_{k}$.
Given  $P_{unique}(L)$ for any such rigid arc $L$, we include the extended pedigree $(P_{unique}(L), e_{\beta})$ to $R_k$. 
 Update the node capacity of any node $u \in \mathcal{V} (N_{k})$ that appears in some rigid pedigree paths, that is, $path(P), P \in R_k,$ by reducing the node capacity of $u$ by the rigid flows ($\mu_P$'s) along these paths. 
Similarly, for any rigid arc, connecting a pedigree $P^{\prime} \in R_{k-1}$, we include the extended pedigree $(P^{\prime}, e_{\beta})$ to $R_{k}$. Set  $\mu_P$ for $P \in R_{k}$ equal to the rigid flow of the defining rigid arcs.
Update the node capacity of any node $u= P^{\prime} \in R_{k-1}$ by reducing the node capacity of $u$ by the sum of the flows along the defining rigid arcs. 
Similarly, we update the arc capacities along the unique pedigree path, $P_{unique}(L)$. 
Update $\mathcal{A}$ by deleting the rigid arcs defining pedigrees in $R_{k}$. We can assume without loss of generality that we have distinct pedigrees in $R_{k}$. Otherwise, we can merge the weights of multiple occurrences of a pedigree and ensure that the pedigrees are all distinct. Finally, obtain $N_{k}$ by including the remaining arcs and the destination nodes in the bipartite network for $F_{k}$, to $N_{k-1}$. Thus, we have $(N_{k}, R_{k}, \mu)$  constructed; Section~\ref{multi} will discuss how we check that it is well-defined.

\begin{example}
    Given  $X = (1, 0 \  \frac{1}{2}  \ \frac{1}{2}, \frac{1}{4}  \ \frac{1}{4}  \  \frac{1}{8}  \ \frac{1}{8}\ \frac{1}{8} \ \frac{1}{8}, 0 \ \frac{1}{2} \ \frac{1}{4} \ 0 \ 0 \ 0 \ 0 \ 0 \  0 \ \frac{1}{8} \ \frac{1}{8}  ) ,$ construct and solve $F_{5}.$
\end{example}\label{ex2}

\begin{figure}[htb]
        \centering
        \includegraphics[width =0.7\textwidth]{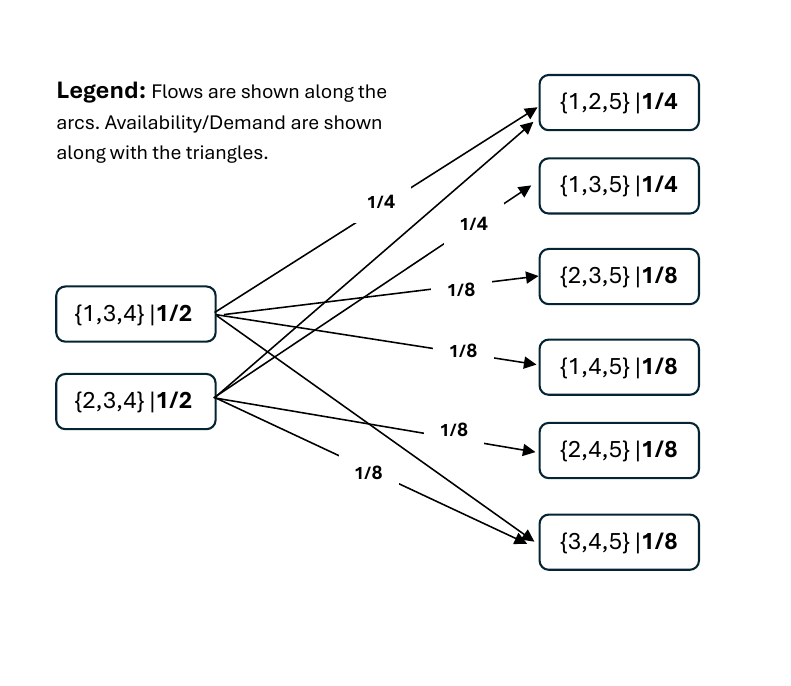}
        \caption{$F_{4}$ feasible \ for \ Example~\ref{ex2}}.
        \label{fig:ex1}
\end{figure}

\begin{figure}[htb]
    \centering
    \includegraphics[]{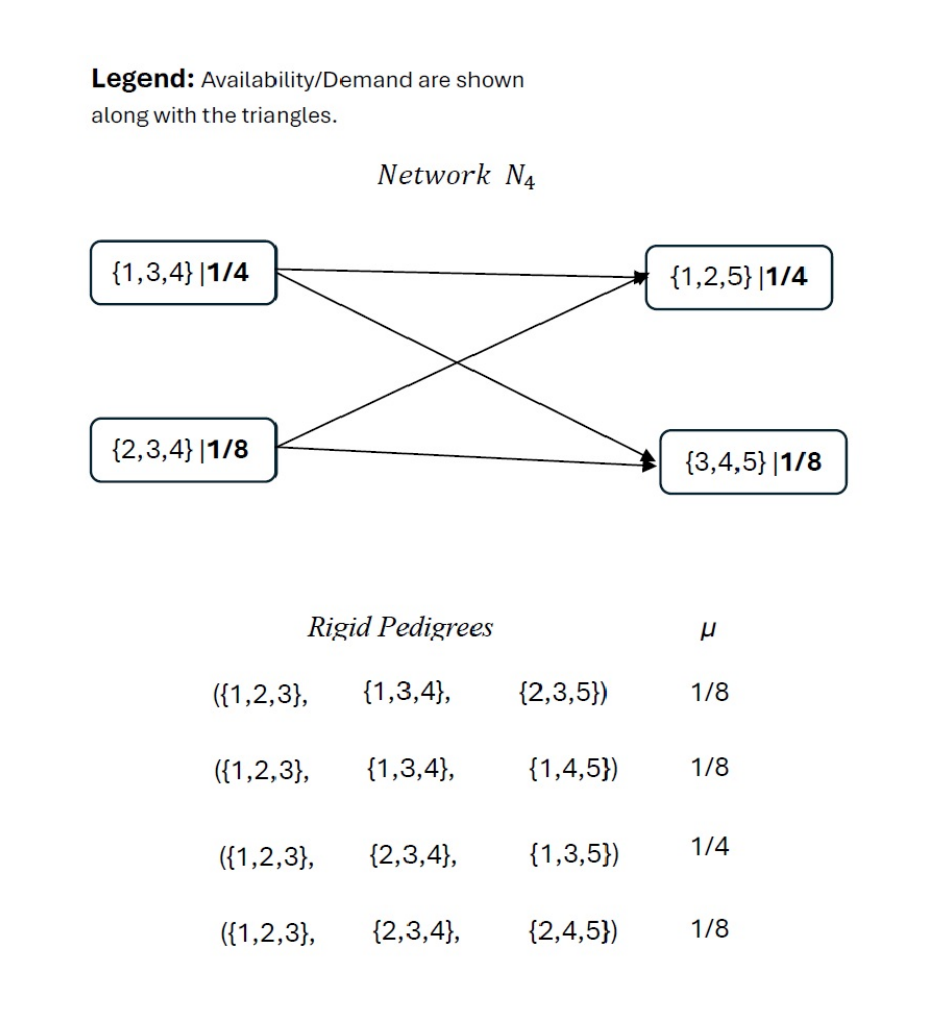}
    \caption{($N_{4}, R_{4}, \mu)$) \ for \ Example~\ref{ex2}. $R_{4}$ has $4$ Pedigrees in $P_{5}.$}
    \label{fig:n4r4muEx1}
\end{figure}

Easy to check that $X \in P_{MI}(6)$.  Now the bipartite graph for $F_{4}$ problem is shown in Figure~\ref{fig:ex1}, along with a feasible flow.
In this example, the rigid arcs are easy to identify, namely, $$\mathcal{R} = \{ ( \{1,3,4 \}, \{2,3,5 \}), (\{1,3,4 \}, \{1,4,5 \}), (\{2,3,4 \}, \{1,3,5 \}), (\{2,3,4 \}, \{2,4,5 \})\}.$$ There are no dummy arcs. And the rest of the arcs form a flow change (fc-cycle), namely, $$(\{1,2,5 \} \frac{+\theta}{\rightarrow} \{1,3,4 \} \frac{-\theta}{\rightarrow} \{3,4,5\} \frac{+\theta}{\rightarrow} \{2,3,4\} \frac{-\theta}{\rightarrow} \{1,2,5 \})$$. 
And there are no bridges in this strongly connected component. (See Appendix A for definitions of fc-cycle, bridge, rigid and dummy arcs.) Figure~\ref{fig:n4r4muEx1} shows $N_{4}, R_{4}$ with $\mu_{P}$ along each rigid pedigree. Since $n =6$, we proceed further.
 So, we have checked $(N_{4},  R_{4}, \mu)$ is well-defined. The origins for $F_{5} $ are  $\{1,2,5\}, \{3,4,5\}$ with availabilities $1/4, 1/8$ respectively and the four rigid pedigrees, $P \in R_{4}$ with availabilities given by respective $\mu_{P}.$ Destinations include $\{1,3,6\}, \{1,4,6\},\{3,5,6\},$ and $ \{4,5,6\}$ with demands $1/2,1/4,1/8$, and $1/8$ respectively. 

For each link $L = (u \in R_{4}, v \in Destinations)$ we check $(P, v) \in P_{6}$ if so we have capacity equal to $\mu_{P}.$ For instance, $P = ( \{1,3,4 \}, \{2,3,5 \})$ with $v = \{1,3,6\}$ is not a pedigree as the common edge $\{1,3\}$ is repeated, so we do not have an arc $(u, v).$ But  $P = ( \{1,3,4 \}, \{2,3,5 \})$ with $v = \{1,4,6\}$ is  a pedigree and we have the corresponding arc $(u, v), $ with capacity $1/8.$ 

\begin{figure}[hbt]
    \centering
    \includegraphics[width=0.9\linewidth ]{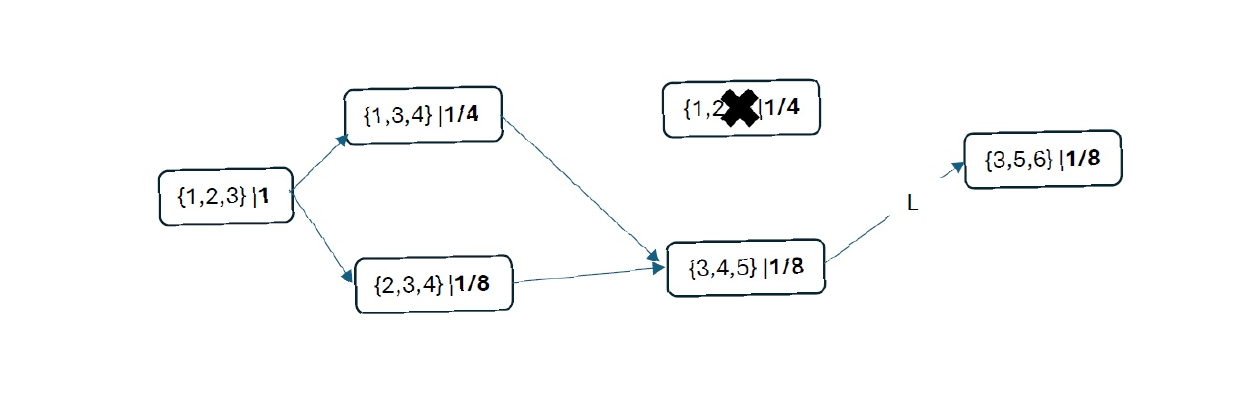}
    \caption{Restricted Network \  for  \ $L =  (\{3,4,5\}, \{3,5,6\})$, Example~\ref{ex2}.}
    \label{fig:N_4_LEx2}
\end{figure}

\begin{figure}[htb]
    \centering
    \includegraphics[width=0.8\textwidth]{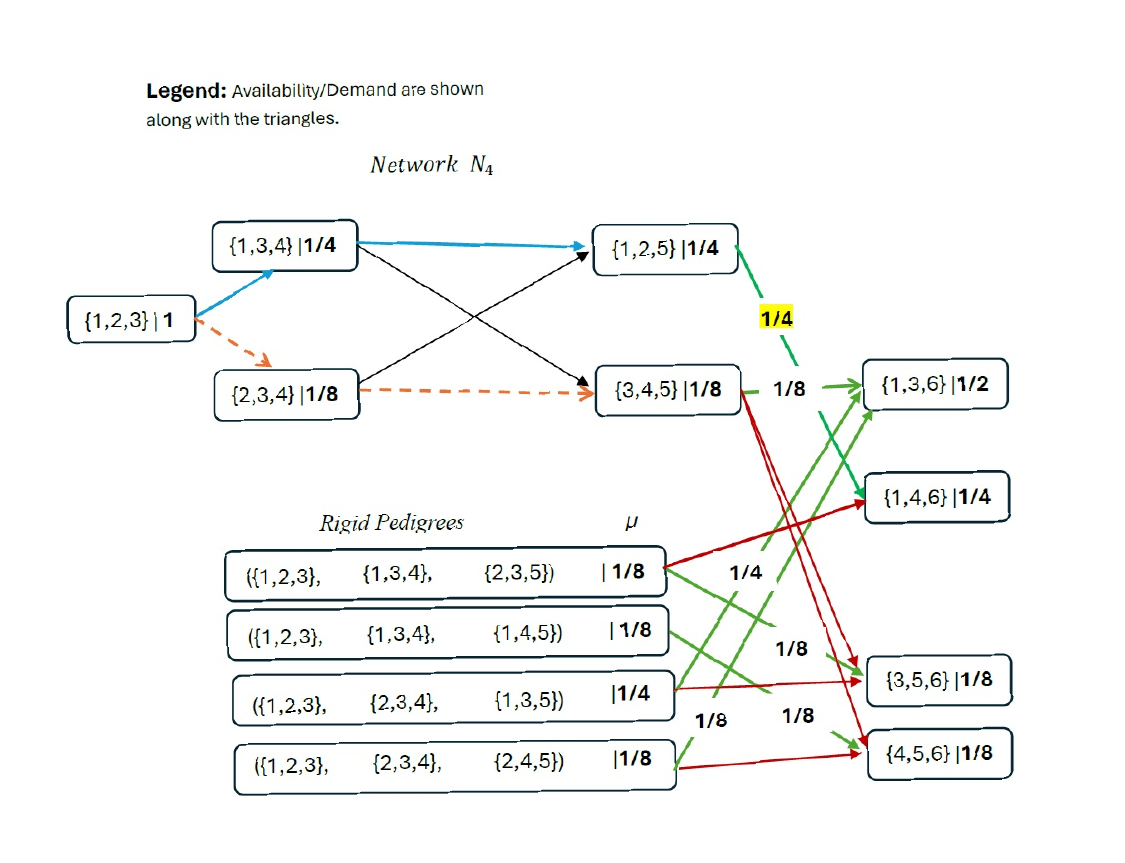}
    \caption{$F_{5}$ \ for \ Example~\ref{ex2}. The unique paths in $N_{4}(L)$ are also shown for \ some $L$ in respective colours.}
    \label{fig:F_5Ex2}
\end{figure}

For other links $l = (u \in \Delta^{5} \ \&  \ \bar{x}(u) > 0, v)$ we construct $N_{4}(L)$ and find the max-flow $C(L)$. For instance, $L = (\{1,2,5\}, \{1,3,6\})$ we have the network, $\{1, 2, 3\} \rightarrow \{2, 3, 4\} \rightarrow \{1,2,5\}$ and maximum flow is $C(L) = 1/8.$  The path is unique, so we store this path $path(L)$ as well for $L$. But $L = (\{1,2,5\}, \{3,5,6\})$ results in an empty network as $\{3,5, 6\}$ has no generator in the restricted network. But $L = (\{3,4,5\}, \{3,5,6\})$ results in the network shown in Figure~\ref{fig:N_4_LEx2},  and maximum flow is $C(L) = 1/8.$ The completed bipartite network for $F_{5}$ with a feasible flow is shown in Figure~\ref{fig:F_5Ex2}.

Since $F_{5}$ is feasible, we could apply the $FFF$ algorithm and identify the set of rigid/dummy arcs in $F_{5}$. However, since the flow $f$ is the unique feasible flow, all arcs are rigid. Arcs with zero flow are dummy arcs, and we can delete them. Other arcs are rigid arcs with frozen flow. Since two of the rigid arcs have a unique path as shown in Figure~\ref{fig:F_5Ex2}, we get two rigid pedigrees in $P_{6}.$ The other rigid arcs are connected to rigid pedigrees in $R_{4}$ and so we get the other rigid pedigrees in $P_{6}$. Thus, we have expressed $X$ as a convex combination of pedigrees in $P_{6}$ or $X \in conv(P_{6})$. In this example, $F_{5}$ feasibility was sufficient for membership in $conv(P_{6})$. In general,  we may have to construct $(N_5, R_5, \mu)$  and solve the multicommodity problem $M CF(6)$ to decide on membership. However, we have
\begin{theorem}\label{nonmemberSuf}
    Given $n > 3, X \in P_{MI}(n)$, if $F_{k}$ is not feasible for $k \in [5, n]$, then $X/k+1 \notin conv(P_{k+1} $.
\end{theorem}

\begin{figure}[htb]
\begin{center}
\includegraphics[angle =90, trim= 0mm 2cm 0mm 5mm, clip=true, scale =.7]{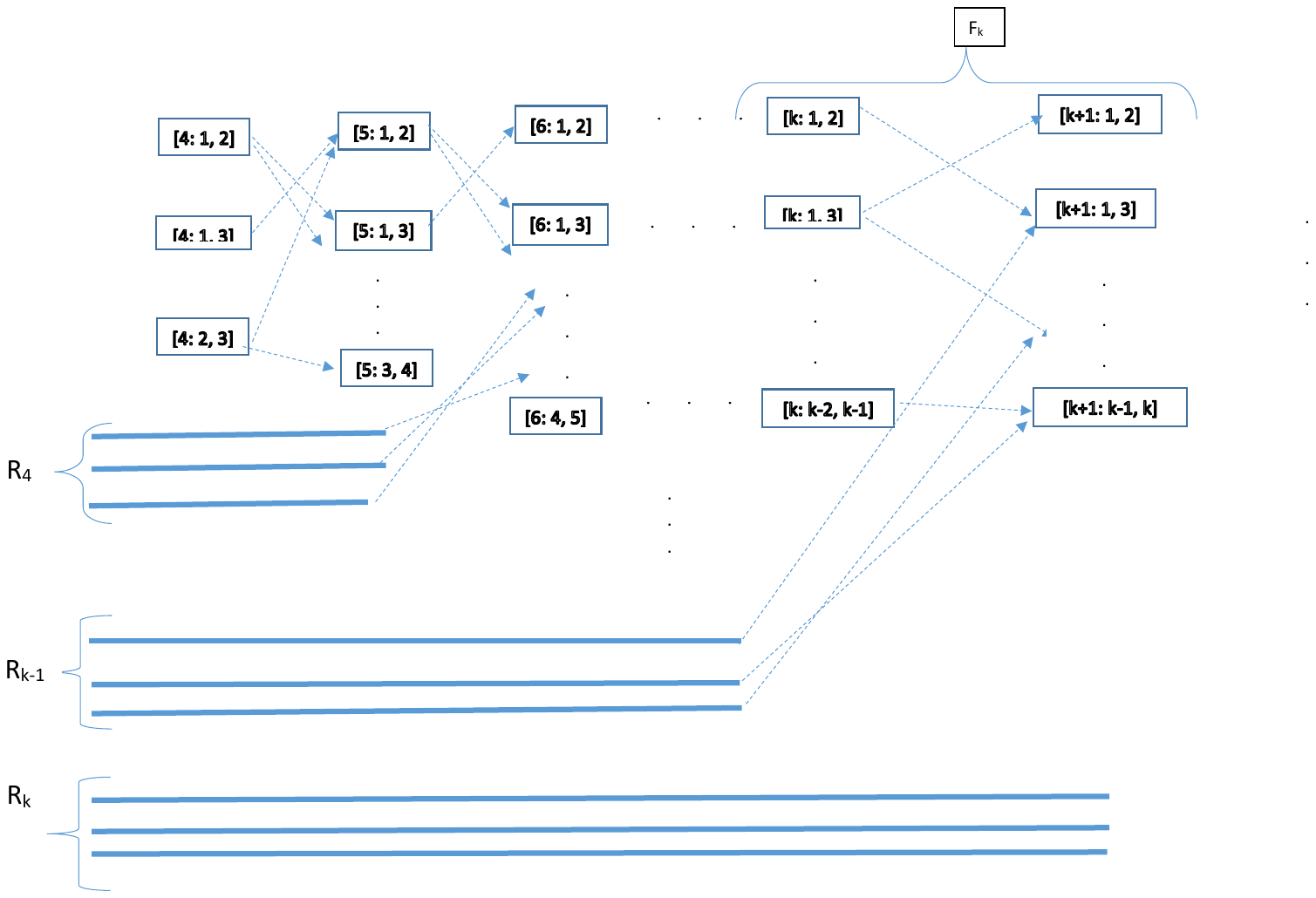}
\end{center}
\caption{Schematic diagram of the network,$N_{k}$}\label{fig:SchemeNk}
\end{figure}

 A schematic diagram of network $N_{k}$ is shown in Figure~\ref{fig:SchemeNk}. In subsection~\ref{ppackability}, we show that  checking whether
 $(N_{k}, R_{k},\mu)$ is well-defined in some special cases, depending on $X/k+1$. 
\section{A Sufficient Condition for Non-membership}\label{Fkfeasibility}
Suppose $X/k+1 \in conv(P_{k+1})$. Consider any active pedigree $P*$ for $X/k+1$. We need to ensure that our deletion rules used in the construction of $N_{k-1}(L)$ for any link $L$ do not discard any such active pedigree.    This section firstly validates our construction of the restricted networks. Then it brings out the importance of studying problem $F_{k}$.

For this purpose, we define another $FAT$ problem ($INST(\lambda,l)$ for which an obvious feasible flow, called the \textit{instant flow}, is available, and this flow can be very useful in proving the feasibility of $F_k$. This problem is used only in the proofs of certain results and is not required to be solved, like $F_{k}$, for a given instance. (see Appendix~\ref{Appendixresults}, Subsection~\ref{resultsinstant})
The main result of this section:
\begin{theorem}{[Theorem on non-Membership]}\label{infeasiblity}
Given $X \in P_{MI}(n)$, and for a  $ k \in V_{n-1} \setminus
V_{3}$, if $X/k \in conv(P_{k}),$ then $$\ F_{k} \textrm{
infeasible implies}  \ X/k+1 \notin conv(P_{k+1}).$$
\end{theorem}

\begin{proof}[Sketch of the proof]
   Proof by contradiction, use Lemma~\ref{pedpath} and Lemma~\ref{oflow}. 
\end{proof}

Next, an example shows that $F_{k}$ feasibility is not sufficient for membership in the pedigree polytope.

\begin{example}\label{exsix}
Consider $X$ given by
$$
\begin{array}{rcl}
{\bf x}_{4} & = & (0, 3/4, 1/4); \\
{\bf x}_{5} & = & (1/2, 0, 0, 0, 0, 1/2); \\
{\bf x}_{6} &  = & (0, 1/4, 1/4, 0, 1/4, 1/4, 0, 0, 0, 0).
\end{array}
$$
It can be verified that  $X \in P_{MI}(6)$, using Lemma~\ref{obvious}, with $n =6$, and $Y = X$. And $F_{4}$ is
feasible, for instance, $f$ given by
$$f_{([4:1,3],[5:1,2])} = 1/4, f_{([4:2,3],[5:1,2])} = 1/4, f_{([4:1,3],[5:3,4])} = 1/2$$ is one such feasible solution.

Or $X/5 \in conv(P_{5})$ and notice that no arc in $F_{4}$ is rigid, that is,  $R_{4} = \emptyset$. Thus, $N_{4}$ is well-defined. 
Next, the restricted networks
$N_{4}(L)$ for the links in $\{ (1, 2) , (3,4) \} \times \{ (1,3),
(2,3)) ,(2,4), (3,4) \}$ are shown in Figure~\ref{fig:r_netE2}, where $L$ is given as the arc connecting a node in $V_{[2]}$ with a node in $V_{[3]}$. 
$F_{5}$ is feasible, with $f$ given by
$$f_{([5:1,2],[6:1,3])} = f_{([5:1,2],[6:3,4])} = 1/4,f_{([5:3,4],[6:2,3])} = f_{([5:3,4],[6:2,4])} = 1/4.$$
\begin{figure}[H]
    \centering
    \includegraphics[scale=.85]{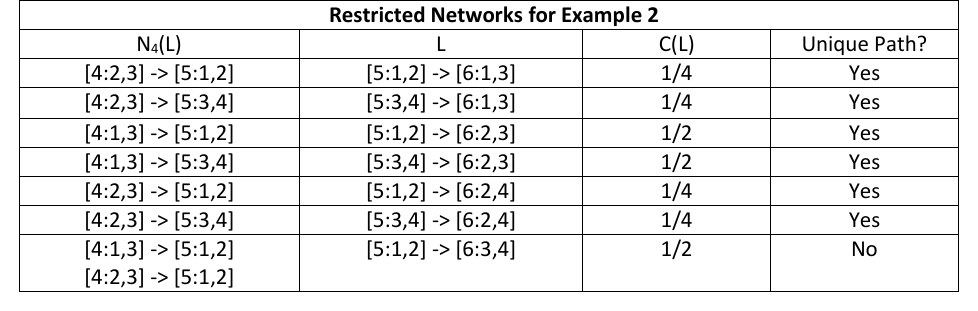}
    \caption{Restricted Networks $N_{4}(L)$ \ for  \ Example~\ref{exsix}.}\label{fig:r_netE2}
\end{figure}
We can have  $N_{5}$ from that. $R_{5}$ is empty. Next, we need to check whether we have evidence for  $N_{5}$ being well-defined. Or  Does $X/6$  in $conv(P_{6}) ?$
Suppose $X/6$ is in $conv(P_{6})$.  Consider any $\lambda\in
\Lambda_{6}(X)$, then  $\lambda$ necessarily assigns on the  pedigrees corresponding to $([4:2,3],[5:1,2]), [6:1,3])$ and  $([4:2,3],[5:3,4]), [6:1,3])$  a total weight of $1/4$ to saturate the capacity of $[6:1,3]$. This implies the capacity of $[4:2,3]$ is saturated as well.  Similarly,there are  pedigrees $X^{s}, s \in I(\lambda
)$ such that $x^{s}_{6}(2,4) = 1$. The total weight for these
pedigrees is  $x_{6}(2,4) =1/4.$ But these pedigrees should have
$x^{s}_{4}(2,3) =1$ as $(2,3)$ is the only generator for $(2,4)$. However, this is not possible as the node capacity of $[4:2,3]$ is already saturated.  Therefore  $X \notin conv(P_{6}).$
\end{example}

 \subsection{Pedigree Packability of  arc flows}\label{ppackability}
Suppose $(N_{k-1}, R_{k-1}, \mu)$ well-defined, so we have evidence for $X/k \in conv(P_{k})$.   Assume that $F_{k}$ is feasible. 
\begin{definition}
Consider any feasible flow, $f$   in $N_{l-1}(L)$ for a link $L
=(e_{\alpha}, e_{\beta}) \in E_{l-1} \times E_{l}$. Let $v_{f}$ be
the value of the flow $f$, that  reaches the sink in
$N_{l-1}(L).$ We say $v_{f}$ is {\em pedigree packable} in case
there exists a subset, $P(L) \subset P_{l}$ such that
\begin{enumerate}[(a)]
\item $\lambda_{r} (\geq 0) $ is the flow along ${\it
path}(X^{r})$ for $X^{r} \in  P(L)$, \item $e^{r}_{l} =
e_{\alpha}, \ X^{r} \in  P(L)$, \item $\sum_{r \ni x^{r}(v) =1}
\lambda_{r} \leq {c}(v), v \in \mathcal{V}(N_{l-1}(L)),$  and \item
$\sum_{X^{r} \in P(L)} \lambda_{r}  = v_{f}.$
\end{enumerate}
where,  ${c}(v)$ is the capacity for node $v$, and ${c}(v) = \bar{x}(v)$ if $v \in  N_{l}$ or ${c}(v) = \bar{\mu}(v)$ if $v \in R_{l}$ for  any $l$.

We refer to $P(L)$ as a {\em pedigree pack} of $v_{f}.$
\end {definition}
\begin{definition}{[Extension Operation]}
Given a  pedigree pack corresponding to a flow $f$ in $N_{l-1}(L)$
for a link $L =(e_{\alpha}, e_{\beta })$ with $v_{f} >0,$ we call
$\overrightarrow{X^{r} L} = (X^{r}, {\bf x}^{r}_{l+1} )$ the {\em
extension}  of $X^{r} \in P(L)$ in case
\begin{equation}\label{ext}
x^{r}_{l+1}(e) = \left\{ \begin{array}{ll}
1 &  \textrm{if} \ e = e_{\beta} \\
0 &  \textrm{otherwise.}
\end{array} \right.
\end{equation}
\end{definition}

That is, $P^{r} = (e^{r}_{4}, \ldots , e^{r}_{l}= e_{\alpha})$ and
this pedigree can be extended to $(e^{r}_{4}, \ldots , e^{r}_{l}=
 e_{\alpha},  e_{\beta} )$. And the corresponding characteristic vector,
$(X^{r}, {\bf x}^{r}_{l+1})  \in P_{l+1}$. Observe that $v_{f}
> 0$ implies a generator of $e_{\beta}$ is in $P^{r}$.

\begin{definition}
Given $l, 4 \leq l < k \leq n, and P \in P_{l}, EXT(P, k)$ denotes the set of pedigrees in $P_{k}$ that are extensions of $P \in P_{l}, 4 \leq l < k \leq n.$ 
That is, 
$$EXT(P, k) = \{ P^{\prime} \in P_k| P = P^{\prime}/l  \}$$.
\end{definition}

We are interested in the pedigree packability of the arc flows in $F_{k}$.
It is easy to see that for any $[k+1:e] \in V_{[k-2]}$ with an arc $(P,e)$ connecting a pedigree $P \in R_{k-1}$, we have a flow
$f_{P, e} \geq 0$ and this flow is pedigree packable, as the $path(P, e)$ is the pedigree path bringing that flow into $e$. 
We wish to check whether this is true for other arcs in $F_k$  as well. Let the total flow from pedigrees in $R_{k-1}$ to $[k+1:e^{\prime}] \in V_{[k-2]}$ be denoted by $\delta_{e^{\prime}}$. 

Explicit use of the fact $X \in P_{MI}(n)$ is made 
to establish the pedigree packability of any node capacity at
layer $k-2$ given that $X/k \in conv(P_{k})$.

\begin{theorem}\label{ppack}
Given $X \in P_{MI}(n)$ and  $(N_{k-1}, R_{k-1}, \mu)$ well-defined, that is, $X/k \in conv(P_{k})$. Consider any $\lambda \in \Lambda_{k}(X)$.  For any $[k:e] \in V_{[k-3]} \textrm{ and } [k+1:e^{\prime}] \in V_{[k-2]}$, we have a flow
$f^{L}$ in $N_{k-1}(L)$  for a link $L
=(e, e^{\prime})$, such that,
\begin{enumerate}[(1)]
\item The value of the flow $f^{L},$ given by $v^{L},$ is pedigree packable for link $L,$ 
\item For any $e^{\prime}\in V_{[k-2],}$ we have
$\sum_{L = (e, e^{\prime}), e \in V_{[k-3]}}v^{L} = x_{k+1}(e^{\prime}) - \delta_{e^{\prime}},$ and 
\item $\sum_{L = (e, e^{\prime}), e \in V_{[k-3]}}\sum_{X^{r} \in \mathcal{P}(L), \ x^{r}(u) =1 }\delta_{r}
\leq {c}(u), u \in \mathcal{V}(N_{k-1}),$
where, $\delta_{r}$ is the flow along the $path(X^{r})$, ${c}(u)$ is the capacity for node $u$.
and $$\mathcal{P}(L) = \{ X^{r}, r \in I(\lambda) \mid x^{r}_{k} (e) = 1, e^{\prime} \in T^{r}\}.$$ 
\end{enumerate}
\end{theorem}

In other words, Theorem~\ref{ppack} assures the existence of
pedigree paths in $N_{k-1}$ bringing in a flow of $v^{L}$
into the sink, $[k: e]$,  and all these paths can be extended to pedigree paths in $N_{k}$ or in $R_{k}$
bringing in a total flow of $x_{k+1}(e^{\prime})$ into $[k+1:e^{\prime}] \in
V_{[k-2]}$. (See Appendix~\ref{Appendixresults} Subsection~\ref{proofppack} for 
proof of this Theorem.)

\begin{remark}\label{pack_rem}
Even though we can apply this theorem for any $e^{\prime}$ with  $x_{k+1}(e^{\prime})> 0$, the simultaneous application of this theorem for more than one
$e^{\prime}$, in general, may not be correct.
 This is so because, for some  paths, the total flow for the different $e^{\prime}$
may violate the node capacity for
a node at some layer, $l$. Example ~\ref{next_ex} illustrates this point. 

Notice that We can obtain Theorem~\ref{packability_cor}  as a corollary to Theorem~\ref{ppack}.
\end{remark}

\begin{corollary}\label{packability_cor2}
Given $X \in P_{MI}(n)$ and $X/k \in conv(P_{k})$, if 
\mbox{$\textbf{x}_{k+1}$ is such that $x_{k+1}(e) = 0 , e \in E_{k-1}$,} that is, 
\begin{equation}\label{eq:condxkplus1}
    \sum_{i =1}^{k-1}x_{k+1}((i,k)) =1,
\end{equation}
 then $F_{k}$ feasible implies $X/k+1 \in
conv(P_{k+1}).$
\end{corollary}
\begin{proof}
Since $F_{k}$ is feasible consider any flow $f =(f_{u,v}), (u , v) \in \mathcal{A}$.  Consider any  origin $u$ such  that either  $u = [k:e' = (i,j)]$ or $u = P$ and $P$ ends in $e'$ for some $1 \leq i < j< k$.  We may have flows from $u$ to $v= [k+1:(i,k)]$ or $v' = [k+1:(j,k)]$  only, as other arcs are forbidden or don't exist under the given condition on $\textbf{x}_{k+1}$. Since $X/k \in conv(P_{k}),$ consider any $\lambda \in \Lambda_{k}(X)$. Pedigrees in $I(\lambda)$ can be partitioned depending on their last component  $e'$. The availability at  $u$ is either  $\bar{x}_{k}(e')$ or $\mu_{P}$ and is equal to $f_{u,v} +f_{u,v'}$ as $F_{k}$ is feasible. Thus, these pedigrees can be extended to pedigrees in $P_{k+1}$ with weights $f_{u,v}, f_{u,v'}$ respectively. We can do this for all $u$.  Hence, we have expressed $X/k+1$ as a convex combination of pedigrees in $P_{k+1},$ as required.
\end{proof}
\begin{remark}
    \item[1] In Corollary~\ref{packability_cor2}, satisfying condition given by equation~\ref{eq:condxkplus1} does not automatically guarantee feasibility of $F_{k}$, as can be seen from Example~\ref{condimpli}.  
    \item[2] Theorem~\ref{packability_cor} and~\ref{packability_cor2}  are useful in concluding $X/k+1 \in conv(P_{k+1})$ with easy to check conditions. In general,  we need to solve a multicommodity flow problem to declare $(N_{k}, R_{k},\mu)$ is well-defined, as shown in Section~\ref{multi}.
\end{remark}

\begin{figure}[!ht]
    \centering
    \includegraphics[scale =0.7]{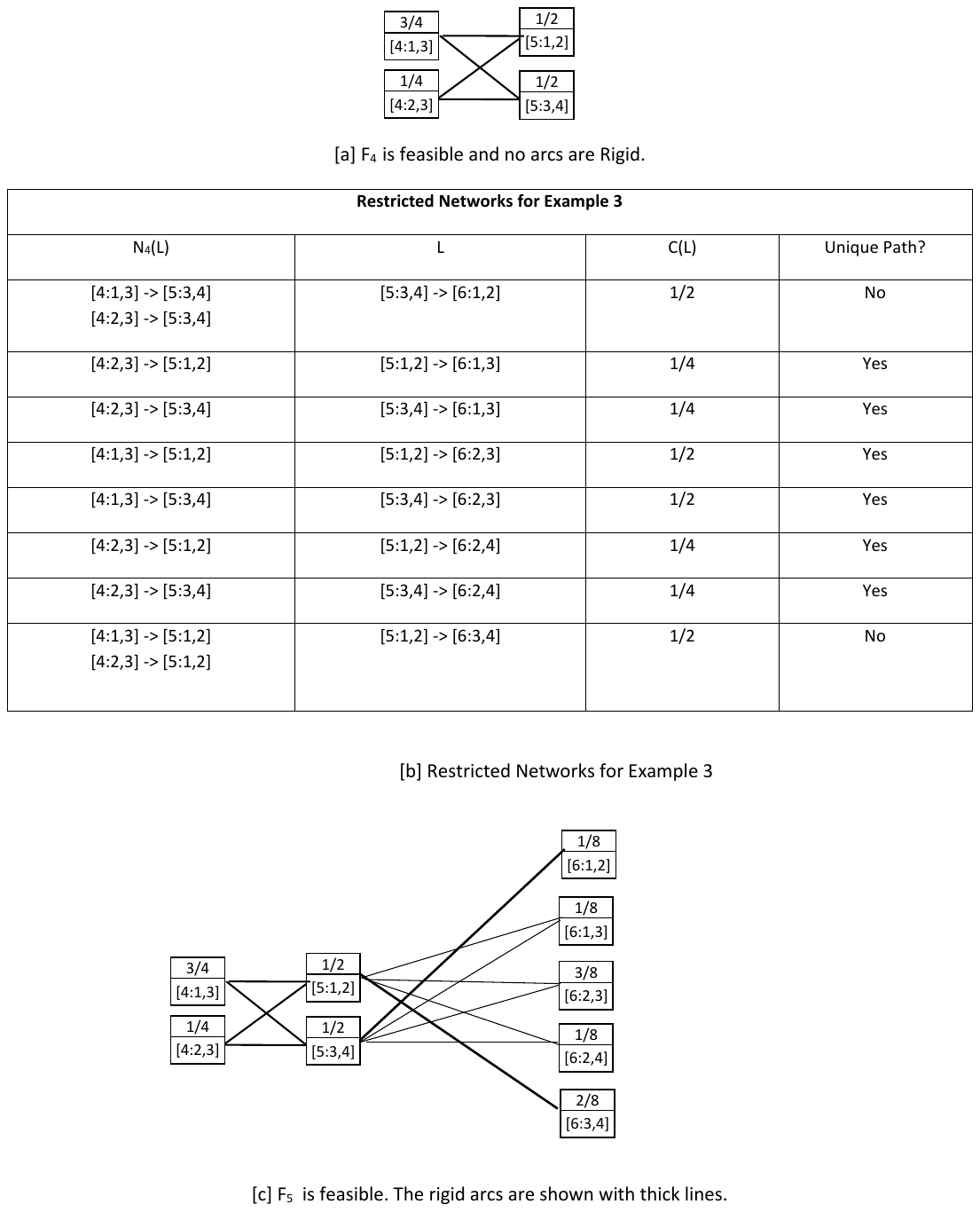}
    \caption{Networks for Example~\ref{next_ex}, [a] $N_{4}$, [b] $N_{4}(L)'s$, and [c]  $N_{5}$.} \label{fig:forEx3}
\end{figure}

\begin{figure}[H] 
    \centering
    \includegraphics[scale =.7]{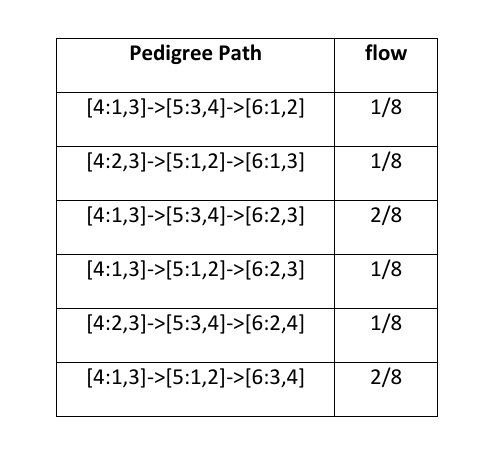}
    \caption{Pedigree Paths showing $X \in conv(P_{6})$ for Example~\ref{next_ex}}\label{fig:pedigreepathsEx3}
\end{figure}

\begin{example}\label{next_ex}
Consider $X$ as given below:
$$
\begin{array}{rcl}
{\bf x}_{4} & = & (0, 3/4, 1/4); \\
{\bf x}_{5} & = & (1/2,0,0,0,0,1/2); \\
{\bf x}_{6} &  = & (1/8,1/8,3/8,0,1/8, 2/8,0,0,0,0).
\end{array}
$$

It can be verified that  $X \in P_{MI}(6)$, using Lemma~\ref{obvious}, with $n =6$, and $Y = X$  and $X/5 \in
conv(P_{5})$ as $F_{4}$ is feasible. See Figure~\ref{fig:forEx3}~[a]. As can be seen from Figure~\ref{fig:forEx3}~[b] flow along each link $L$ is pedigree packable. Also, $F_{5}$ is feasible. Arcs $[5:3,4] -> [6:1,2]$ and $[5:1,2] -> [6:3,4]$ are rigid with frozen flows, $1/8$ and $2/8$ respectively. See Figure~\ref{fig:forEx3}~[c]. However, it is not obvious whether $X \in conv(P_{6}),$ or not. Figure~\ref{fig:pedigreepathsEx3}  gives the pedigree paths with flows, and it can be checked that these flows saturate all the node capacities. Hence, $X \in conv(P_{6}).$ 

But this may not be the case in general. For instance, consider, $X^{*}$ given by
$$
\begin{array}{rcl}
{\bf x^{*}}_{4} & = & {\bf x}_{4}; \\
{\bf x^{*}}_{5} & = & {\bf x}_{5} ; \\
{\bf x^{*}}_{6} &  = & (0,1/4,1/2,0,1/4, 0,0,0,0,0).
\end{array}
$$
It can be checked, $F_{5}$ is feasible but $X/6 \notin conv(P_{6}).$

Instead, consider  $X'$  as given below:
$$
\begin{array}{rcl}
{\bf x}'_{4} & = & (0, 3/4, 1/4); \\
{\bf x}'_{5} & = & (1/2,0,0,1/2,0,0); \\
{\bf x}'_{6} & = & (0,0,0,0,0, 1,0,0,0,0).
\end{array}
$$
We can  check that the paths
$$[4:2, 3] \rightarrow [5:1, 2] \rightarrow [6:3, 4],$$
$$[4:1, 3] \rightarrow [5:1, 2] \rightarrow [6:3, 4],$$
and $$[4:1, 3] \rightarrow [5:1, 4] \rightarrow [6:3, 4]$$ bring
the flows of $1/4$, $1/4$   and $1/2$, respectively to $[6:3, 4]$,
adding to $1 = x_{6}(3, 4)$. And 
so $X' \in conv(P_{6})$, as assured by
Theorem~\ref{ppack}.
\end{example}
\begin{example}\label{condimpli}
Consider  $X$  as given below:
$$
\begin{array}{rcl}
{\bf x}_{4} & = & (0, 1/2, 1/2); \\
{\bf x}_{5} & = & (0,0,0,0,0,1); \\
{\bf x}_{6} & = & (1/2,0,0,0,0, 0,0,0,1/2,0) \\
{\bf x}_{7} & = & (0,0,0,0,0,0,0,0,0,0,1/4,0,1/2,0,1/4).
\end{array}
$$
We can  check that $X$ is in $P_{MI}(7)$, and the paths
$$[4:2, 3] \rightarrow [5:3, 4] \rightarrow [6:1, 2],$$
$$[4:1, 3] \rightarrow [5:3, 4] \rightarrow [6:3, 5],$$
bring a flow of $1/2$ to $[6:1, 2]$,
and $[6:3, 5]$ respectively.  And so $X/6 \in conv(P_{6})$. Though ${\bf x}_{7}$ satisfies the equation~\ref{eq:condxkplus1}, $F_{6}$ is not feasible, as the total demand at $[7:3,6] \& [7:5,6]$ is $3/4$ and only $[6:3,5]$ with availability $1/2$ is connected to these nodes. 
\end{example}
\section{A Multicommodity flow Problem to Check Membership}\label{multi}

Our next task is to check that $(N_{k}, R_{k},\mu)$ is well-defined. To declare the network $(N_{k}, R_{k}, \mu)$ to be well-defined, we need to have evidence that $X/k+1 \in conv(P_{k+1})$. This was easy for $k=4$.
As seen in the Example ~\ref{exsix}, even though there are
pedigree paths bringing the flow along each arc in $F_{k}$, there could be conflicts arising out of the simultaneous capacity restrictions on these flows in the network $N_{k-1}$. This calls for ensuring that these restrictions are not violated. When $N_{k}$ is empty, $R_{k}$ and $ \mu$ give the evidence, as $\sum_{P \in R_{k}}\mu_{P}P  = X/k+1.$   Otherwise,  the multicommodity flow problem defined in this section is required to be solved to check the membership status of $X/k+1$.
\subsection{Defining the Multi-Commodity Flow Problem}

\begin{definition}[Commodities]
Consider the network $(N_{k}, R_{k}, \mu)$. For every arc $a \in \mathcal{A}(N_{k}) \setminus \mathcal{A}(N_{4})$ designate a unique {\em commodity} $s$. Let $\mathcal {S}$ denote the set of all commodities. We write $a \leftrightarrow s$ and read $a$ designates $s$. Let $\mathcal{S}_{l}$ denote the set of all commodities designated by arcs in  $F_{l}$.  
\end{definition}

\begin{definition}
Given a pedigree $X^{r} \in P_{k+1}$, we say that it agrees with
an arc $a = (u, v) \in F_{l}, 4 \leq l \leq  k$ in case either (1) $u = [l:
e^{r}_{l}] ,v = [l+1: e^{r}_{l+1}]$ or (2) $u = P^{r}/l \in R_{l-1}$ with $v = [l+1: e^{r}_{l+1}]$. We denote this by $X^{r}
\parallel a$ and read  $X^{r}$ agrees with $a$. For any
$\lambda\in \Lambda_{k+1}(X)$ and $r \in I(\lambda )$, let
$I_{s}(\lambda )$ denote the subset of $I(\lambda )$ such that the corresponding pedigrees agree with arc $a$  which designates $s$.
\end{definition}

\begin{sidewaysfigure}
    \centering
    \includegraphics[width = \textwidth]{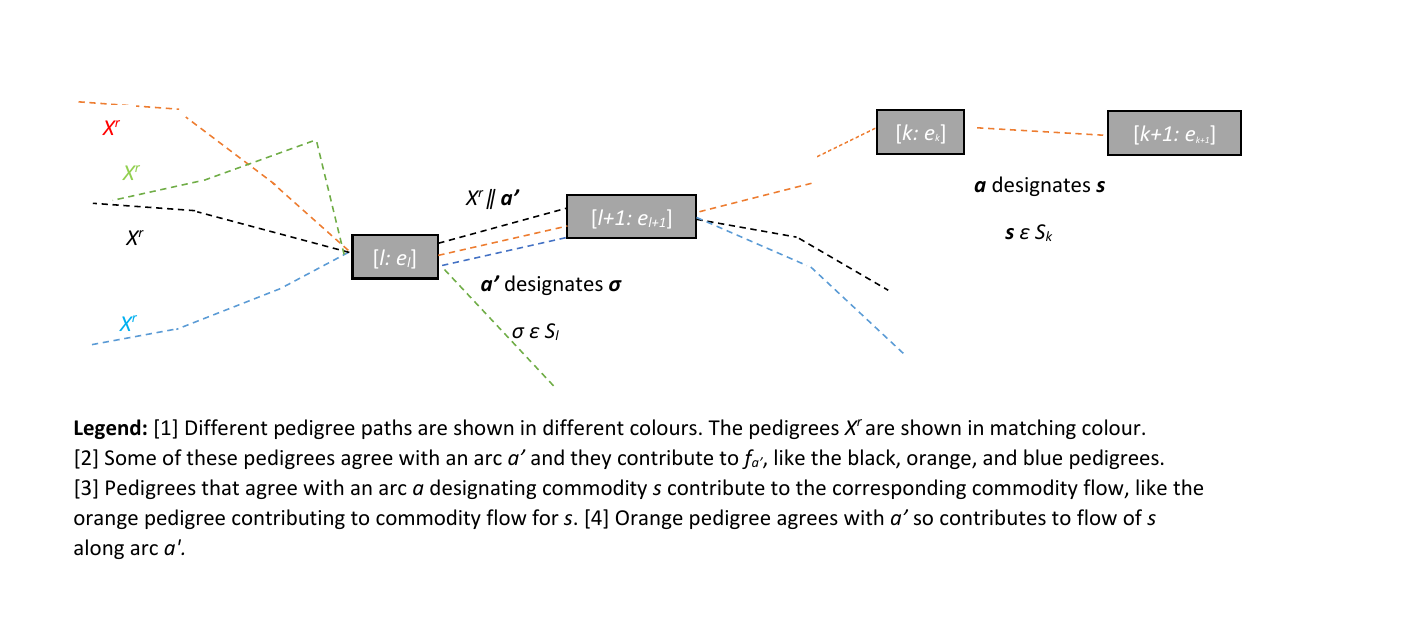}
    \caption{Pedigrees $X^{r}  \parallel a^{\prime}$, arc $a$ designating $s$}
    \label{fig:parallel}
\end{sidewaysfigure}

Figure~\ref{fig:parallel} explains the concepts (1) a pedigree `agreeing' with an arc and (2) an arc `designating' a commodity.

Next, we define the multicommodity flow problem. 

Let $c_a$ be the capacity of any arc $a$ in the network. Let $f_{a}, f^{s}_{a} \geq 0$ be the flow
through arc $a$ and the flow of commodity $s$ through arc $a$, respectively. For any $s$, we allow this
flow to be positive only for the arcs in the restricted network
corresponding to $s$, denoted by $N_{l-1}(s)$, for $s \in \mathcal{S}_{l}$. 
For each $l, 5\leq l \leq k$, we have the following capacity restrictions on the flow through any arc:
 
\begin{equation}\label{ub}
    c_{a} \geq  f_{a}  \geq 0, a \in \mathcal{A}(N_{k})
\end{equation}
\begin{equation}\label{upperbound}
f^{s}_{a} \leq \left\{ \begin{array}{ll}
c_{a} &  \textrm{if} \ a \in  \mathcal{A}(N_{l-1}(s))\\
0 &  \textrm{otherwise.}
\end{array} \right.
\end{equation}

For each commodity, $s$  at each intermediate node $v$ we conserve the flow:
\begin{equation}\label{conserve}
\sum_{u \ \ni \ a = (u, v) }f^{s}_{a} \    =   \ \sum_{w \ \ni \ a
= (v, w) }f^{s}_{a}, \   v \in \mathcal {V}(N_{l-1}) \setminus \{ V_{[1]} \cup V_{[l-2]} \}, s \in \mathcal{S}_{l}.
\end{equation}
We require the total flow of all commodities in $\mathcal{S}_{l}$ through an arc to be the flow $f_{a}$:
\begin{equation}\label{define f}
\sum_{s \in \mathcal {S}_{l}}f^{s}_{a}   =   f_{a},\ a \textrm{ an
arc in } \mathcal{A}(N_{l-1}).
\end{equation}
We define the total flow of commodity $s$ through the arc $a$ designating $s$ to be $v_{s}$ and is equal to the flow along the arc:
\begin{equation}\label{vs def}
\sum_{a^{\prime} \in \delta^{-}(a)} f^{s}_{a^{\prime}}   =   f_{a} \triangleq  v_{s}, a \leftrightarrow s, s \in \mathcal {S}_{l}, 
\end{equation}
where $\delta^{-}(a) \textrm{ is the set of arcs entering } a.$ In addition, at each node $v \in \mathcal{V}(N_{l})$ we have the
node capacity restriction on the total flow through the node as
well: 
\begin{equation}\label{nodecap}
\sum_{s \in \mathcal {S}_{l}}\sum_{u \ni  a = (u, v) }f^{s}_{a}
\leq   \bar{x}(v), \ v \in \mathcal{V}(N_{l}).
\end{equation}

We also have the availability restrictions at the source nodes in $V_{[1]}$ and in $R_{q}$:
\begin{equation}\label{avail1}
    \sum_{w \ni a = (u, w)} f_{a} \leq \bar{x}(u), u \in V_{[1]}.
\end{equation}
\begin{equation}\label{avail2}
    \sum_{w \ni a = (u, w)} f_{a} \leq \bar{\mu}(u), u \in R_{q}, q \in 4 \leq q \leq l-1.
\end{equation}
We also have the node capacity  restrictions at the nodes in $V_{[l-2]}$:
\begin{equation}\label{nodecapk}
    \sum_{u \ni a = (u, v)} f_{a} \leq \bar{x}(v), v \in V_{[l-2]}.
\end{equation}
The objective is to maximise the total flow along the arcs in $F_{k}$:
\begin{equation}\label{objective}
z = \sum_{a \in F_{k}}f_{a} = \sum_{s \in \mathcal{S}_{k}}v^{s}.
\end{equation}

The multi-commodity flow problem can now be stated as:
\begin{problem}{[Problem: $MCF(k)$]}\label{mcflow}
\begin{eqnarray}\label{mulcom}
maximise  & &  z = \sum_{s \in \mathcal{S}_{k}}v^{s}  \\ 
subject \ to \ & &     \\ \nonumber 
for \  all & & l, 5 \leq l \leq k \\ \nonumber
constraints    & & ~(
\ref{ub}) \ \textrm{ through } \ ~(\ref{nodecapk}). \\ \nonumber
\end{eqnarray}

\end{problem}

Next, we observe  some easy-to-prove facts about the multicommodity
network and Problem~\ref{mcflow}.

Let $z^{*}$ denote  the objective function value for  an optimal
solution to Problem~\ref{mcflow}. Let $f=(f^{1}, \ldots ,
f^{|\mathcal {S}|})$ denote any feasible multicommodity flow for
Problem~\ref{mcflow}, where  $f^{s}$ gives the flow vector for
commodity $s \in \mathcal{S}_{l}$, for a fixed ordering of the arcs in the network, $N_{l-1}, 5 \leq l \leq k.$

\begin{remark}\label{mc}
\item[1] $z^{*}$ can be at most $z_{max} = 1 -\sum_{P \in R_{k}} \mu_{P}$, as shown in Remark~\ref{z_max} using Lemma~\ref{extension Lemma}.
\item[2] If $z^{*} = z_{max}$, then for any optimal solution to
Problem~\ref{mcflow}, the bundle capacity $\bar{x}(v)$ at each  node $v
\in V_{[l-3]}, 4 \leq l \leq k+1$ is saturated. \item[3] Every
feasible path bringing a positive flow for any commodity $s \in \mathcal{S}_k$ passes
through each layer of the network, satisfying
the commodity flow restrictions for $s$.  \item[5] If $z^{*} = z_{max}$ then any
optimal solution to Problem~\ref{mcflow}, is such that the
solution restricted to the portion of the network corresponding to
arcs in $F_{l}$, constitutes a feasible solution to $F_{l}$. This
follows from Remark~\ref{mc}.3. Now, any such flow saturates the node
capacities at layers $l-3$ and $l-2$ and $f^{s}_{a}, a \in F_{l}$
can be positive only for the arc $a$ designating $s$. So letting
$\sum_{s \in \mathcal {S}_{l}}f^{s}_{a} = g_{a}$, we can check that $g$  is feasible for
$F_{l}$ as all the restrictions of problem $F_{k}$ are also
present in Problem~\ref{mcflow}. 
\end{remark}
The following lemma helps in determining an upper bound on the optimal value attainable in the multicommodity flow problem.
\begin{lemma}\label{extension Lemma}
Suppose $X/k+1 \in conv(P_{k+1}).$ Consider any  $\lambda \in \Lambda_{k+1}(X)$, and $P \in R_{l-1}$ for some $l \leq k.$
$$\sum_{r \in I(\lambda) \cap EXT(P, k+1)} \lambda_{r}= \mu_{P}.$$
\end{lemma}
\begin{proof}
Without loss of generality, assume $R_{l-1} \neq \emptyset$. 
Proof by induction on $w, l \leq w \leq  k.$

Consider  $w = l$.
Since $\lambda \in \Lambda_{k+1}(X)$ there exists a $\bar{\lambda} \in \Lambda_{l}(X)$ such that, $INST(\bar{\lambda}, l)$ is feasible.  In every $\lambda \in \Lambda_{l}(X)$, any $P \in R_{l-1}$ receives a positive weight equal to $\mu_{P}$, since $R_{l-1}$ is the set of rigid pedigrees. Therefore, 
$$\sum_{r \in I_{\bar{\lambda}} \cap EXT(P, l+1)} \bar{\lambda}_{r} \leq  \mu_{P}$$ for the total flow out of the origin, $P$ cannot exceed the availability at $P$. However, strict inequality is not possible in the above inequality, as that would leave the demand at some destination not met, and so, the $INST(\bar{\lambda}, l)$ problem will not be feasible, leading to a contradiction.  Hence, we have a basis for induction on $w$.
Assume the lemma is true for all $w$ up to $k-1$. 
We shall show that the result is true for $w = k$.

Consider $\lambda \in \Lambda_{k+1}(X)$.
Let $P^{\prime}$ receive a positive weight as per $\lambda$. Let  $P^{r} = P^{\prime}/k$.
Let $\bar{\lambda}_{r} = \sum_{P^{\prime}  \ni  P^{\prime}/k =P^{r}}\lambda_{P^{\prime}}$.
Note that, $\bar{\lambda} \in \Lambda_{k}(X)$. By induction hypothesis, for $\bar{\lambda} \in \Lambda_{k}(X)$ and $P \in R_{l-1}$ we have,
$$\sum_{r \in I_{\bar{\lambda}} \cap EXT(P, k)}\bar{\lambda}_{r} = \mu_{P}.$$
Noticing $P^{r} \in EXT( P, k)$ and $P^{\prime}$ is an extension of some $P^{r}$, we have the result for $w=k$. Hence the result.
\end{proof}
\begin{remark}\label{z_max}
From the above lemma we can observe that the node capacities at the last
layer in $N_{k}$ add to $1 -\sum_{P \in R_{k}} \mu_{P}$, and so  $z^{*}$ is at most $1 -\sum_{P \in R_{k}} \mu_{P}$, call it $z_{max}$. Also, the node capacities at other layers of $N_{k}$ add up to  $z_{max}$ as well.
\end{remark}
\subsection{Necessary and Sufficient Condition for Membership}
Next, we prove the necessity of the multicommodity flow problem having a feasible solution with an optimum of $z_{max}$, for membership of $X/k+1$  in $conv(P_{k+1})$. 
\begin{theorem}\label{imptheorem}
Given $X/k+1 \in conv(P_{k+1})$ then  there exists   a 
solution for the multicommodity flow problem (Problem
\ref{mcflow}),  with  $z^{*} = \sum_{s \in \mathcal {S}_{k}}v^{s} = z_{max}.$
\end{theorem}
{\bf Proof:}
Since $X/k+1$ is in $conv(P_{k+1})$,  choose any $\lambda\in \Lambda_{k+1}(X)$.
Define $f_{a}$ as follows:
\begin{equation}\label{flowdef}
 f_{a} =  \left\{ \begin{array}{ll}
\sum_{r \in I(\lambda ) |  X^{r} \ \parallel \ a}\lambda_{r} &  \textrm{ if } \ \ a \in N_{k} \\
        0     &  \textrm{otherwise.}
\end{array} \right.
\end{equation}
For $ 5 \leq l \leq k$  and $s \in \mathcal{S}_{l}$ define $f^{s}_{a}$ as follows:
\begin{equation}\label{f_sdef}
 f^{s}_{a} =  \left\{ \begin{array}{ll}
\sum_{r \in I_{s}(\lambda ) |  X^{r} \ \parallel \ a}\lambda_{r} &  \textrm{ if } \ \ a \in N_{l-1}(s) \\
        0     &  \textrm{otherwise.}
\end{array} \right.
\end{equation}

We shall show that  $f, f^{s}$ as defined are feasible for Problem~\ref{mcflow} and the objective value is $z_{max}$.
The nonnegative restrictions on $f_{a}$ and $f^{s}_{a}$ are met as $\lambda_{r}$ are positive for all $ r \in I(\lambda)$.  Using arguments similar to those in the proof of Theorem~\ref{nonmemberSuf}, we can observe that the flow $f$ defined by equation~\ref{flowdef} is feasible for $F_{l}$,  all the capacity restrictions on arcs and nodes are met. For each commodity $s \in \mathcal{S}_{l}$,  for $r \in I_{s}(\lambda )$, note that $X^{r}/l$ represents a path in $N_{l-1}$ or the corresponding path is in $R_{l-1}$.  For any of these paths, commodity flow is conserved at each node along the path. A node not in any of these paths does not have a positive flow of this commodity through that node.
\newline
We have, 
\begin{equation}
\sum_{u \ \ni \ a = (u, v) }f^{s}_{a} =  \sum_{u \ \ni \ a =(u,v)}
\  \sum_{r \in I_{s}(\lambda )|X^{r} \ \parallel \ a}\lambda_{r}
\end{equation}
\begin{equation}
\ \ \ \  =  \sum_{r \in I_{s}(\lambda ) |  x^{r}_l(e) =1 }\lambda_{r},\\
\end{equation}
for $v=[l:e] \in \mathcal {V}(N_{l}), s \in \mathcal {S}_{l}.$

Similarly
\begin{equation}
\sum_{w \ \ni \ a = (v, w) }f^{s}_{a}  =  \sum_{w \ \ni \ a = (v,
w)} \ \sum_{r \in I_{s}(\lambda ) |  X^{r} \ \parallel \
a}\lambda_{r}
\end{equation}
\begin{equation}
 \ \ \  \ =  \sum_{r \in I_{s}(\lambda ) |  x^{r}_l(e) =1 }\lambda_{r},\\
\end{equation}
for $v=[l:e] \in \mathcal {V}(N_{l}), s \in \mathcal {S}_{l}.$

Hence, commodity flow conservation restrictions are all met. 
Moreover, for $l = k$, the flow $v^{s}$ along the arc $a$ in $F_{k}$, defining commodity $s$ in $S_{k}$,
are nonnegative and they add up to $\sum_{v \in \Delta^{k+1} }\bar{x}_{k+1}(v)$. And every $\lambda \in \Lambda_{k+1}(X)$ 
is such that, there exists  $r \in I(\lambda)$ corresponding to $P \in R_{k}$ with $\lambda_{r} = \mu_{P}$. Hence, $1 - \sum_{P \in R_{k}}\mu_{P}$ is the maximum possible flow in Problem~\ref{mcflow} and that equals $\sum_{v \in \Delta^{k+1} }\bar{x}_{k+1}(v)$.
Hence the total flow in the network is $\sum_{s \in \mathcal
{S}_{k}}v^s = \sum_{s \in \mathcal {S}_{k}}f^s_a =\sum_{s \in \mathcal
{S}_{k}}\sum_{r \in I_{s}(\lambda )}\lambda_r =z_{max}$. 
Thus, we have verified that $f, f^{s}$ is feasible and the objective function value is $z_{max}$. \qed

An implication of Theorem~\ref{imptheorem} is what we can conclude when the Problem~\ref{mcflow}  has a maximal flow less than $z_{max}$, which is stated as Corollary~\ref{suffnonmem}.
\begin{corollary}\label{suffnonmem}
If Problem~\ref{mcflow} has a maximal flow $z^{*} < z_{max}$ then $X/k+1 \notin conv(P_{k+1})$. 
\end{corollary}

The most important result of this article is the converse of Theorem~\ref{imptheorem}. First, we need two lemmas that are used to show the converse. Assume that we have an optimal solution $(f,f^{s})$ for the multicommodity flow problem with $z* = z_{max}.$ Recall that any $a  \in  \mathcal{A}(F_{l})$ is such that, either  $ a = [u, v], u \in V_{[l-3]}$,  with $v \in V_{[l-2]}$ or  $a =(P,v)$,  with $(P,v)$  not rigid, $P \in R_{l-1}$, $v \in V_{[l-2]}.$ (If $[P, v]$ is rigid then the extended pedigree would be included in $R_{l}$ and the arc $a$ will not be included in $\mathcal{A}(F_{l})$). 

When we have a solution to the multi-commodity flow problem, essentially we are apportioning the node or arc capacity to the different commodities in $\mathcal{S}_k$ that flow through a node $u$. But $u$ may be a node in layer $l-3$ of $N_{k}$, (if it has not been deleted because the updated capacity became zero), and it may occur in some of the rigid paths corresponding to pedigrees in   $R_{l-1}, \ldots, R_{k-1}.$ 

{\bf Caution:} Hereinafter, without loss of generality, we  assume $\mathcal{S}_l, l \in [5, k]$ is updated to contain only $s$ with positive $ v_s.$ 
\begin{definition}\label{Ys_def}
For any  $s \in \mathcal{S}_{k}$,  let  $Y^{s} \in R^{\tau_{k}},$ denote the vector  $(y^{s}_{4}((1,2)), \ldots, y^{s}_{k}((k-2,k-1))),$ 
with $y^{s}_{l}(e), 4 \leq l\leq k, e \in E_{l-1}$ given by
\begin{equation}\label{ys_sdef}
 y^{s}(u)  \equiv y^{s}_{l}(e) =  \left\{ \begin{array}{ll}
\sum_{a =[u,v] \in F_{l}} f^{s}_{a} + \sum_{q = l}^{k}\sum_{a \in \mathcal{A}(u,q)} f^{s}_{a} \ \textrm{ if }  u = [l, e] \in V_{[l-3]}, \\
        0       \textrm{              otherwise.       }
\end{array} \right.
\end{equation}
where $ \mathcal{A}(u,q) = \{ a \mid a =[P, v] \in F_{q}, \  P \in R_{q-1}  \ \&  \  u \textrm{ occurs in } path(P)\}.$
\end{definition}


$ \mathcal{A}(u,q)$ is the set of all arcs $[P, v] \in F_{q}$ such that $P$ is a rigid pedigree in $R_{q-1}$,and the $path(P)$ has $u$ in it. Thus, $\sum_{q = l}^{k}\sum_{a \in \mathcal{A}(u,q)} f^{s}_{a}$ gives the total flow of commodity $s$ along the rigid paths that contain $u$. So,  $y^{s}(u)$ gives the capacity used up for the flow of commodity  $s$ through node $u$, as per solution $f^{s}$.  

Note that,  $\sum_{u \in V_{[l-3]}}y^{s}(u) = \sum_{e \in E_{l-1}}y^{s}_{l}(e) = v^s, 4 \leq l \leq k$. Also, $\sum_{s \in \mathcal{S}_{k}}y^{s}(u) = \bar{x}(u) = \bar{x}_{l}(e), 4 \leq l \leq k, e \in E_{l-1}.$
Therefore $\sum_{s \in \mathcal{S}_{k}}Y^{s} = \bar{X}/k,$
where $\bar{X}/k$ gives the vector of residual capacities at the nodes, after adjusting for the fixed flow $\mu_{P}$  for rigid pedigree $P \in R_{k}$.

\begin{lemma}\label{YsinMI}
Given that we have an optimal solution  for Problem~\ref{mcflow} with $z^{*} = z_{max}$, consider the optimal commodity flow  $f^{s}, s \in \mathcal{S}_{k}.$ Let $Y^{s}$ be as defined above in equation~(\ref{ys_sdef}), then,
$1/v^s (Y^{s}) \in P_{MI}(k)$, for any $s \in \mathcal{S}_{k}.$
\end{lemma}

{\bf Proof:} 
Let $Y = 1/v^s (Y^{s}), s \in S_{k}.$ Feasibility of $Y$ for $P_{MI}(k)$ can be verified by checking that  non-negativity of $Y$ is met by definition. And the equality constraints~\ref{eq:7}  are easily met by $Y$, as $\sum_{e \in E_{l-1}}y^{s}_{l}(e) = v^{s}, \ 4 \leq l \leq k.$ 
For $l=3$ as $y^{s}_{4}(e) \geq 0, e \in E_{3}$, $\sum_{e \in E_{3}} y^{s}_{4}(e) = v^{s}.$ So $y_{4}(e) = 1/v^{s}(y^{s}_{4}(e)) \leq 1, e \in E_{3}.$  Or  $Y4 \in P_{MI}(4).$ 

If the result is false, there is a smallest $q$ such that $Y/q $ is not feasible for the problem $MIR(q)$.   That is,  the corresponding slack variable vector $U$ has a negative entry, say  $U_{i^{o} j^{o}}  < 0$ for some $e =(i^{o}, j^{o}).$ And since $Y/q-1$ is feasible for $MIR_{q-1}$ the slack variable vector, say $V \geq 0.$

Since $Y^{s} \leq X/k,$ the support of $X/k$ also supports $Y$ (however, some $y^{s}_{l}(e)$ could be $0$).
Now consider the construction of the network $N_{q}.$
In the network $(N_{q-1}, R_{q-1}, \mu),$ the maximum possible flow of commodity $s$ into $[q:e]$ can not be more than the flow through the generator nodes of $e =(i^{o}, j^{o}).$ That is, the maximum flow is less than or equal to $\sum_{e^{\prime} \in G(e)} y^{s}_{j^{o}}(e^{\prime}). $ Moreover, out of this,  nodes $[l:(i^{o}, j^{o})], j^{o} +1 \leq l \leq q-1$ receive $\sum_{l = j^{o} +1}^{q-1}y^{s}_{l}(e).$ But no part of this flow of commodity $s$ can flow into $[q: e]$ with $e =(i^{o}, j^{o}).$ Thus the maximum possible flow is less than or equal to $\sum_{e^{\prime} \in G(e)} y^{s}_{j^{o}}(e^{\prime}) - \sum_{l = j^{o} +1}^{q-1}y^{s}_{l}(e),$ which by definition is $ V_{i^{o} j^{o}}$ and is nonnegative. 
On the other hand, $y^{s}_{q}(e)$  by definition is the total flow of commodity $s$ passing through node $[q:e]$ and so it cannot be larger than $V_{i^{o} j^{o}}$ as the flow conservation is met for commodity $s$ at each intermediate node of $N_{k-1}(s).$
So, we  have $U_{i^{o} j^{o}} = V_{i^{o} j^{o}} - y^{s}_{q}((i^{o} j^{o})) \geq 0,$ leading to a contradiction. \qed

Next, we shall show that for $k =5$, we have Lemma~\ref{converse_for_five}, which will form the basis for a proof by induction to show that  $1/v^{s}(Y^{s}/k) \in conv(P_k), \  s \in \mathcal{S}_{k}$.

\begin{lemma}\label{converse_for_five}
Let $k =5$, and $\sigma$ be a commodity in $\mathcal {S}_{5},$ designated by $a =([5:e_{\sigma}^{\prime}],[6: e_{\sigma}]).$ 
If the multi-commodity flow problem (Problem
\ref{mcflow}),  is feasible with  $z^{*} = \sum_{\sigma \in \mathcal {S}_{5}}v^{\sigma} = z_{max}$,  then
\begin{itemize}[label=]
    \item{ 1)} $1/v^{\sigma}(Y^{\sigma}) \in conv(P_{5}), \sigma \in \mathcal{S}_{5},$  and so
    \item{ 2)} $(1/v^{\sigma}(Y^{\sigma}), ind(e_{\sigma})) \in conv(P_{6}), \sigma \in \mathcal{S}_{5}$. Hence, $X/6 \in conv(P_{6}).$ 
 \end{itemize}
\end{lemma}

{\bf Proof:}
When the conditions of the lemma are satisfied, we have $z_{max} = 1 -\sum_{P \in R_{4}}\mu_{P}.$ Without loss of generality assume  $N_{4} \neq \emptyset$, for any $\sigma \in \mathcal {S}_{5},$ all the commodity flow paths in $N_{4}(\sigma)$ are pedigree paths as they correspond to arcs in $N_{4}.$ Hence, we have statement [1] of the lemma. Since any such pedigree path ends in the $tail(a)$, where $a \leftrightarrow \sigma$, and can be extended to a pedigree path corresponding to a pedigree in $P_{6}$ using $head(a).$ And we can do this for each $\sigma \in \mathcal {S}_{5}$. As $ \sum_{ \sigma \in \mathcal{S}_{5}} v^{\sigma}(1/v^{\sigma}(Y^{\sigma}), ind(e_{\sigma})) = X/6,$ we have $X/6 \in P_{6}$ proving statement [2]. 
\qed

\begin{remark}\label{remarkconvlemma}
Observe that for $k >5, 1/v^{s}(Y^{s}/5) \in conv(P_{5}), s \in \mathcal{S}_{k}.$ (similar to the proof of statement [1]). Now $\sum_{s \in \mathcal{S}_{k}} f^{s}_a = v^{\sigma}, a \leftrightarrow \sigma.$ Therefore, for each $s \in \mathcal{S}_{k}$ and $a \in F_{5},  f^{s}_a$ is pedigree packable. So $1/v^{s}(Y^{s}/6) \in conv(P_{6}).$ 
\end{remark}

\begin{lemma}\label{Ysinconv}
Given that we have an optimal solution  for Problem~\ref{mcflow} with $z^{*} = z_{max}$, consider the optimal commodity flow  $f^{s}, s \in \mathcal{S}_{k}.$ Let $Y^{s}$ be as defined  in equation~\ref{ys_sdef}, then,
$1/v^{s}(Y^{s}/k) \in conv(P_k), \  s \in \mathcal{S}_{k}$.
\end{lemma}
{ \bf Sketch of the proof:}
We have $1/v^{s}(Y^{s}/5) \in conv(P_5), s \in \mathcal{S}_{k},$ from Remark~\ref{remarkconvlemma}. Assume the result  for all $l \leq k-1,$ that is, $1/v^{s}(Y^{s}/l) \in conv(P_l), $ for  $s \in \mathcal{S}_{k}.$ 
Show that the result is true for $l = k$. We have pedigree packability of any arc $a^{\prime}$ in $F_{k-1}$ that has a flow $f^{s}_{a^{\prime}}.$ (An application of Theorem~\ref{ppackability}, with  induction hypothesis for $l=k-1$) and  $1/v^{s}(Y^{s}/k) \in P_{MI}(k)$ (from Lemma~\ref{YsinMI}) we have this. However, we need to show that all such arcs in $F_{k-1}$ carrying the flow for commodity $s$ are simultaneously pedigree packable. But commodity $s$ is carried only by arcs $(w,u)$ ending in $u = [k:e^{\prime}]$, for some $e^{\prime}$ and arc $(u, v)$ in $F_{k}$  designates $s$. This implies $y^{s}_k(e^{\prime}) = \sum_{a^{\prime}=(w,u)} f^{s}_{a^{\prime}} = v^{s}$. Therefore, $1/v^{s}(Y^{s}/k)$ has the last component as a unit vector corresponding to $e^{\prime}$. Thus, we have conditions of Theorem~\ref{packability_cor} satisfied. \qed

\begin{theorem}\label{impconvtheorem}
If the multi-commodity flow problem (Problem
\ref{mcflow}),  is feasible with  $z^{*} = \sum_{s \in \mathcal {S}_{k}}v^{s} =z_{max}$, then  $X/k+1 \in conv(P_{k+1})$
\end{theorem}

{\bf Proof:}
Assume that the multi-commodity flow problem (Problem
\ref{mcflow}),  is feasible with  $z^{*} = \sum_{s \in \mathcal{S}_{k}}v^{s} = z_{max}$.
Let $a = ([k:e_{s}^\prime],[k+1:e_{s}])  \leftrightarrow s.$
From Lemma~\ref{Ysinconv} we have $1/v^{s}(Y^{s}/k) \in conv(P_k), \  s \in \mathcal{S}_{k}$.
Consider any such $s$, there exists $\gamma_{P} > 0$, and $\sum_{P \in \mathcal{P}_{s}}\gamma_{P} = 1$, such that $1/v^{s}(Y^{s}/k)  = \sum_{P \in \mathcal{P}_{s}}\gamma_{P}X_{P},$ where $\mathcal{P}_{s}$ is a subset of pedigrees in $P_{k}$, and $X_{P}$ is the characteristic vector of $P$.
Each such $P$ ends in $e_{s}^\prime.$ Therefore, $P$ can be extended to a pedigree in $P_{k+1}$ using $e_{s}.$ Notice that the weight of any such $P$ is $v^{s}\gamma_{P}.$ And these weights add up to $v^{s}.$
We can do this for each $s \in \mathcal{S}_{k}.$ Thus we have shown that $\sum_{s \in \mathcal {S}_{k}}1/v^{s}(Y^{s}, ind(e_{s})) \in conv(P_{k+1}).$ But $\sum_{s \in \mathcal {S}_{k}} v^{s}(1/v^{s}(Y^{s}, ind(e_{s}))) = \bar{X}/k+1.$  Now all $P \in R_{k}$ are pedigrees in $P_{k+1}$. As $\sum_{s \in \mathcal {S}_{k}} v^{s} = z_{max} = 1 - \sum_{P \in R_{k}}\mu_{P},$ and  $\bar{X}/k+1 + \sum_{P \in R_{k}}\mu_{P}X_{P} = X/k+1$, we have $X/k+1 \in conv(P_{k+1})$. \qed

This theorem is a vital converse result, showing that having a feasible solution with an objective function value equal to $z_{max}$ for Problem~\ref{mcflow} is sufficient for membership in the pedigree polytope.

Thus, we can state the main result needed for the Framework,  giving a necessary and sufficient condition for membership in the pedigree polytope as Theorem~\ref{maintheorem}.

\begin{theorem}\label{maintheorem}
Given $n >5$, $X \in P_{MI}(n),$ and $X/n-1 \in conv(P_{n-1}),$ then $X\in conv(P_{n})$ if and only if there exists a solution  for the multi-commodity flow problem (Problem
\ref{mcflow}), $MCF(n-1)$,  with  $z^{*} = \sum_{s \in \mathcal {S}_{k}}v^{s} =z_{max}.$
\end{theorem}

The next section is devoted to estimating the computational complexity of checking this necessary and sufficient condition given by Theorem~\ref{maintheorem}.

\section{Computational complexity of checking the necessary and sufficient condition}\label{compcomplex} In this section we explore how easy it is to check the condition given by Theorem~\ref{maintheorem} for $X/k+1$ to be in $conv(P_{k+1})$. Properties of the pedigree polytope known are useful in carrying out this estimation.
Especially, we require [1] the dimension of the pedigree polytope and [2] the geometry of the adjacency structure of the polytope; and we recall results from~\cite{TSADMpaper} and other related works \citep{TSACompoly, arthanari2013DOpaper}.

\subsection{On the Mutual Adjacency of Pedigrees in the Set of Rigid Pedigrees}\label{AdjacencyR}
In this subsection, given $k \geq 5,$ and  $(N_{k-1},  R_{k-1}, \mu)$, is well-defined, we shall show an important theorem on the mutual adjacency of pedigrees in $R_{k-1}$ in $conv(P_{k})$. Also, we discuss how the the size of the set of rigid pedigrees, $|R_{k}|$ grows as a function of $k$.
In general if either $R_{k-1} = \emptyset$ or a singleton, we have nothing to prove. Assume without loss of generality  $|R_{k-1}| > 1.$
\begin{lemma}\label{forR4}
First we notice that   for $k=5$ when  $(N_{k-1},  R_{k-1}, \mu)$ is well-defined,  $R_{4}$ can not have $P^{[i]} =(e^{i}_{4},e^{i}_{5}), i = 1,2$, pedigrees from $P_{5}$ that are are non-adjacent. Or equivalently, pedigrees in $R_{4}$ are mutually adjacent.
\end{lemma}

\begin{figure}[hbt]
    \centering
    \includegraphics[scale = 0.6]{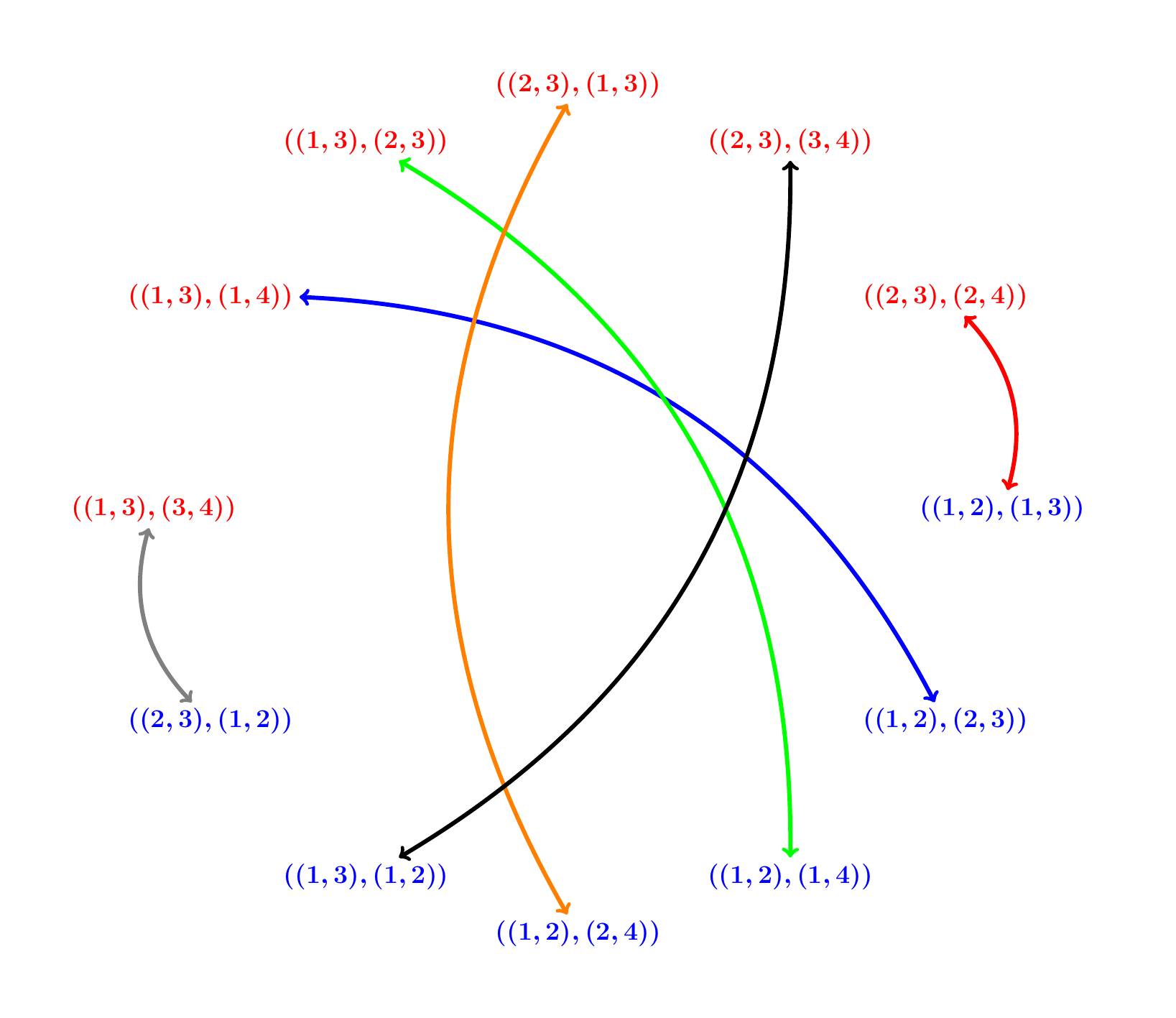}
    \caption{A pair of non-adjacent pedigrees in  $conv(P_{5})$ are shown by a coloured arrows.}\label{fig:P5nadj}
\end{figure}

\begin{remark}\label{onnandc1}
     \item[1] Though $R_{4}$ is constructed differently, can we expect similar result for   $R_{k},k>4$? This is answered affirmatively by Theorem~\ref{adjacency theorem}.
 \item[2] What can we say about the maximum size of $R_{4}$? Or how many of the pedigrees in $P_{5}$ can be mutually adjacent? For this, we need to look at the adjacency structure of the $1-skeleton$ of $conv(P_{5})$. However, it is revealing to present the non-adjacency among pedigrees and then conclude on the mutual adjacency of pedigrees.
In Figure~\ref{fig:P5nadj} we have non-adjacent pedigrees are connected by a coloured arrow. One from a pair of nodes connected by each coloured arrow forms a set of mutually adjacent pedigrees. For instance, Pedigrees printed in red are mutually adjacent. We can see that the maximum size of $R_{4}$ cannot be more than $6$.   In general, we need an upper bound on $|R_{k}|.$ This is done in Theorem~\ref{cardinality theorem}. 
\end{remark}

\begin{theorem} \label{adjacency theorem}
Pedigrees in $R_{k-1}$ are mutually adjacent in $conv(P_{k}).$ That is, the corresponding characteristic vectors of any two rigid pedigrees are adjacent in the $1-skeleton$ or graph of $conv(P_{k})$.
\end{theorem}

Next, we discuss the size of  $R_{k-1}$.
\begin{theorem}\label{cardinality theorem}
Given $k \geq 5$ and  $(N_{k-1}, R_{k-1}, \mu )$, is well-defined, we have
$|R_{k-1}| \leq   dim(\Lambda_{k}(X))+1,$ 
where $\Lambda_{k}(X) = \{ \lambda : \sum_{i \in I(\lambda)} \lambda_{i} P_{i} = X/k,  \sum_{i \in I(\lambda)} \lambda_{i} =1, \lambda_{i} > 0 \}$ (defined in section~\ref{preli}).
\end{theorem}

\begin{proof}[Sketch of the proof]
    The mutual adjacency of the pedigrees in $R_{k-1}$ implies they form a simplex. So there can be at most $ dim(\Lambda_{k}(X)) + 1$ pedigrees in $R_{k-1}$.
\end{proof}

\begin{corollary}\label{CordinalityR}
Given $k \geq 5$ and $(N_{k-1}, R_{k-1}, \mu )$, is well-defined, we have
$|R_{k-1}| \leq  \tau_{k} - k + 4.$
\end{corollary}
\begin{proof}
From the above theorem, $|R_{k-1}| \leq dim(\Lambda_{k}(X))+1.$ But $dim(\Lambda_{k}(X))$ can be at most equal to the dimension of $conv(P_{k})$, which is $\tau_{k} -(k-3).$ (See for proof in \cite{TSADMpaper}.)
\end{proof}
The above results are useful in estimating the computational complexity of the proposed approach to check the membership of a given $X$ in the pedigree polytope.

\subsection{Estimating the Computational Burden at Different Steps of the Framework}\label{estcomp}
In this subsection, we will go into estimating upper bounds on the computational burden at each of the steps of the algorithmic framework given by Figure~\ref{fig:framework}. 
\begin{figure}[H]
    \textbf{FRAMEWORK:} \ \ \textbf{\textit{Membership Checking Steps and Procedures}}
    
    \textbf{Purpose:} \textrm{ Given the number of cities,$n \geq 5$, the framework outlines the stages and conditions to be met for membership of  $X \in Q^{\tau_{n}}$ in pedigree polytope $conv(P_{n})$.} 
    \begin{enumerate}
        \item  \textrm{It requires a certificate that $X$ belongs to $P_{MI}(n)$,} \item \textrm{It requires a certificate for $F_{4}$ feasibility,} \item \textrm{For each $5 \leq  k  \leq n,$ it sequentially checks for  a certificate for $F_{k}$ feasibility, before it checks for a certificate for $X/k \in conv(P_{k})$,}
        \item \textbf{until finally a certificate for $X \in conv(P_{n})$ is produced.} \item \textrm{Otherwise, $X \notin conv(P_{n})$ is established when we fail to produce any of the earlier mentioned certificates.}
    \end{enumerate}
        \textbf{Step:1a} \textrm{If} $X \in P_{MI}(n)$ \textrm{ Proceed further,   otherwise } \textbf{Stop}. \\ 
    \textbf{\bigskip Output:} \textrm{ Evidence exits for} $X \notin conv(P_{n})$,  $X \notin P_{MI}(n).$ \\
    \textbf{Step:1b} \textrm{If} $F_{4},$  \textrm{is feasible, identify }, all $P \in R_{4}$ \textrm{using }  $FFF$ \textrm{ algorithm.  Set } $k =5,$ \textrm{proceed further, otherwise } \textbf{Stop}.\\ 
    \textbf{\bigskip Output:} \textrm{  Evidence exits for$X \notin conv(P_{n})$,}  $F_{4}$ \textrm{ is not feasible.}\\
    \textbf{Step:2a} Find capacity $C(L)$ of each link, $L \in V_{[k-3]} \times V_{[k-2]}$,  that is,  the maximum flow in the restricted network, $N_{k-1}(L).$  Construct $F_{k}.$\\
    \textbf{Step:2b} \textrm{If} $F_{k},$ \textrm{ is feasible, proceed further, } \textrm{ otherwise} \textbf{Stop}. \\
    \textbf{\bigskip Output:} \textrm{ Evidence exits for$X \notin conv(P_{n})$,}  $F_{k}$ \textrm{ is not feasible.}\\
    \textbf{Step:3a} \textrm{Use} $FFF$ procedure and identify rigid and dummy arcs in $F_{k}.$\\
    \textrm{ \bigskip Find rigid paths } $P \in R_{k}$ \textrm{ along with respective rigid flow} $\mu_{P}.$ \\
    \textrm{ Update capacities of nodes and arcs in } $N_{k}.$ \textbf{Construct} $(N_{k},R_{k},\mu).$ \\
    \textbf{Step:4} \textrm{Construct and solve the multi-commodity flow problem, $Problem~\ref{mcflow}$.} \\
    \textbf{Step:5} \textrm{If the optimum total flow, $z^{*}$ is equal to the maximum possible, that is,} $z_{max},$ \textrm{ then, $X/k \in conv(P_{k})$, increment $k$ and proceed further,
    \textrm{ otherwise } \textbf{Stop}. \\ 
    \textbf{\bigskip Output:} \textrm{ Evidence exits for$X \notin conv(P_{n})$,} $X/k \notin conv(P_{k})$.\\
    \textbf{Step:6} \textrm{ If $k < n,$ \textbf{Repeat Step:2a}  on wards, otherwise }\textbf{Stop}.} \\
   \textbf{\bigskip Output:}  $X \in conv(P_{n}).$ \caption{\textbf{Algorithmic Framework for Checking Membership in Pedigree Polytope}}\label{fig:framework}
\end{figure}
Though there may exist better estimates on the worst-case complexity in many of the problems encountered, we do not go into such discussions at this stage.

\textbf{Step:1a} \textrm{ Producing  a certificate for} $X \in P_{MI}(n):$
This involves 
[1] checking $X \geq 0,$
[2] checking $ \sum_{e \in E_{k-1}}x_{k}(e) = 1$  for $ 4 \leq k \leq  n,$
[3] checking sequentially $U^{(l+1)} \geq 0, 3 \leq l \leq n-1,$ where $U$ is as defined in the reformulated Problem~\ref{reform} of $MI- relaxation.$
This can be done in time strongly polynomial in $n$, involving only additions and subtractions.
\textbf{Step:1b} Constructing and checking  $F_{4},$ \textrm{ feasibility.} 
Solving $F_{4}$ is equivalent to solving a $FAT$ problem, with maximum $3$ origins and $6$ destinations with arcs connecting $[4:e]$ with $[5:e^{\prime}],$ if $e$ is a generator of $e^{\prime}$. 
And applying the $FFF$ algorithm for this case is equally easy.

\textbf{Step:2a} Find the capacity of each link, $L$,  that is, the maximum flow in the restricted network, $N_{k-1}(L).$  Construct $F_{k}.$

Firstly, constructing the restricted network $N_{k-1}(L)$ involves applying the deletion rules~\ref{restrict}. Each of these rules is either [i] deleting a certain set of nodes or [ii] checking whether for a node all its generators have been deleted. These can be done in time linear in the number of nodes in the network. And [ii] may be repeated at most $|\mathcal{V}N_{k-1}|$ times, for each node in the network. Finally, finding the sub-network induced by the remaining nodes,
$\mathcal{V}(N_{k-1}(L))$,  can also be done efficiently, to obtain $N_{k-1}(L).$

Secondly, we need to estimate the complexity of finding the capacity for each link, $L.$ Each of these problems can be solved in time polynomial in the number of nodes and arcs in the network  $N_{k-1}(L).$\citep{Networks} 

 This is where the results proved in subsection~\ref{AdjacencyR} become relevant. From Theorem~\ref{cardinality theorem} and corollary~\ref{CordinalityR}, we have a polynomial bound on the cardinality of $R_{k-1}$, namely, $|R_{k-1}| \leq  \tau_{k} - k + 4.$
In $N_{k-1}$   in the last layer we have a maximum of $p_{k-1}$ nodes. In each of the other layers, we have  $ \leq p_{l} + \tau_{l}$ nodes in layer $l >4 $. Therefore, at most $\sum_{l = 5}^{k-1}(p_{l-1} +\tau_{l}) + p_{3} + p_{k-1} = \sum_{l=5}^{k}\tau_{l}$ nodes are possible.  Therefore, the number of nodes is  $\leq (k-5)\times \tau_{k}.$ 
 The number of arcs in $N_{k-1}$ can be at most $\sum_{l = 5}^{k-1}(p_{l-1} +\tau_{l})p_{l} + 3\times 6$.
 The number of links, $L$ we have are at most $ p_{k-1} \times p_{k} \ (< k^{4}),$ as we consider links between nodes in $V_{[k-3]}$ and $V_{[k-2]}.$
Thus, both the number of nodes and arcs in the network are bounded above by a polynomial in $k$ for each link $L$.Or each of these flow problems can be solved in time strongly polynomial in $k$.

\textbf{Step:2b} \textrm{ Checking  $F_{k},$  feasibility.} 
 Though solving a $FAT$ problem can also be done in time polynomial in $|\mathcal{O}|$, $|\mathcal{D}|$ and $|F|$, where $F$ is the set of forbidden arcs, as shown in ~\cite{SandW2018}. Construction of $F_{k}, k > 4$ involves nodes corresponding to pedigrees in $R_{k-1}$ as origins, apart from the nodes in $V_{[k-3]}$. Therefore, it is crucial to ensure the number of origins in $F_{k}$ does not grow to be exponentially large. 
 
 As seen earlier, we have a polynomial bound on the cardinality of $R_{k-1}$, namely,$|R_{k-1}| \leq  \tau_{k} - k + 4.$

 Thus, the number of origins in $F_{k}$ will never grow exponentially.
 
\textbf{Step:3} \textrm{Use} $FFF$ procedure and identify rigid and dummy arcs in $F_{k}.$ 

Frozen Flow Finding ($FFF$) algorithm to identify rigid and dummy arcs, as noted earlier in \ref{preli}, can be done in linear time in the
size of the graph $G_f$ corresponding to the problem $F_{k}$.  The
size of $G_f$ is at most $(p_{k-1}+ |R_{k-1}|) \times p_{k} + p_{k-1} +|R_{k-1}| + p_{k}$.

\textrm{ Find rigid paths for } $P \in R_{k}$ \textrm{ along with respective rigid flow } $\mu_{P}.$
Rigid paths in $R_{k}$ either come from a unique path in $N_{k-1}$ extended by a rigid arc or from a rigid path in $R_{k-1}$  extended using a rigid arc. We store the unique path $P_{unique}(L)$ corresponding to a link $L,$ in case it was obtained while finding the max flow in $N_{k-1}(L)$.

Thus $(N_{k},R_{k},\mu).$  can be constructed in time polynomial in $k$. 

\textbf{Step:4} The next task is to find whether a feasible multi-commodity flow exists with
an objective value equal to $z_{max}.$ 
Since we have only a \emph{linear} and not an integer multi-commodity problem to be solved, we need to solve only a linear programming problem. This linear programming problem is a combinatorial linear programming problem as defined earlier. Theoretically problem can be solved using Tardos' algorithm in strongly polynomial time in the input size of the dimension of the Problem~\ref{mcflow}. And the input size is polynomial in $k$
And we have to solve at most $n-4$ such problems in all.
Thus using the Framework~\ref{fig:framework} we can check in strongly polynomial time  whether $X \in conv(P_{n}).$ 

Here we have not gone for tight bounds for the computational requirements, as the purpose is to prove a crucial result that the necessary and sufficient condition given by Theorem~\ref{maintheorem} can be checked efficiently and is stated as Theorem~\ref{compexity}.
\begin{theorem}\label{compexity}
Given $n, X \in P_{MI}(n),$ and $X/n-1 \in conv(P_{n-1}).$ Then  checking whether there exists a solution  for the multi-commodity flow problem (Problem
\ref{mcflow}),  with  $z^{*} = \sum_{s \in \mathcal {S}_{k}}v^{s} =z_{max}$, can be done efficiently.
\end{theorem}

\section{Concluding remarks}\label{conclusions}
A pedigree defined and studied by the author is a combinatorial object in $1-1$ correspondence with Hamiltonian cycles~\cite{TSADMpaper}. A necessary and sufficient condition for membership in the pedigree polytope is shown as the existence of a feasible multi-commodity flow with maximum possible value in a recursively constructed layered network. We give an algorithmic framework for checking membership, and the computational complexity is strongly polynomial in $n.$  Therefore, the separation problem for the pedigree polytope can be solved efficiently using the construction given in~\cite{Maurras}. \cite{GLS} have shown the connections between separation and linear optimisation over polytopes under certain conditions. After checking that the conditions are satisfied (as done in Chapter 7 of ~\cite{arthanari2023pedigree}), we have the necessary consequence that pedigree optimisation is efficiently solvable. Since the $MI$-formulation of the STSP is a particular instance of the pedigree optimisation problem, the implication of the paper's main result is far-reaching, as this implies we can have a polynomial algorithm to solve the symmetric travelling salesman problem. These results appear in the recent book, \emph{Pedigree Polytopes}~\citep{arthanari2023pedigree}. Primary objective of the article is presenting the efficient solvability of $M3P$ so that experts can vet the results obtained.

 This article opens up new research that applies similar ideas to other combinatorial optimisation problems. Observe that a pedigree is a perfect, pure simplicial complex. And the boundary complex gives us the corresponding Hamiltonian cycle. (See books on algebraic topology for these terms, like~\cite{Munkers}.) We have a new possibility of solving integer optimisation in general and combinatorial optimisation problems, using deeper results from algebraic topology.[See ~\cite{Billera} for supporting view.] Since the implication of efficient solvability of the membership problem for the pedigree polytope is $NP = P$, we have immense potential for future research. A recent article by the author gives the details needed to install the Python Package written by the author, that implements the suggested Framework of solving $M3P$~\cite{arthanari2025py}.



\section{Appendix on Relevant Results and Proofs}\label{Appendixresults}
In this appendix, we collect together results and definitions that are relevant for understanding the proofs used in showing the validity of the Framework that is being implemented.
\subsection{Results on Instant Flow}\label{resultsinstant}
\begin{definition}\label{inducedFAT}
Given $X/k+1 \in conv(P_{k+1}), \lambda \in \Lambda_{k+1}(X)$ we
define a $FAT$ problem obtained from $I(\lambda )$  for a given $l
\in \{4, \ldots , k \}$ called, $INST(\lambda, l)$ as follows:

Partition $I(\lambda )$ in two different ways according to ${\bf
x}_{l}^{r}, {\bf x}_{l+1}^{r}$, resulting in  two partitions
$S_{O}$ and $S_{D}$. We have
$$S_{O}^{q} = \{ r \in I(\lambda ) | x_{l}^{r}(e_{q}) =1 \} , e_{q} \in E_{l-1}  \textrm{ and  } x_{l}(e_{q}) > 0$$
and
$$S_{D}^{s} = \{ r \in I(\lambda ) | x_{l+1}^{r}(e_{s}) =1 \} , e_{s} \in E_{l}  \textrm{ and  } x_{l+1}(e_{s}) > 0.$$
Let $|y|_{+}$ denote the number of positive coordinates of any
vector $y$. Let $n_{O} = |{\bf x}_{l}|_{+}$ and $n_{D} = |{\bf
x}_{l+1}|_{+}.$ Let $a_{q} = \sum_{r \in S_{O}^{q}} \lambda _{r}
= x_{l}(e_{q}), q = 1, \ldots , n_{O}.$ Let $b_{s} = \sum_{r \in
S_{D}^{s}} \lambda _{r} = x_{l+1}(e_{s}), s = 1, \ldots  , n_{D}.$
Let the set of forbidden arcs, $F$,  be given by
$$ F = \{ (q, s) | S_{O}^{q} \cap S_{D}^{s} = \emptyset \}.$$

The problem with an origin for each $q$ with availability
$a_{q}$, a sink for each $s$ with demand $b_{s}$ and the forbidden
arcs given by $F$ is called the $INST(\lambda, l)$ problem.
\end{definition}
\begin{remark}
 \item[1] From Lemma~\ref{theorem:lemma3} we know that such a
problem is feasible and a feasible flow is given by $$ f_{q,s} =
\sum_{r \in S_{O}^{q} \cap S_{D}^{s}} \lambda _{r}.$$ We call such
an $f$ the {\em instant flow} for the $INST(\lambda, l)$ problem. 
\item[2] The sinks and  the demands in the
above problem are same as that of $F_{l}$. \item[3] But in the problem $F_l$ ,  for an origin $[l:e_\alpha]$ we have availability as the updated $x_l(e_\alpha)$, that is, $\bar{x}([l:(e_\alpha)]$, and for $P \in R_{l-1}$ we have $\mu_{P}$ as the availability. 
But there are no capacity restrictions on the arcs of the $INST(\lambda, l)$ problem.
\end{remark}

We shall show that the  instant flow   for the $INST(\lambda, l)$ problem is indeed feasible for the problem
$F_{l}$. 

\begin{definition}
Given $\lambda \in \Lambda_{k+1}$, and $l$, for $r \in I(\lambda )$,  let $L_{l}^{r}$ denote $(e, e')$ such that $x_{l}^{r}(e) = x_{l+1}^{r}(e') =1$. That is, $L_{l}^{r} $ is the $(l-3)^{rd}$ and $(l-2)^{nd}$ elements of the pedigree $P^{r}$ (given by $X^{r}$).
\end{definition}

\begin{lemma}\label{oflow}
Given $\lambda \in \Lambda_{k+1}$,  and $l$, if for  $r \in I(\lambda )$,  either $path(X^{r}/l)$ is available in $N_{l-1}(L_{l}^{r})$, or corresponds to a pedigree in $R_{l-1}$, 
then the  instant flow   $f$ for the $INST(\lambda, l)$ problem is feasible for $F_{l}$.
\end{lemma}

\begin{proof}[Sketch of the proof]
Given the instant flow for $INST(\lambda, l)$,  consider the flow $f_{q,s}$ along an arc $L =  ([l:e_{q}], [l+1:e_{s}])$. Show that this flow can be split between different arcs connecting $[l+1:e_{s}]$ and origins in $F_{l}$, that have $[l:e_{q}]$ in common, without violating any arc capacities.

In this proof, firstly, we partition $I(\lambda)$ into $\mathbb{S}$ and the rest. $\mathbb{S}$ contains $r \ni X^{r}$ agrees with a pedigree $P$ in $R_{l-1}$.   Then while considering shrunk nodes corresponding to $P \in R_{l-1}$, we partition $R_{l-1}$  into $R_{l-1, q}, q \ni P \in R_{l-1}$ and ends in $e_{q}.$ We also partition $\mathbb{S}$ into ${\mathbb{S}}_{P}$ for $P \in R_{l-1}.$ Thus, $$\bigcup_{q}\bigcup_{P \in R_{l-1, q}} {\mathbb{S}}_{P}= \mathbb{S}.$$
This way, every $r$ in $I(\lambda)$ is accounted for.
First, we consider origins in $F_{l}$ not corresponding to $P \in R_{l-1}$, and then, we consider origins in $F_{l}$ corresponding to $P \in R_{l-1}$ and prove the result.  
\end{proof}

We can also show the conditions stated by  Lemma~\ref{oflow} are true for all pedigrees active for $X/k+1 \in conv(P_{k+1})$.

\begin{lemma}\label{path5}
Every   $P^{*} = (e^*_4, \ldots, e^*_{k+1})$  active for $X/k+1$,
is such that, either  [a] $P^{*}/5$ is available  in $R_{4}$ or [b] $path(P^{*}/5)$ is available in $N_{4}(L_{5}^{*})$, where $L_{5}^{*} = (e^{*}_5, e^{*}_6)$.
\end{lemma}
\begin{proof}[Sketch of the proof]
$Path(P^{*}/5)$ is given by $[4:e_{4}^{*}]\rightarrow [5:e_{5}^{*}]$. We have $[5:e_{5}^{*}]$ in
$N_{4}(L_{5}^{*})$ as the lone sink. If the result is not true, we have either [a] $[4:e_{4}^{*}]$ is not a node in $N_{4}(L_{5}^{*}),$ or [b] $[4:e_{4}^{*}]$ exists but $Path(P^{*}/5)$ is not
 in $N_{4}(L_{5}^{*}).$
In each case, we apply Proof by contradiction and establish the result by examining all the deletion rules applicable in constructing the restricted networks.
\end{proof}

Lemma~\ref{path5} forms the basis to prove Lemma~\ref{pedpath} that is crucial in showing that the infeasibility of $F_{k}$ implies
that $X/k+1 \notin conv(P_{k+1})$.

\begin{lemma}{[Existence of Pedigree Paths]}\label{pedpath}
Every   $X^{*}$   active for $X/k+1$,  is such that, for $4 \leq l \leq k$ either [a] 
$P^{*}/l$ is in $R_{l-1}$, or [b] $path(P^{*}/l)$ is available in $N_{l-1}(L_{l}^{*})$, for $5 \leq
l \leq k$, where $L_{l}^{*}$ denotes $(e, e')$ such that
$x_{l}^{*}(e) = x_{l+1}^{*}(e') =1$
\end{lemma}

\begin{proof}[Sketch of the proof]
 If any
$P^{*}$ active for $X/k+1 \in conv(P_{k+1})$ corresponds to $(P, e_\beta)$ for some $P \in R_{k-1}$, then there is nothing to prove.
We prove the result by induction on l.
From Lemma \ref{path5} we have the result for $l = 5.$    
\end{proof}

\begin{remark}
This is a significant result as it shows our construction of the restricted network does not destroy any pedigree path corresponding to a pedigree that is active for $X/k+1$. Also, it aids us in proving the Theorem~\ref{infeasiblity} that gives a sufficient condition for non-membership in the pedigree polytope. 
\end{remark}

\subsection{Proof of Pedigree Packability}\label{proofppack}
\begin{proof}
Note that $X/k+1 \in P_{MI}(k+1)$. Recall the definition of
$U^{(k-3)}$ in the reformulation of $MI-$ relaxation, we have (Equation~\ref{UAk+1X}), that is,
$$U^{(k-3)} - A^{(k+1)}{\bf x}_{k+1} = U^{(k-2)} \geq 0.$$
In fact, $U^{(k-3)}$ is the slack variable vector corresponding to
$X/k.$ So,  $x_{k+1}(e^{\prime}) \leq U^{(k-3)}(e^{\prime}).$ Now $[k+1:e^{\prime}] \in
V_{[k-2]}$ means $x_{k+1}(e^{\prime}) >0 ,$ and so  $U^{(k-3)}(e^{\prime}) >0.$

Since  $X/k \in conv(P_{k})$, consider any weight vector
$\lambda\in \Lambda_{k}(X).$   Every $X^{r}$ corresponding to a pedigree $P^{r} \in  R_{k-1}$  appears in each of these $\lambda$s with the \emph{fixed weight} $\lambda_{r} =\mu_{P^{r}}$. 

We have
\begin{equation}\label{e1}
\sum_{X^{r} \in I(\lambda ), \ x^{r}(u) =1 }\lambda_{r} = c(u), u
\in \mathcal{V}(N_{k-1}).
\end{equation}
From Lemma~\ref{oneone} Lemma 4.2 from ~\cite{Alt} we have, if $X$ is an integer
solution to the $MI$- relaxation problem for $n$ cities, then $U$, the corresponding slack variable vector, is the edge-tour incident
vector of the $n$-tour corresponding to $X$. Applying this for $n= k$ and noticing $X/k \in conv(P_k)$, we find the same $\lambda$
can be used to write $U^{(k-3)}$ (the slack variable vector for $X/k$) as a convex combination of
$\dot{T}^{r}, \ r \in I(\lambda )$, where $T^{r}$ is the $k-tour$
corresponding to $X^{r}$ and $\dot{T}$ denotes the edge-tour incident vector of $T$.

Let
$$ J_{(e^{\prime})} = \{ r | X^{r} \in P_{k} \textrm{ and } e^{\prime} \in T^{r} \}.$$ Thus,
\begin{equation}
U^{(k-3)}(e^{\prime}) = \sum_{r \in I(\lambda ) \cap J_{e^{\prime}}} \lambda
_{r}\dot{T}^{r}. \label{Ueq}
\end{equation}

Now partition $I(\lambda )$ with respect to ${\bf x}_{k}^{r}$ as
follows: Let $I_{e}$ denote the subset of $I(\lambda )$ with
\mbox{ $x^{r}_{k}(e) =1.$}  $$\sum_{r \in
I_{e}}\lambda_{r} = x_{k}(e), \textrm{ for }
e \in E_{k}.$$ We have,
\begin{equation}
I(\lambda )  \cap J_{e^{\prime}} = \bigcup_{e | I_{e}  \neq \emptyset }
I_{e} \cap J_{e^{\prime}}. \label{parti}
\end{equation}
As   $I_{e}$'s are disjoint, we have a partition  of
$I(\lambda ) \cap J_{e^{\prime}}.$ The  $path(X^{r})$ corresponding
to $ r \in I_{e} \cap J_{e^{\prime}}$, is available in $N_{k-1}(L),$ since any  $X^{r} \in P_{k}$ active for $X/k$ is such that the $path(X^{r})$ is available in $N_{k-1}$ or it corresponds to a rigid pedigree.

Now let $$\mathcal{P}(L) = \{ X^{r} | r \in I_{e} 
\cap J_{e^{\prime}} \}.$$ Or equivalently, $$\mathcal{P}(L) = \{ X^{r}, r \in I(\lambda) \mid x^{r}_{k} (e) = 1, e^{\prime} \in T^{r}\}.$$ 

Any $path(X^{r})$ for an $X^{r} \in
\mathcal{P}(L)$ can  be extended to  a path 
ending in $[k+1: e^{\prime}]$, using  the  arc $([k:e],
[k+1:e^{\prime}]).$ Thus we have a subset of pedigrees in $P_{k+1}$,
corresponding to these extended paths. We see from
equations~\ref{Ueq}~and~\ref{parti} that we can do this for each
$e$, and a maximum of $U^{(k-3)}(e^{\prime})$ can flow into
$[k+1:e^{\prime}]$.

Of which, along the rigid pedigrees in $R_{k-1}$, we have a total flow of $\delta_{e^{\prime}}$ into $[k+1:e^{\prime}]$. 
Since $x_{k+1}(e^{\prime}) \leq U^{(k-3)}(e^{\prime}),$ we can choose nonnegative
$\delta_{r} \leq \lambda_{r}$, so that we have exactly a flow of
$x_{k+1}(e^{\prime})-\delta_{e^{\prime}}$ into $[k+1:e^{\prime}]$ along the paths corresponding to non-rigid pedigrees
$X^{r} \in \cup_{e} \mathcal{P}(L) $. 
Now we have part $3$ of the theorem, from
\begin{enumerate}[a)]
\item $\bigcup_{e} \mathcal{P}(L)$ is a subset of
$\{X^{r} | r \in I(\lambda ) \}$, \item $\delta_{r} \leq \lambda_{r}$
and \item the expression for $c(u)$, given by equation ~\ref{e1}.
\end{enumerate}
Letting $v_{L } = \sum_{X^{r} \in \mathcal{P}(L)}\delta_{r}$ we have  the parts $1$ and $2$ of the  result.  
\end{proof}
\subsection{Proofs of the Lemmas (~\ref{YsinMI}, ~\ref{converse_for_five}) needed in showing the Sufficiency}\label{proofnands}

\begin{proof}[Lemma~\ref{YsinMI}]
Let $Y = 1/v^s (Y^{s}), s \in S_{k}.$ Feasibility of $Y$ for $P_{MI}(k)$ can be verified by checking that we have $U \geq 0,$ satisfying the conditions of Lemma~\ref{obvious}, namely, Equations~\ref{EnX} through \ref{Xnonneg} are satisfied with $n = k.$

Non-negativity of $Y$ is met by definition. And the constraints  $E_{[k]}Y =\mathbf{1}_{k-3}$ in Equation~\ref{EnX} are easily met by $Y$, as $\sum_{e \in E_{l-1}}y^{s}_{l}(e) = v^{s}, \ 4 \leq l \leq k.$ 
Recall that  $U^{(3)}
= \left( \begin{array}{c}
     {\bf 1}_{3}  \\
     \bf{0}
  \end{array}
  \right)$, and  $U^{(l)}$ satisfies:

\begin{equation}
   U^{(l)} - A^{(l+1)}y_{l+1} = U^{(l+1)}, \ 3 \leq l \leq k
\end{equation}
    (See Subsection~\ref{MI_Ins} for definitions.) We shall check that $U^{(l+1)} \geq 0 $ for $3 \leq l\leq k.$

For $l=3$ as $y^{s}_{4}(e) \geq 0, e \in E_{3}$, $\sum_{e \in E_{3}} y^{s}_{4}(e) = v^{s}.$ So $y_{4}(e) = 1/v^{s}(y^{s}_{4}(e)) \leq 1, e \in E_{3}$ and so $U^{(4)} \geq 0.$ 

Suppose the result is not true for all $ l > 3$, then there exists a $q \geq 4$  for which, $U^{(q+1)}_{i^{o} j^{o}}  < 0$ for some $e =(i^{o}, j^{o}).$ 

Since $Y^{s} \leq X/k,$ the support of $X/k$ also supports $Y$ (however, some $y^{s}_{l}(e)$ could be $0$).
Now consider the construction of the network $N_{q}.$
In the network $(N_{q-1}, R_{q-1}, \mu),$ the maximum possible flow of commodity $s$ into $[q:e]$ can not be more than the flow through the generator nodes of $e =(i^{o}, j^{o}).$ That is, the maximum flow is less than or equal to $\sum_{e^{\prime} \in G(e)} y^{s}_{j^{o}}(e^{\prime}). $ Moreover, out of this,  nodes $[l:(i^{o}, j^{o})], j^{o} +1 \leq l \leq q-1$ receive $\sum_{l = j^{o} +1}^{q-1}y^{s}_{l}(e).$ But no part of this flow of commodity $s$ can flow into $[q: e]$ with $e =(i^{o}, j^{o}).$ Thus the maximum possible flow is less than or equal to $\sum_{e^{\prime} \in G(e)} y^{s}_{j^{o}}(e^{\prime}) - \sum_{l = j^{o} +1}^{q-1}y^{s}_{l}(e),$ which by definition is $ U^{(q)}_{i^{o} j^{o}}$ and is nonnegative. 
On the other hand, $y^{s}_{q}(e)$  by definition is the total flow of commodity $s$ passing through node $[q:e]$ and so it cannot be larger than $U^{(q)}_{i^{o} j^{o}}$ as the flow conservation is met for commodity $s$ at each intermediate node of $N_{k-1}(s).$
So, we  have $U^{(q+1)}_{i^{o} j^{o}}  \geq 0,$ leading to a contradiction. 
$\Box$ 
\end{proof}
\begin{proof}[Lemma~\ref{converse_for_five}]
When the conditions of the lemma are satisfied, we have $z_{max} = 1 -\sum_{P \in R_{4}}\mu_{P}.$ If $N_{4} = \emptyset,$ there is nothing to prove as all the commodity flows are using rigid pedigree paths.  If $N_{4} \neq \emptyset$, for any $\sigma \in \mathcal {S}_{5},$ all the commodity flow paths in $N_{4}(\sigma)$ are pedigree paths as they correspond to arcs in $N_{4}.$ Hence, we have statement [1] of the lemma. Since any such pedigree path ends in the $tail(a)$, where $a \leftrightarrow \sigma$, and can be extended to a pedigree path corresponding to a pedigree in $P_{6}$ using $head(a).$ And we can do this for each $\sigma \in \mathcal {S}_{5}$. As $ \sum_{ \sigma \in \mathcal{S}_{5}} v^{\sigma}(1/v^{\sigma}(Y^{\sigma}), ind(e_{\sigma})) = X/6,$ we have $X/6 \in P_{6}$ proving statement [2]. 
\end{proof}
\subsubsection{Complexity of Checking Membership}
\begin{proof}[Lemma~\ref{forR4}]
Suppose $P^{[i]} =(e^{i}_{4},e^{i}_{5}), i = 1,2$, are non-adjacent in the pedigree polytope for $k=5$. Since, $conv(P_{k}), k > 3$ are combinatorial polytopes (\cite{TSACompoly}), from the definition of combinatorial polytopes we have two other pedigrees $P^{[3]}, P^{[4]} \in P_{5}$  such that, $ X^{*} = 1/2(X^{[1]} + X^{[2]}) = 1/2(X^{[3]} + X^{[4]}),$ where $X^{[i]}, $ is the characteristic vector of the corresponding pedigree, $P^{[i]}, i=1, \ldots, 4.$
Since the support of $X^{*}$,  supports $X^{[3]},$ as well as  $X^{[4]},$  one of these (without loss of generality say $X^{[3]}$) will have the first component of $X^{[1]},$ and the second component, different from that of $X^{[1]}$ and so necessarily, it is the second component of $X^{[2]}.$ Similarly $X^{[4]}$ has the first component of $X^{[2]},$ and the second component  of $X^{[1]}.$ Thus, we can have a positive flow, $ \epsilon < min(\mu(P^{[1]}), \mu(P^{[2]})$), along the arcs  for  $P^{[3]}$ and  $P^{[4]}$ and decrease the same along the arcs  for $P^{[1]}$ and  $P^{[2]}$ in $F_{4}$, this contradicts the rigidity of the arcs for $P^{[1]}$ and  $P^{[2]}$.
Hence any two $P^{[i]} \in R_{k-1}, i = 1,2$, are adjacent in the pedigree polytope for $k=5$, namely $conv(P_{5})$.
\end{proof}

\begin{proof}[Theorem~\ref{adjacency theorem}]
Suppose  $P^{[i]} \in R_{k-1}, i = 1,2$ are  non-adjacent in $conv(P_{k})$.
Recall the construction of $R_{k-1}$ from subsection~\ref{DefRk}. Any pedigree $P^{[i]}$ that is included in $R_{k-1}$ is such that it is constructed by either [i] a unique pedigree path in $N_{k-2}$ or [ii] a unique rigid pedigree from $R_{k-2}$ extended by $L_{i} = (e^{[i]}_{k-1}, e^{[i]}_{k})$  that corresponds to a rigid arc in $F_{k-1}$. We have a fixed positive weight $\mu ( P^{[i]})$ attached to $P^{[i]}$, making it a rigid pedigree.  

Let $P^{[i]} =(e^{[i]}_{4}, \ldots, e^{[i]}_{ k-1}, e^{[i]}_{k}), i=1,2.$ Therefore, $P'^{[i]} = (e^{[i]}_{4}, \ldots, e^{[i]}_{ k-1})$ is the unique pedigree path ending in $e^{[i]}_{k-1}$ for the link, $L_{i} = (e^{[i]}_{k-1}, e^{[i]}_{k})$ that is rigid.

We prove the theorem by considering different cases.
\begin{caseof}
\case $P'^{[1]} = P'^{[2]}$

In this case, $L_{i}$ must be different for different $i$, otherwise the pedigrees $P^{[i]}, i = 1,2, $ will not be different.  As the two pedigrees differ in only in one component (the last component), applying Lemma~5.3 in \cite{TSADMpaper}, we assert that $X^{[i]}, i = 1,2,$ are adjacent in $conv(P_{k})$. Hence the result. 
\case $P'^{[1]} and P'^{[2]}$ are different

Since, $conv(P_{k})$ is a combinatorial polytope ( \cite{TSACompoly}), if $X^{[1]}$ and $X^{[2]}$ are non-adjacent, there exist $X^{[3]}$ and $X^{[4]}$ such that $\frac{1}{2}X^{[1]}+X^{[2]} = \frac{1}{2}X^{[3]}+X^{[4]}$. And the support for $X^{[3]}$ and $X^{[4]}$ is same as that of $X^{[1]}$ and $X^{[2]}$. Therefore, a flow $\epsilon$ can be rerouted through the paths corresponding to $X^{[3]}$ and $X^{[4]}$, where $ 0 < \epsilon \leq$ minimum of $\{ \mu(X^{[i]}), i =1,2 \} $. 

\emph{Sub-case a:}  $e^{[1]}_{k-1} = e^{[2]}_{k-1}$

In this case, we should have $e^{[1]}_{k}$ not equal to $e^{[2]}_{k}$. And $P^{[3]}$ or $P^{[4]}$ should end in $e^{[1]}_{k}$ or $e^{[2]}_{k}$. Without loss of generality, let $P^{[3]}$ end in $e^{[1]}_{k}$.
$P^{[3]}$ is different from $P^{[1]}$ but corresponds to the same rigid link, $(e^{[1]}_{k-1}, e^{[1]}_{k})$, this contradicts the uniqueness of the path for this link. Hence the result.

\emph{Sub-case b:} $e^{[1]}_{k} = e^{[2]}_{k}$

[The proof is similar to Sub-case a.]
In this case, we should have $e^{[1]}_{k-1}$ not equal to $e^{[2]}_{k-1}$. And $P^{[3]}$ or $P^{[4]}$ should have $e^{[1]}_{k-1}$ or $e^{[2]}_{k-1}$. Without loss of generality, let $P^{[3]}$ have $e^{[1]}_{k-1}$.
$P^{[3]}$ is different from $P^{[1]}$ but corresponds to the same rigid link, $(e^{[1]}_{k-1}, e^{[1]}_{k})$, this contradicts the uniqueness of the path for this link. Hence the result.

\emph{Sub-case c:} $e^{[1]}_{k-1}$ not equal to $e^{[2]}_{k-1}$ \ \& \ $e^{[1]}_{k}$ not equal to $e^{[2]}_{k}$ 

If either $P^{[3]}$ or $P^{[4]}$ ends with $L_{1}$ or $L_{2}$, then this contradicts the uniqueness of the respective path for these links. Hence the result.

If one of $P^{[3]}$ or $P^{[4]}$ ends with ($e^{[1]}_{k-1}$, $e^{[2]}_{k})$, then  the other should end with ($e^{[2]}_{k-1}$, $e^{[1]}_{k})$. Then consider rerouting a flow of $\epsilon > 0$ along the paths $P^{[i]},1= 1, \ldots , 4$, by increasing the flow along $P^{[3]}$ and $P^{[4]}$ and reducing the same along the paths corresponding to $P^{[1]}$ and $P^{[2]}$. Thus, the links $L_{i}, i=1,2$ are not rigid. This contradicts the fact that links are rigid, as $P^{[i]} \in R_{k}, i=1,2$. Hence the result.
\end{caseof}
This completes the proof of the theorem.
\end{proof} 
\section{Appendix on Frozen Flow Finding Algorithm}\label{Appendixfrozen}
\begin{definition}{[Diconnected Components]}
Given a digraph $G = (V, A)$, we say $u, v \in V$ are {\em
diconnected} if there is a directed path from $u$ to $v$ and also
there is a directed path from  $v$  to $u$ in $G.$ We write $u
\leftrightarrow v$. $\leftrightarrow$ is an equivalence relation
and it partitions $V$ into equivalence classes, called diconnected
components of $G$.
\end{definition}
Diconnected components are also called {\em strongly connected}
components.

\begin{definition}{[Interface]}
Given a digraph $G = (V, A)$, consider the diconnected components
of $G$. An arc $e = (u,v) $ of $G$ is called an {\em interface} if
there exist two different diconnected components $C_{1}, C_{2}$
such that $u \in C_{1} $ and $v \in C_{2}.$

The set of all interfaces of $G$ is denoted by $I(G).$
\end{definition}
\begin{definition}{[Bridge]}
Given a graph $G =(V, E)$, an edge of $G$ is called a {\em bridge} if $G-e$ has more components than $G$, where by component of a graph we mean a maximal connected subgraph of the graph.
\end{definition}

\begin{definition} {[Mixed Graph]}
A {\em mixed graph} $G = (V, E \cup A)$ is such that it has both
directed and undirected edges. $A$ gives the set of directed edges
(arcs). $E$ gives the set of undirected edges(edges). If $A$ is
empty $G$ is a graph, and if $E$ is empty, $G$ is a digraph.
\end{definition}
In general, we call the elements of $E \cup A$, edges. Finding
the disconnected components of a digraph can be achieved using a
depth-first search method in $O(|G|)$ where $|G|$ is the size of G
given by $|V| + |A|.$ Similarly, bridges in a graph can be found in
$O(|G|)$.

\begin{definition}{[$G_{f}$]}\label{Gf}
With respect to a feasible flow $f$ for a $FAT$ problem, we define
a mixed graph $G_{f} = (\mathcal{V}, A \ \cup  \ E)$ where
$\mathcal{V}$ is as given in the $FAT$    problem, and
$$A = \{ (O_{\alpha} , D_{\beta} ) | f_{\alpha , \beta} = 0, (\alpha , \beta ) \in \mathcal{A} \} \  \cup \ \{ (D_{\beta} , O_{\alpha})  | f_{\alpha , \beta} = c_{\alpha , \beta}, (\alpha , \beta ) \in \mathcal{A} \}$$,
$$E = \{ (O_{\alpha} , D_{\beta} )  | 0 < f_{\alpha , \beta}  < c_{\alpha , \beta} , (\alpha , \beta ) \in \mathcal{A}\}.$$
\end{definition}

\begin{definition}{[Flow change Cycle]}
A simple cycle in $G_{f}$ is called a {\em flow change cycle }
(fc-cycle), if it is possible to trace it without violating the
direction of any of the arcs in the cycle. Undirected edges of
$G_{f}$ can be oriented in one direction in one fc-cycle and 
the other direction in another fc-cycle.
\end{definition}

\begin{theorem}{[Characterisation of Rigid Arcs]}
Given a feasible flow, $f$, to a $FAT$ problem, an arc  is rigid
if and only if its corresponding edge  is not contained in any
fc-cycle in $G_{f}$.
\end{theorem}
The proof of this is straightforward from the definitions (see
\cite{Gus} ). It is also proved there that   the set $\mathcal{R}$
is given by the algorithm that we call  {\bf Frozen Flow Finding
($FFF$)}algorithm stated below:

\newcommand{\keyw}[1]{{\bf #1}}
\begin{algorithm}[Frozen Flow Finding]
\begin{tabbing}
\\
\quad \=\quad \=\quad \kill \keyw{Given:} A Forbidden arcs
transportation problem with a feasible flow $f$.\\
\quad \=\quad \=\quad \kill \keyw{Find:} The set of rigid arcs $\mathcal{R}$, in the bipartite graph of the problem. \\
\> \> \keyw{Construct} The mixed graph $G_f$ as per
Definition~\ref{Gf}.\\
 \> \> \keyw{Find} The diconnected components of $G_{f}$ (say
$C_{1}, \ldots , C_{q}$).\\
\> \> \keyw{Find} The set of interfaces, $I(G_{f})$. \\
\> \> \keyw{Find} The set of all bridges $B(G_{f})$ in the
underlying graphs, \\
\> \> treating each $C_{r}$ as an undirected graph.\\
\> \> \keyw{Output} $\mathcal{R} = I(G_{f}) \cup B(G_{f})$. Stop.
\end{tabbing}
\end{algorithm}

We have a linear-time algorithm, as each of the steps 1 - 3
can be done in $O(|G_{f}|)$.
The capacities of the nodes in $N_{k}$ are updated as above.
\newpage
\section{Bibliography}
Bibliographical references includes all cited references and a few related ones. Uncited ones have bibliographic notes as well.
{}
\end{document}